\documentclass[trsc,nonblindrev,fleqn]{informs3} 

\OneAndAHalfSpacedXI 

\graphicspath{{img/}}
\usepackage{algorithm}
\usepackage[noend]{algpseudocode}
\usepackage{multirow}
\usepackage{graphicx}
\usepackage{float}
\usepackage{lscape}
\usepackage{longtable}
\usepackage{hhline}
\usepackage{subcaption}
\usepackage{enumitem}
\usepackage{mathtools}
\DeclarePairedDelimiter{\ceil}{\lceil}{\rceil}

\usepackage{natbib}
 \bibpunct[, ]{(}{)}{,}{a}{}{,}%
 \def\bibfont{\small}%
\TheoremsNumberedThrough     

\EquationsNumberedThrough    

\MANUSCRIPTNO{123456} 

\begin{document}


\RUNAUTHOR{Hasan, Van Hentenryck, and Legrain} 

\RUNTITLE{The Commute Trip Sharing Problem}

\TITLE{The Commute Trip Sharing Problem}

\ARTICLEAUTHORS{%
\AUTHOR{Mohd. Hafiz Hasan}
\AFF{University of Michigan, Ann Arbor, Michigan 48109, \EMAIL{hasanm@umich.edu} \URL{}}
\AUTHOR{Pascal Van Hentenryck}
\AFF{Georgia Institute of Technology, Atlanta, Georgia 30332, \EMAIL{pascal.vanhentenryck@isye.gatech.edu} \URL{}}
\AUTHOR{Antoine Legrain}
\AFF{Polytechnique Montr\'{e}al, Montr\'{e}al, Qu\'{e}bec H3T 1J4, \EMAIL{antoine.legrain@polymtl.ca} \URL{}}
} 

\ABSTRACT{%
  Parking pressure has been steadily increasing in cities as well as
  in university and corporate campuses. To relieve this pressure, this
  paper studies a car-pooling platform that would match riders and
  drivers, while guaranteeing a ride back and exploiting spatial and temporal
  locality. In particular, the paper formalizes the Commute Trip
  Sharing Problem (CTSP) to find a routing plan that maximizes ride
  sharing for a set of commute trips. The CTSP is a generalization of
  the vehicle routing problem with routes that satisfy time-window,
  capacity, pairing, precedence, ride-duration, and driver
  constraints. The paper introduces two exact algorithms for the CTSP:
  A Route-Enumeration Algorithm and a Branch-and-Price
  Algorithm. Experimental results show that, on a high-fidelity,
  real-world dataset of commute trips from a mid-size city, both
  algorithms optimally solve small and medium-sized problems and
  produce high-quality solutions for larger problem instances. The
  results show that car pooling, if widely adopted, has the potential
  to reduce vehicle usage by up to 57\% and decrease vehicle miles
  traveled by up to 46\% while only incurring a 22\% increase in
  average ride time per commuter for the trips considered.  }%


\KEYWORDS{ride sharing; vehicle routing with time windows; column generation; branch and price; mixed-integer programming}

\maketitle

%


\section{Introduction}

Parking occupies a significant portion of our cities. In the United
States, for instance, there are at least 800 million parking spaces
and, in Los Angeles County, 14\% of the city space is devoted to
parking \citep{ParkingLACounty}. Parking also contributes to
congestion: Based on a sample of 22 studies in the United States, the
average share of traffic cruising to find a parking spot is 30\% and
the average cruising time is just under 8 minutes in downtown areas
\citep{Shoup2005,Shoup2006}.

\begin{figure}[!t]
\includegraphics[width=0.32\linewidth]{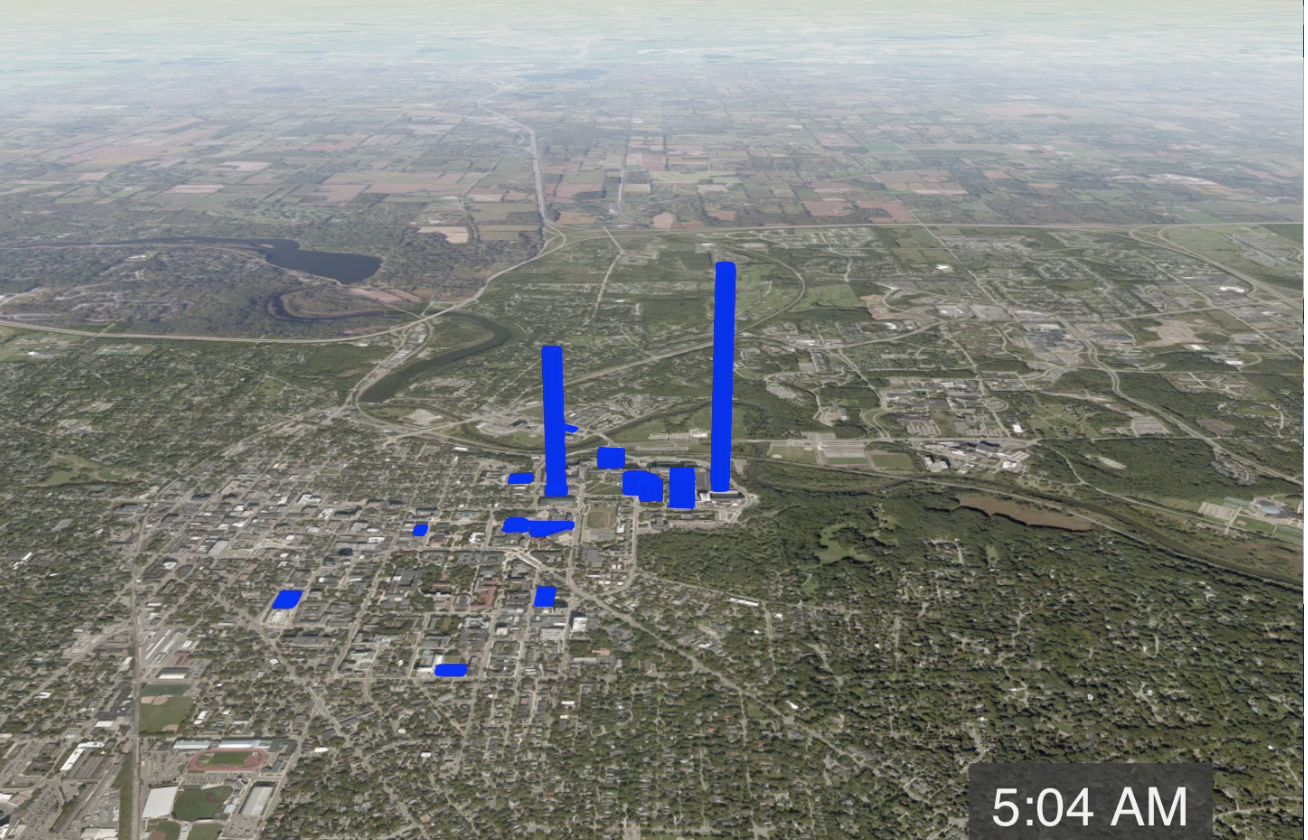}
\includegraphics[width=0.32\linewidth]{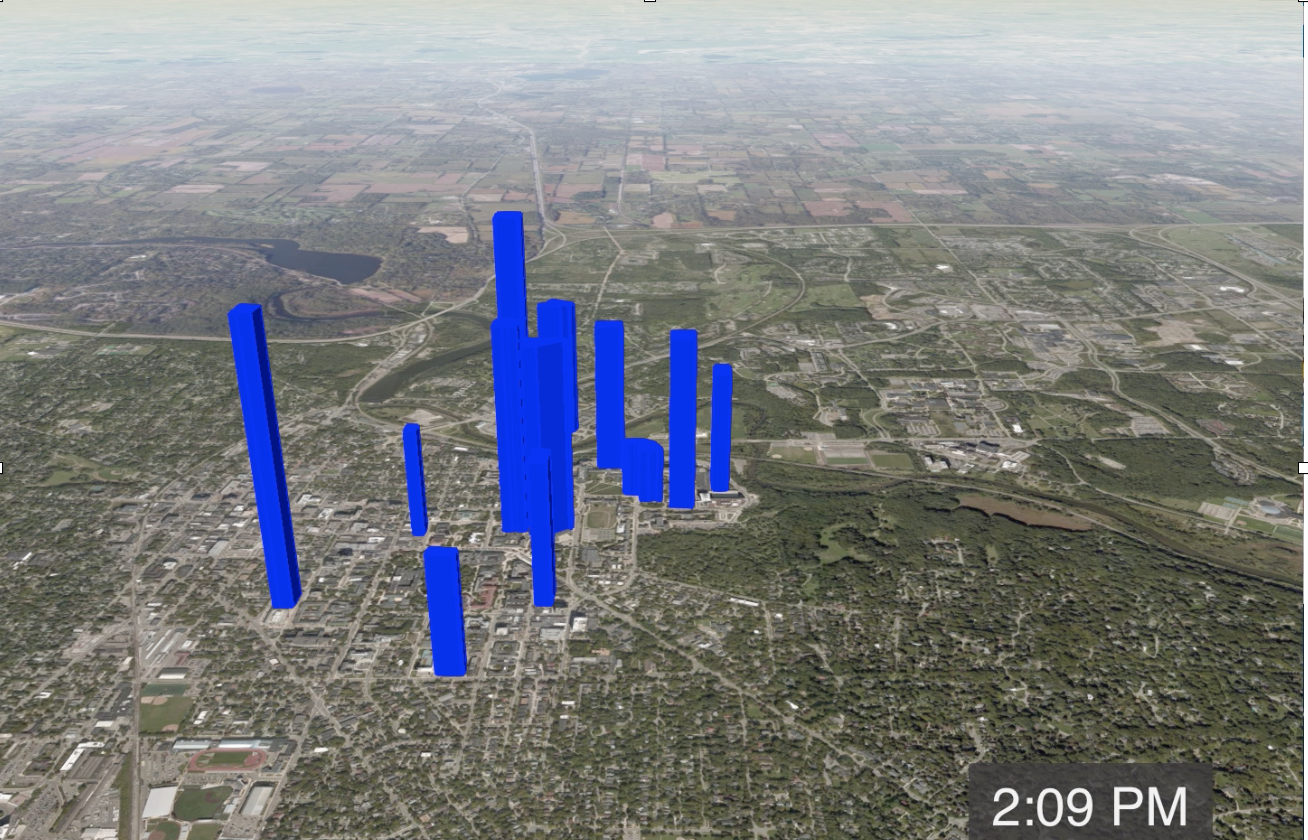}
\includegraphics[width=0.32\linewidth]{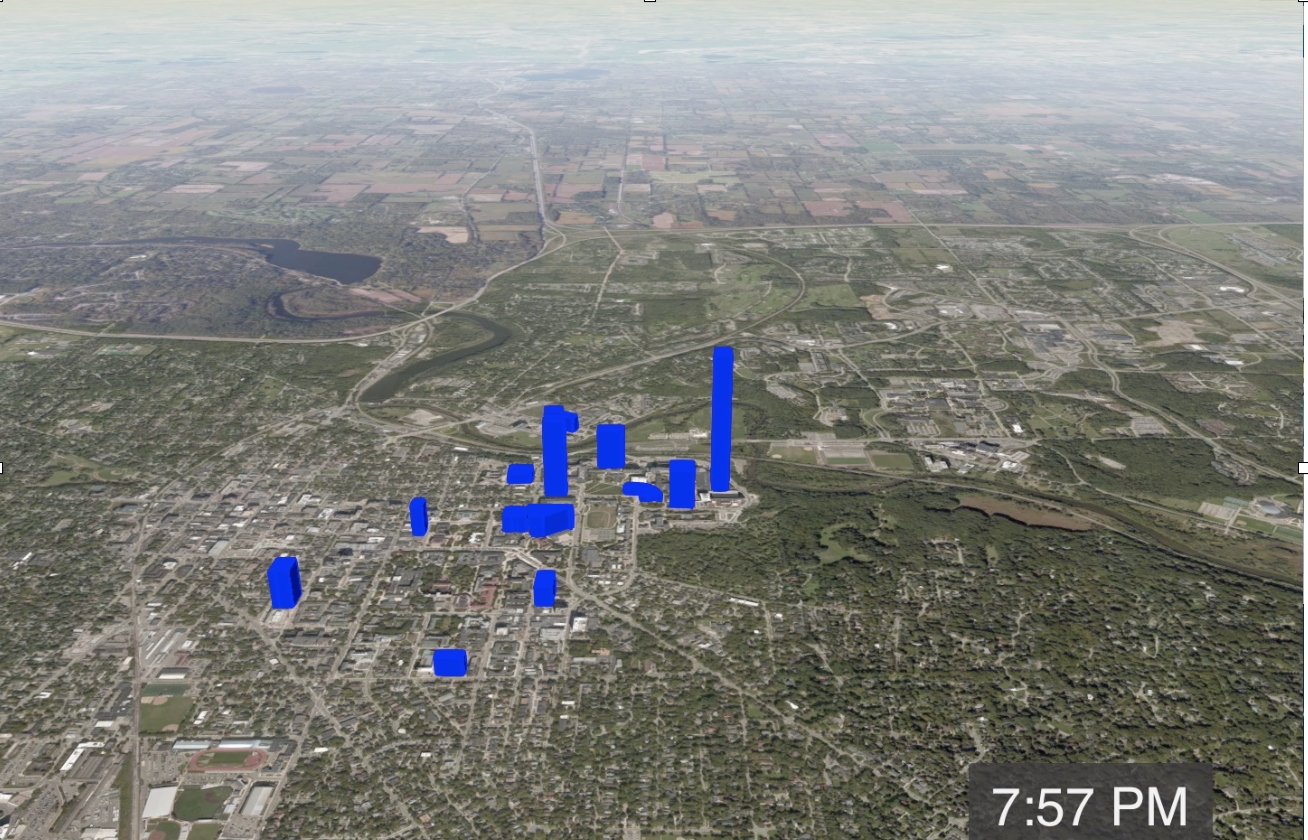}
\caption{The Main University Parking Lots in Downtown Ann Arbor.}
\label{fig:ParkingAnnArbor}
\end{figure}

Parking pressure has also been steadily increasing in cities,
university campuses, and corporations, alongside other concerns 
such as traffic congestion, fuel prices, and greenhouse gas 
emissions. In the city of Buffalo, New
York, the overall supply of parking space has remained constant
for the last 20 years, while the downtown population and the workforce
have increased by 70\% and 30\% respectively
\citep{ParkingBuffalo}. These parking shortages are perceived as an
impediment to future economic developments, as corporations may elect
to move elsewhere when growing their operations. University campuses
feel similar parking pressures. For instance, Stanford University
suffers from a lack of parking spaces due to construction and a growth
in population \citep{ParkingStanford}. The research underlying this
paper was originally motivated by parking pressure at the University
of Michigan in Ann Arbor. Figure \ref{fig:ParkingAnnArbor} depicts the
parking utilization of the 15 most used parking lots in downtown Ann
Arbor. They show a typical parking usage: Cars arrive in the morning,
park in the lot for 6 to 10 hours, and leave the lot in the evening.

To address the increasing demand on these lots, we started to
investigate the potential of a community-based car-pooling program
\citep{hasan2018}. The idea was to implement a car-pooling program
organized around the communities commuting to the university,
exploiting the knowledge of when employees were arriving in the
morning and leaving in the evening. However, while car-pooling has
long been proposed as a solution to reduce peak-hour congestion and
parking utilization, its adoption in the US remains poor as 76.4\% of
American commuters chose to drive alone according to the 2013 American
Community Survey \citep{mckenzie2015}. A study on factors influencing
carpool formation by \cite{li2007} revealed difficulty in finding
people with the same location and schedule as the primary reason for
not carpooling. As a result, we investigated how to alleviate this
burden and studied the feasibility of a matching platform that would
automatically identify commuting groups based on factors determined to
be consequential to individuals’ commuting decisions. One of the
results of our study was the recognition that an effective car-pooling
platform will need to accommodate different sharing patterns for every 
weekday and, as a result, the platform will need to optimize trip matching
on a daily basis to allow significant car pooling to occur \citep{hasan2018}.

The goal of this paper is to propose, and analyze, scalable
optimization algorithms for powering such a platform. A meta-analysis
of related work reveals that car-pooling and car-sharing platforms
should at least implement the following three guiding principles:
\begin{enumerate}
\item Spatial proximity of riders \citep{Richardson1981,Buliung2009};
\item Temporal proximity of riders \citep{Tsao1999,Buliung2010,Poulenez1994};
\item Guaranteed ride back home \citep{correia2011}.
\end{enumerate}
The first two guidelines are natural since car-pooling is unlikely to
occur for riders who are not close spatially or whose schedules are
not compatible. The third guideline is critical: It is unlikely that
many riders will use a platform that does not guarantee a ride back
home in the evening. The guarantee of a ride back home is one of the
main contributions of this work: For instance, the car-pooling
platform {\sc Scoop} provides only weak guarantees for ``ride back''
and with monthly limits on how much auxiliary services can be used
when a ride back is not available. In contrast, this paper approaches
the matching of riders in two steps. In the first step, riders are
grouped spatially into neighborhoods using a clustering algorithm. In
the second step, an optimization algorithm selects drivers and matches
riders to minimize the number of cars and the total travel distance.
The approach follows the three guiding principles listed above and, in
particular, ensures that every rider has a guaranteed ride back.  The
contributions of this paper are threefold.
\begin{enumerate}
\item It first defines the Commute Trip Sharing Problem (CTSP) that
  formally captures the matching problem described previously.

\item It proposes two algorithms, a Route-Enumeration Algorithm and a
  Branch-and-Price Algorithm, to solve the CTSP for a cluster of
  riders.

\item It analyzes the scalability of the two algorithms along
  different dimensions, including the capacity of the vehicles, the
  size of the clusters, and their ability to be deployed in real situations. 
\end{enumerate} 

\noindent
The CTSP can be viewed as a generalization of the Vehicle Routing
Problem (VRP) with routes satisfying time-window, capacity, pairing,
precedence, ride-duration, and driver constraints. In addition to
picking up and dropping off riders within desired time windows while
ensuring vehicle capacities are not exceeded, routes in the CTSP must
also ensure their ride durations are not excessively long to limit
user inconvenience. In this sense, the CTSP shares some similarities
with the Dial-a-Ride Problem (DARP). It differs from the DARP in that
it relies on the use of personal vehicles to serve all trip requests,
which come in pairs for each commuter as each rider makes a trip to
the workplace and another back home. The drivers of these vehicles
therefore belong the set of riders, and their routes to the workplace
and back home must be carefully constructed and balanced to ensure
that every rider is covered on their way to work and guaranteed a ride
back home. These additional requirements make the CTSP unique and
particularly challenging. The CTSP also uses a lexicographic objective
function that first minimizes the number of cars and then the total
travel distance.

The paper proposes two exact algorithms for the CTSP: A Route-Enumeration Algorithm (REA) which exhaustively searches for feasible routes from all possible trip combinations before route selection is optimized with a mixed-integer program (MIP), and a Branch-and-Price algorithm (BPA) which uses column generation and a pricing algorithm based on dynamic programming. On top of the two algorithms, the paper highlights key characteristics of the CTSP differentiating it from the DARP that allow its routes to be enumerated by the REA. While the BPA builds on conventional techniques to solve the CTSP via column generation, it introduces a wait-time relaxation technique, which is a novel alternative to the weak and strong dominance relations proposed by \cite{gschwind2015} for finding feasible routes that simultaneously satisfy time-window and ride-duration constraints in the pricing problem. The paper also proposes a time-limited, root-node heuristic which is derived from the BPA and demonstrates its capability to produce high-quality solutions for medium to large problem instances within a 10-minute time span, \emph{making it well suited for time-constrained scenarios within an operational setting}. Finally, the paper proposes a clustering algorithm to decompose large-scale problems by spatially grouping commuters based on their home locations, which is then used to generate problem instances for evaluating the algorithms from a real-world dataset of commute trips from the city of Ann Arbor, Michigan.

The remainder of this paper is organized as follows. The next section
first briefly reviews relevant literature, and it is followed by
Section \ref{sec:notation_prelims} which introduces the terminology
and assumptions used throughout this work. Section
\ref{sec:comm_trip_share_prob} then provides a formal definition and a
mathematical formulation of the CTSP, followed by Section
\ref{sec:route_enum_algo} which introduces the first algorithm to
solve the problem, the REA. Next, Section \ref{sec:branch_price_algo}
describes the second algorithm, the BPA together with its derived root-node heuristic, while the clustering
algorithm is presented in Section \ref{sec:clustering}. Lastly,
computational results are reported in Section \ref{sec:exp_results},
after which concluding remarks are provided in Section
\ref{sec:conclusion}.

\section{Related Work}\label{sec:related_work}
The Vehicle Routing Problem with Time Windows (VRPTW) seeks a set of minimum-cost routes for a fleet of vehicles that start and end at a central depot. The routes serve a set of customers with specific demands and time windows describing allowable service times, and they must ensure that each customer is served exactly once and the capacity of the vehicles are not exceeded. The problem has been widely studied in the literature and is known to be NP-hard, since finding a solution for a fixed fleet has been shown to be NP-complete \citep{savelsbergh1985}. Nevertheless, various methods, from metaheuristics like \cite{taillard1997} and \cite{braysy2005} to exact solution approaches based on Lagrangian relaxation \citep{kohl1997, kallehauge2006} or column generation \citep{desrochers1992, kohl1999}, have been proposed to efficiently solve it. An extensive review of the problem may be obtained from \cite{cordeau2002}.

\cite{dumas1991} introduced the Pickup and Delivery Problem with Time Windows (PDPTW) to satisfy transportation requests requiring both pickup and delivery. It generalizes the VRPTW by introducing additional pairing and precedence constraints that require each route to serve the pickup location before the delivery location of the same customer, and it is solved using a column generation algorithm which utilizes dynamic programming to solve its pricing subproblem. The DARP, which is commonly used to model door-to-door transportation services for the disabled and the elderly, builds upon the PDPTW by introducing ride-duration limit constraints for each customer. As humans are being transported in the DARP instead of goods, customer ride time becomes an essential quality of service criterion. Methods proposed to solve small and medium-sized instances of the problem include heuristics \citep{jaw1986, bodin1986}, metaheuristics \citep{cordeau2003a, ritzinger2016}, and exact algorithms \citep{cordeau2006, gschwind2015}. The problem has also been extensively surveyed by \cite{cordeau2003b} and \cite{cordeau2007}.

When column generation is used to solve the various generalizations of the VRPTW, its pricing subproblem typically involves solving an Elementary Shortest Path Problem with Resource Constraints (ESPPRC). Since the problem has been proven to be NP-hard in the strong sense by \cite{dror1994}, most works have resorted to relaxing the elementary path requirement to result in a Shortest Path Problem with Resource Constraints (SPPRC) which admits a pseudo-polynomial algorithm. Non-elementary paths are then tackled using various methods; for instance, \cite{desrosiers1984} and \cite{dumas1991} have them eliminated in the integer solution of the restricted master problem, while \cite{desrochers1992} and \cite{irnich2006} opted for a middle ground approach by performing $2$- and $k$-cycle elimination respectively. Exact algorithms have also been proposed for solving the ESPPRC, e.g. by \cite{feillet2004} and \cite{chabrier2006}. Popular methods for solving the SPPRC and ESPPRC utilize dynamic programming, for instance the label-correcting algorithm of \cite{desrosiers1983} which is based on the Ford-Bellman-Moore algorithm, the label-setting algorithm of \cite{desrochers1988a} which generalizes Dijktra's algorithm, or the generalized label-setting algorithm for multiple resource constraints of \cite{desrochers1988b}. Methods using Lagrangian relaxation (e.g., \cite{beasley1989, borndrfer2001}) or constraint programming (e.g., \cite{rousseau2004}) have also been explored. An in-depth overview of the SPPRC is provided by \cite{irnich2005}.

More recently, the availability of large-scale datasets like the New York City (NYC) Taxi and Limousine Commission (TLC) trip record, which contains data for over 1 billion taxi trips in NYC recorded since January 2009, has driven research towards ride sharing for on-demand transportation. For instance, \cite{santi2014} introduced the notion of shareability graphs in an attempt to quantify the benefits of sharing these taxi rides. \cite{alonso2017} then built upon the shareability graph idea to mathematically model the on-demand ride sharing problem and propose an anytime optimal algorithm to solve it. \cite{agatz2012} provided an overview of planning considerations and issues of dynamic ride-sharing, classified different variations of ride-sharing problems, including single- and multi-modal versions, and reviewed related optimization models and approaches to address them. On the other hand, \cite{mourad2019} took a broader view of shared mobility and surveyed optimization approaches for a wider range of applications, from prearranged to real-time problem settings that even includes combined transportation of people and freight. To our knowledge, \cite{hasan2018} are the first to focus vehicle routing optimization on commute trips. They introduced the CTSP, explored the performance of various optimization models that enforce different sets of driver and commuter matching constraints, and discovered that commuter matching flexibility, i.e., their willingness to be matched with different drivers and passengers daily, is key for an effective ride-sharing platform. This paper extends their work by first refining their best performing ride-sharing model, introducing two algorithms to optimize the model, and performing an extensive comparative analysis of both algorithms using a high-fidelity, real-world dataset.

\section{Notation and Preliminaries}\label{sec:notation_prelims}
A trip $t=<o,dt,d,at>$ consists of an origin $o$, a desired departure time $dt$, a destination $d$, and a desired arrival time $at$. On any day, a commuter $c$ makes two trips: a trip to the workplace, $t_c^+$, and a trip back home, $t_c^-$. These trips are referred to henceforth as inbound and outbound trips respectively. A route $r$ is a sequence of origin and destination locations from a set of inbound or outbound trips whereby each origin and destination from the set is visited exactly once. For instance, a possible route for trips $t_1=<o_1,dt_1,d_1,at_1>$ and $t_2=<o_2,dt_2,d_2,at_2>$ is $r=o_2\rightarrow o_1\rightarrow d_1\rightarrow d_2$. An inbound route covers only inbound trips and an outbound route covers only outbound trips. Each route $r$ serves a set of riders $\mathcal{C}_r$ and has a driver $D_r\in \mathcal{C}_r$. The driver must be the rider residing at the start location of the route. For instance, commuter 2 must be the driver of route $o_2\rightarrow o_1\rightarrow d_1\rightarrow d_2$. The total number of riders in the vehicle at any point along a route cannot exceed its capacity. 

\begin{definition}[Valid Route]
	A valid route $r$ visits $o_c$ before $d_c$ for every rider $c\in \mathcal{C}_r$, starts at $o_{D_r}$ and ends at $d_{D_r}$, and respects the vehicle capacity.
\end{definition}

The paper assumes that commuters sharing rides are willing to tolerate
some inconvenience in terms of deviations to their trips' desired
departure and arrival times as well as in terms of extensions to the 
ride durations of their individual trips. Therefore, a time window $[a_i,b_i]$ is
constructed around the desired times and is associated with each
pickup or drop-off location $i$, where $a_i$ and $b_i$ denote the
earliest and latest times at which service may begin at $i$
respectively, and a duration limit $L_c$ is associated with each
commuter $c$ to denote her maximum ride duration. In this paper, $T_i$
denotes the time at which service begins at location $i$, $s_i$ is the
service duration at $i$, $pred(i)$ denotes the location on a route
visited just before $i$, and $\tau_{(i,j)}$ is the estimated travel
time for the shortest path between locations $i$ and $j$.

\begin{definition}[Feasible Route]
A feasible route $r$ is a valid route that has pickup and drop-off
times $T_i\in [a_i,b_i]$ for each location $i\in r$ and ensures the
ride duration of each commuter $c\in \mathcal{C}_r$ does not exceed
$L_c$.
\end{definition}

\noindent\begin{minipage}{\textwidth}
\begin{flalign}
\min & \quad T_{d_{D_r}} - T_{o_{D_r}} & \label{eqn:minimize_trip_duration} \\
\text{s.t.}\nonumber & \\
      & \quad a_{o_c}\leq T_{o_c} \leq b_{o_c} \qquad \forall c\in \mathcal{C}_r  \label{eqn:time_window_origin} \\
& \quad T_{d_c} \leq b_{d_c} \qquad \forall c\in \mathcal{C}_r \label{eqn:time_window_destination} \\
& \quad T_{pred(o_c)} + s_{pred(o_c)} + \tau_{(pred(o_c),o_c)} \leq T_{o_c}  \qquad \forall c\in \mathcal{C}_r\setminus\{D_r\}  \label{eqn:travel_time_origin} \\
& \quad T_{pred(d_c)} + s_{pred(d_c)} + \tau_{(pred(d_c),d_c)} = T_{d_c}  \qquad \forall c\in \mathcal{C}_r \label{eqn:travel_time_destination} \\
& \quad  T_{d_c} - (T_{o_c} + s_{o_c}) \leq L_c \qquad \forall c\in \mathcal{C}_r \label{eqn:ride_duration_limit} 
\end{flalign}
\vspace{\parskip}
\end{minipage}

Determining if a valid route $r$ is feasible amounts to solving the
route-scheduling problem of \eqref{eqn:minimize_trip_duration}--\eqref{eqn:ride_duration_limit}. Its objective is to minimize the total
duration of the route. Constraints \eqref{eqn:time_window_origin} and
\eqref{eqn:time_window_destination} are time-window constraints for
pickup and drop-off locations respectively, while constraints
\eqref{eqn:travel_time_origin} and \eqref{eqn:travel_time_destination}
describe compatibility requirements between pickup/drop-off times and
travel times between consecutive locations along the route. Finally,
constraints \eqref{eqn:ride_duration_limit} specify the ride-duration
limit for each rider. Note that constraints
\eqref{eqn:travel_time_origin} allow waiting at pickup locations, and
constraints \eqref{eqn:time_window_origin} and
\eqref{eqn:time_window_destination} implicitly limit the trip duration
of rider $c$ by $(b_{d_c} - a_{o_c})$.

The route validity requirement specifies route structural constraints
which enforce pairing and precedence of origins and
destinations, vehicle capacity, and the driver role, whereas the
feasibility requirement specifies time-window and ride-duration limit
constraints which are temporal in nature in addition to those for
route validity. Lastly, this work assumes utilization of a homogeneous
fleet of vehicles with capacity $K$ to serve all trips, and that all
travel times and distances satisfy the triangle inequality.

\section{The Commute Trip Sharing Problem}\label{sec:comm_trip_share_prob}

The CTSP aims at finding a set of minimum-cost feasible routes to
cover all inbound and outbound trips of a set of commuters
$\mathcal{C}$ while ensuring the set of drivers for inbound and
outbound routes are identical. Let $\Omega^+$ and $\Omega^-$ denote
the set of all feasible inbound and outbound routes respectively, and
$c_r$ denote the cost of route $r$. The CTSP formulation uses a binary
variable $X_r$ to indicate whether a route $r\in\Omega^+\cup\Omega^-$
is selected, a binary constant $\alpha_{r,i}$ which is equal to 1 iff
route $r$ serves rider $i$ (i.e., $\alpha_{r,i}=1$ iff
$i\in\mathcal{C}_r$), and a binary constant $\beta_{r,i}$ which is
equal to 1 iff rider $i$ is the driver of route $r$
(i.e., $\beta_{r,i}=1$ iff $i=D_r$). The problem formulation
is given by \eqref{eqn:minimize_total_cost}--\eqref{eqn:route_selection}.

\noindent\begin{minipage}{\textwidth}
\begin{flalign}
\min & \quad \sum_{r\in\Omega^+\cup\Omega^-} c_r X_{r} & \label{eqn:minimize_total_cost} \\
\text{s.t.} \nonumber & \\
& \quad \sum_{r\in\Omega^+} \alpha_{r,i} X_{r} = 1 \qquad \forall i\in\mathcal{C} \label{eqn:inbound_coverage} \\
& \quad \sum_{r\in\Omega^-} \alpha_{r,i} X_{r} = 1 \qquad \forall i\in\mathcal{C} \label{eqn:outbound_coverage} \\
& \quad \sum_{r\in\Omega^+} \beta_{r,i} X_{r} - \sum_{\hat{r}\in\Omega^-} \beta_{\hat{r},i} X_{\hat{r}} = 0 \qquad \forall i\in\mathcal{C} \label{eqn:driver_balance} \\
& \quad X_{r} \in \{0,1\} \qquad \forall r\in\Omega^+\cup\Omega^- \label{eqn:route_selection}
\end{flalign}
\vspace{\parskip}
\end{minipage}

The model features a lexicographic objective that first minimizes the
number of cars and then the total distance. It is rewritten into a
single objective by appropriate weighting of the two sub-objectives.
The cost $c_r$ penalizes the total distance of route $r$ and heavily
penalizes its selection. Let $\delta_{(i,j)}$ denote the {\em
  distance} of the shortest path between nodes $i$ and $j$. $c_r$ is
  then given by the addition of variable and fixed costs of the route:
\begin{equation}
c_r = \hat{c}_r + \bar{c}
\end{equation}
where the variable and fixed costs, $\hat{c}_r$ and $\bar{c}$, are given by:
\begin{equation}
\hat{c}_r = \sum_{(i,j)\in r}\delta_{(i,j)}
\end{equation}
\begin{equation}
\bar{c} = M \max_{r\in\Omega^+\cup\Omega^-}{\sum_{(i,j)\in r} \delta_{(i,j)}}
\end{equation}
%
%
where $M$ is a large number. In practice, $M$ is set to 1000, which is sufficiently large to ensure that the number of selected routes is first minimized followed by their total distance. Constraints \eqref{eqn:inbound_coverage} and \eqref{eqn:outbound_coverage} enforce coverage of each rider's inbound and outbound trips by exactly one route each, while constraints \eqref{eqn:driver_balance} ensure drivers of inbound and outbound routes are identical. The set-partitioning problem of \eqref{eqn:minimize_total_cost}--\eqref{eqn:route_selection} is referred to as the master problem (MP) from this point forth.

The CTSP is essentially a vehicle routing problem with driver, capacity, time-window, pairing, precedence, ride-duration, and driver constraints, making it most similar to the DARP. However, the key distinctions of the CTSP are:
\begin{enumerate}[label={(\alph*)}]
	\item Drivers in the CTSP are members of the set of riders, i.e., $D_r\in\mathcal{C}_r$. This leads to driver constraints which require routes to start and end at the drivers' origins and destinations respectively, whereas requests in the DARP are served by shared vehicles whose routes begin and end at a central depot.
	\item The set of drivers for inbound and outbound routes needs to be balanced, leading to constraints \eqref{eqn:driver_balance} in the MP. These constraints add another layer of complexity which is not present in the DARP.
\end{enumerate}
Therefore, the CTSP can also be seen as a DARP with additional constraints.

\section{The Route-Enumeration Algorithm}\label{sec:route_enum_algo}
One approach to solve the CTSP is by enumerating all routes in $\Omega^+\cup \Omega^-$ before solving the MP with a MIP solver. The REA supports this approach by exhaustively searching for these routes from all possible trip combinations. Let $\mathcal{T}^+$ and $\mathcal{T}^-$ denote all inbound and outbound trips taken by the set of commuters $\mathcal{C}$ respectively, i.e., $\mathcal{T}^+ = \{\,t_c^+\,:\,c\in\mathcal{C}\,\}$ and $\mathcal{T}^- = \{\,t_c^-\,:\,c\in\mathcal{C}\,\}$. Without loss of generality, Algorithm \ref{alg:AllFeasibleInboundRoutes} summarizes how $\Omega^+$ is obtained from $\mathcal{T}^+$ using a homogeneous fleet of vehicles with capacity $K$. 

Routes of all individual trips from $\mathcal{T}^+$ are first added to $\Omega^+$ (lines 2--3). To obtain feasible routes covering more than 1 trip, an index $k$ is first set to the number of shared trips desired, after which all $k$-combinations of trips from  $\mathcal{T}^+$ (denoted by $\mathcal{Q}_k$) are enumerated (lines 4--5). For each trip combination $q\in\mathcal{Q}_k$, the set of valid routes for the combination is then enumerated. For instance, let $k=2$, $q=\{t_1,t_2\}$, and $t_1=<o_1,dt_1,d_1,at_1>$ and $t_2=<o_2,dt_2,d_2,at_2>$. The set of valid routes for $q$ is $\{\,o_1\rightarrow o_2\rightarrow d_2\rightarrow d_1,\,o_2\rightarrow o_1\rightarrow d_1\rightarrow d_2\,\}$. $\Omega^v_q$ denotes the set of all valid routes for $q$ and $\mathcal{C}_q$ denotes the set of all riders making the trips in $q$.

\begin{algorithm}[!t]
	\caption{Route-Enumeration Algorithm for $\Omega^+$}\label{alg:AllFeasibleInboundRoutes}
	\begin{algorithmic}[1]
		\Require $\mathcal{T}^+,K$
		\State $\Omega^+ \leftarrow \text{\O}$
		\For{each $t_c^+\in\mathcal{T}^+$}
		\State $\Omega^+ \leftarrow \Omega^+ \cup \{o_c^+\rightarrow d_c^+\}$
		\EndFor
		\For{$k=2$ to $K$}
		\State $\mathcal{Q}_k \leftarrow \{$all $k$-combinations of $\mathcal{T}^+\}$
		\For{each $q\in \mathcal{Q}_k$}
		\State $\Omega^v_q \leftarrow \{$all valid routes of $q\}$
		\For{each $c\in \mathcal{C}_q$}
		\State $\Omega_{\text{temp}} \leftarrow \text{\O}$
		\For{each $r\in\Omega^v_q\,:\,D_r=c$}
		\If{$feasible(r)$}
		\State $\Omega_{\text{temp}} \leftarrow \Omega_{\text{temp}} \cup \{r\}$
		\EndIf
		\EndFor
		\State $\Omega^+ \leftarrow \Omega^+ \cup \{\argmin_{r\in\Omega_{\text{temp}}}\sum_{(i,j)\in r}\delta_{(i,j)}\}$
		\EndFor
		\EndFor
		\EndFor
		\State \textbf{return} $\Omega^+$
	\end{algorithmic}
\end{algorithm}

The algorithm then iterates over every rider $c\in\mathcal{C}_q$ and considers only routes in $\Omega^v_q$ where $c$ is the driver, i.e., $\{r\in\Omega^v_q\,:\,D_r=c\}$ (lines 8--10). A function $feasible(r)$, which solves the route-scheduling problem of \eqref{eqn:minimize_trip_duration}--\eqref{eqn:ride_duration_limit} on route $r$ and returns a Boolean value indicating whether $r$ is feasible, is then utilized to identify feasible routes to be stored in a temporary set $\Omega_{\text{temp}}$. Only the route with the shortest travel distance from $\Omega_{\text{temp}}$ is then added to $\Omega^+$ (line 13). Note that this step is optional. It is done to reduce the size of $\Omega^+$ with the knowledge that only one route may be selected for each driver covering $\mathcal{C}_q$ in a feasible solution to the MP. The route with minimal travel distance is chosen knowing that the secondary objective of the MP is to minimize the total distance of selected routes.

In practice, the search procedure in lines 7--13 may be executed more efficiently via a depth-first search implementation which uses the length of the best feasible route to prune the search space and through parallel execution of the search procedure for all $q\in\mathcal{Q}_k$ since they are independent of each other. The procedure of exploring all $k$-combinations is repeated with increasing values of $k$ from 2 up to the vehicle capacity $K$ to completely enumerate $\Omega^+$. $\Omega^-$ is obtained by repeating the algorithm on $\mathcal{T}^-$.


The fact that drivers are commuters themselves is the key characteristic of the CTSP: It allows its routes to be exhaustively enumerated. Indeed, drivers must complete their trips within their time windows and, most importantly, their trips are subject to ride-duration constraints. As a result, in general, a route typically consists of three phases: A pickup phase where the driver picks up passengers, a driving phase where the vehicle travels to the destination, and a drop-off phase where the driver drops off all the passengers before ending her trip. After the drop-offs, the driver has no time  to go back and pick up another set of passengers due to her time-window and ride-duration constraints. This permits the REA to consider only routes that contain up to $k \leq K$ passengers (line 4 in Algorithm \ref{alg:AllFeasibleInboundRoutes}) to enumerate all possible routes, and $K$ is typically small. In contrast, the DARP uses dedicated drivers who are not subject to any ride-duration constraints and can serve riders throughout the day. Therefore it cannot restrict attention to routes with only $k \leq K$ passengers: The number of passengers in a route is not limited by the capacity of the vehicle, but by the total number of travelers. 

\section{The Branch-and-Price Algorithm}\label{sec:branch_price_algo}

The BPA combines existing techniques with some novel elements to solve the CTSP. At its core is a conventional column-generation algorithm which utilizes a restricted master problem (RMP)---the linear relaxation of the MP defined on a subset of all feasible routes $\Omega^{+\prime}\cup\Omega^{-\prime}$---and solves a pricing subproblem (PSP) to identify new feasible routes with negative reduced costs. The PSP solves several dynamic programs that search for resource-constrained shortest paths representing the feasible routes. A bi-level branching strategy tailored specifically for the CTSP is then employed for obtaining integer solutions.

\emph{This work introduces a novel wait-time relaxation technique that not only obtains feasible routes that simultaneously satisfy time-window and ride-duration constraints in the PSP, but also guarantees elementarity of the routes.} It proposes utilization of a resource that models trip durations excluding wait times which allow the dynamic programs to produce preliminary routes with minimal reduced costs that first satisfy a set of constraints necessary for route feasibility. The feasibility of the preliminary routes are then evaluated with the inclusion of wait times, and infeasible ones are added to a set of forbidden paths whose members are prevented from subsequent discovery via the dynamic-programming approach of \cite{pugliese2013b}.

\subsection{The Pricing Subproblem}
The PSP is responsible for finding new feasible routes with negative
reduced costs. Letting $\pi_i^+$, $\pi_i^-$, and $\sigma_i$ denote the
optimal duals of constraints \eqref{eqn:inbound_coverage},
\eqref{eqn:outbound_coverage}, and \eqref{eqn:driver_balance} of the RMP
respectively, the reduced cost of an inbound route $r^+$ is given by:
\begin{equation}
rc_{r^+} = c_{r^+} - \sum_{i\in\mathcal{C}_{r^+}} \pi_i^+ - \sigma_{D_{r^+}}
\end{equation}
while that of an outbound route $r^-$ is given by:
\begin{equation}
rc_{r^-} = c_{r^-} - \sum_{i\in\mathcal{C}_{r^-}} \pi_i^- + \sigma_{D_{r^-}}
\end{equation}

\noindent
These routes are obtained by considering each rider $d\in\mathcal{C}$ as the 
driver of an inbound route $r_d^+$ and an outbound route $r_d^-$, and then
finding such routes with minimum reduced costs. To obtain
these routes, the algorithm builds a pair of graphs $\mathcal{G}_d^+$
and $\mathcal{G}_d^-$ for each $d\in\mathcal{C}$. In the following,
$\mathcal{G}$ denotes the set of all constructed graphs, i.e.,
$\mathcal{G} = \{\mathcal{G}_d^+ \,:\,d\in\mathcal{C}\} \cup
\{\mathcal{G}_d^-\,:\,d\in\mathcal{C}\}$.  Without loss of generality,
the presentation outlines how a route $r_d^+$ with minimal reduced
cost is found from $\mathcal{G}_d^+$.

First, let $n=|\mathcal{C}|$, and $\mathcal{O}^+=\{1,\cdots,n\}$ and
$\mathcal{D}^+=\{n+1,\cdots,2n\}$ denote the sets of all
origin and destination nodes respectively. The origin and
destination of rider $i$ are then represented by nodes $i$ and $n+i$
respectively. The graph $\mathcal{G}_d^+ =
\{\mathcal{N}_d^+,\mathcal{A}_d^+\}$ is built with nodes
$\mathcal{N}_d^+=\mathcal{O}^+\cup\mathcal{D}^+$ and fully-connected
edges $\mathcal{A}_d^+$. A ride-duration limit $L_i$ and a demand
$\kappa_i$, representing the number of riders to be picked up at node
$i$, are then associated with each node $i\in\mathcal{O}^+$, a time
window $[a_i,b_i]$ and a service duration $s_i$ are associated with
each node $i\in\mathcal{N}_d^+$, and a travel time $\tau_{(i,j)}$ and
a reduced cost $c_{(i,j)}$ are associated with each edge
$(i,j)\in\mathcal{A}_d^+$. Letting $\gamma^+(i)$ and $\gamma^-(i)$
denote the set of outgoing and incoming edges of node $i$, the edge
costs are defined as follows so that the total cost of any path from
$d$ to $n+d$ is equivalent to $rc_{r_d^+}$:
\begin{equation}\label{eqn:inbound_edge_costs}
c_{(i,j)}=
\begin{cases}
\bar{c} + \delta_{(i,j)} - \pi_i^+ - \sigma_d&\qquad\forall (i,j)\in\gamma^+(d)\\
\delta_{(i,j)} - \pi_i^+&\qquad\forall i\in\mathcal{O}^+\setminus\{d\},\, \forall (i,j)\in\gamma^+(i)\\
\delta_{(i,j)}&\qquad\forall i\in\mathcal{D}^+,\, \forall (i,j)\in\gamma^+(i)
\end{cases}
\end{equation}

\begin{figure}[!t]
	\FIGURE
	{\centering
		\includegraphics[width=0.45\linewidth]{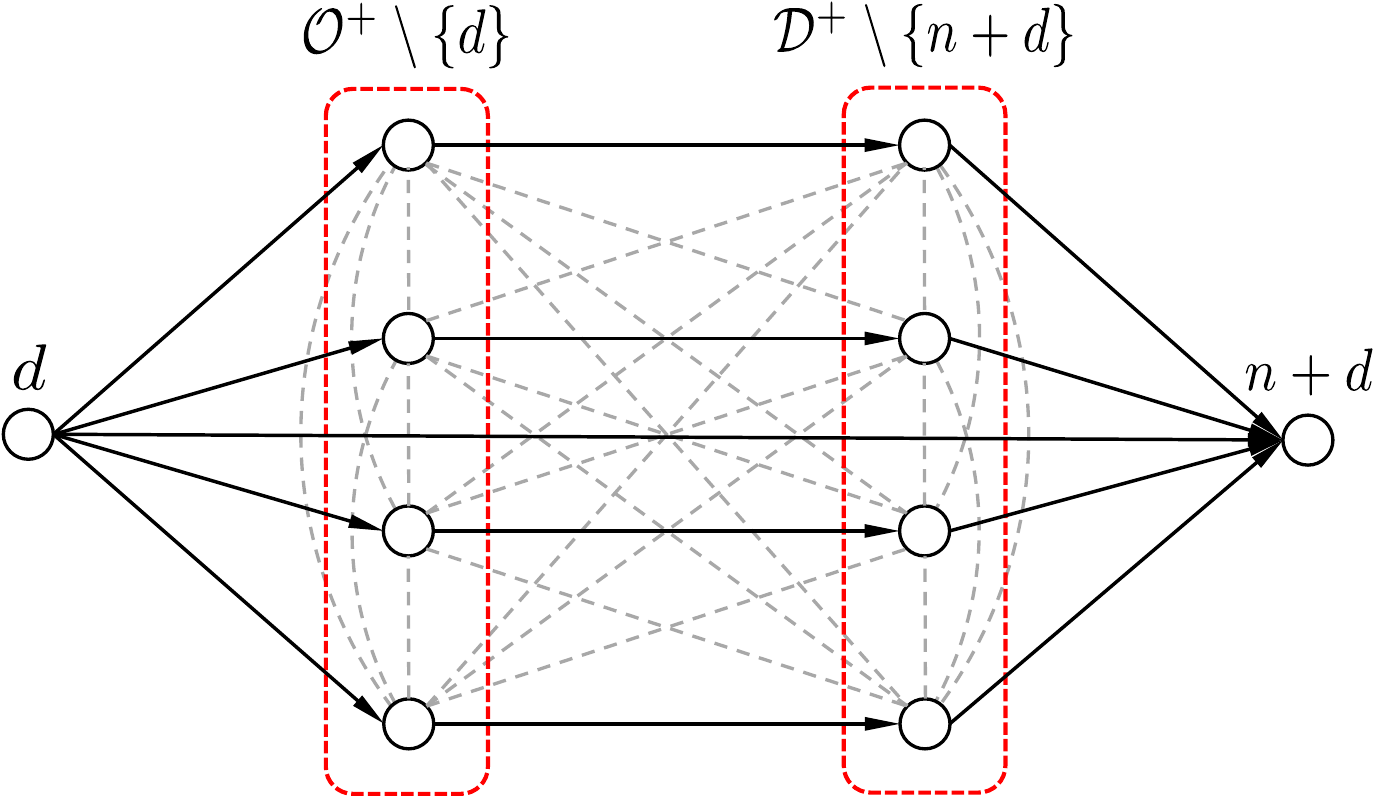}}
	{Graph $\mathcal{G}_d^+$ After Application of Edge Elimination Rules (a) and (b) from Section \ref{sec:edge_elimination} (Each Dotted Line Represents a Pair of Bidirectional Edges).\label{fig:Graph}}
	{}
\end{figure}

Similarly, letting $\mathcal{O}^-$ and $\mathcal{D}^-$ denote the sets of all outbound origin and destination nodes, the graph $\mathcal{G}_d^- =
\{\mathcal{N}_d^-,\mathcal{A}_d^-\}$ is built with nodes
$\mathcal{N}_d^-=\mathcal{O}^-\cup\mathcal{D}^-$ and fully-connected
edges $\mathcal{A}_d^-$, and the costs of edges $(i,j)\in\mathcal{A}_d^-$ are defined as follows to ensure the total cost of any path from $d$ to $n+d$ in $\mathcal{G}_d^-$ is equal to $rc_{r_d^-}$:
\begin{equation}\label{eqn:outbound_edge_costs}
c_{(i,j)}=
\begin{cases}
\bar{c} + \delta_{(i,j)} - \pi_i^- + \sigma_d&\qquad\forall (i,j)\in\gamma^+(d)\\
\delta_{(i,j)} - \pi_i^-&\qquad\forall i\in\mathcal{O}^-\setminus\{d\},\, \forall (i,j)\in\gamma^+(i)\\
\delta_{(i,j)}&\qquad\forall i\in\mathcal{D}^-,\, \forall (i,j)\in\gamma^+(i)
\end{cases}
\end{equation}
A priori feasibility constraints, further detailed in Section \ref{sec:edge_elimination}, are then applied to identify and eliminate edges that cannot belong to any feasible route. Figure \ref{fig:Graph} provides a sketch of $\mathcal{G}_d^+$ after application of several of these edge elimination rules.

The minimum-reduced-cost $r_d^+$ is then obtained by finding the least-cost feasible path from $d$ to $n + d$ in $\mathcal{G}_d^+$. Recall that for the path to be feasible, it must satisfy the time-window, capacity, pairing, precedence, ride-duration and driver constraints. The problem is therefore an Elementary Shortest Path Problem with Resource Constraints (ESPPRC) which is known to be NP-hard \citep{dror1994}. While the driver constraint is enforced by construction by making $d$ the source and $n+d$ the target of the shortest-path problem, the remaining constraints are implemented by introducing and enforcing constrained resources in a resource-constrained shortest path algorithm (RCSPA) which is further elaborated in Section \ref{sec:rcspa}. On the whole, the PSP involves solving $2n$ independent ESPPRCs to produce up to $2n$ feasible routes with negative reduced costs.

\subsection{Time Windows Tightening and Edge Elimination}\label{sec:edge_elimination}

Pre-processing of the time-window, precedence, pairing, capacity,
ride-duration limit, and driver constraints makes it possible to
identify edges that cannot belong to any feasible route which may 
then be removed from $\mathcal{G}$. Without loss of
generality, the following description focuses on edge elimination for
$\mathcal{G}_d^+$.

Prior to determining infeasible edges, the time windows of all
nodes are tightened by sequentially reducing their upper and lower
bounds using the following rules introduced by \cite{dumas1991}.

\begin{itemize}
	\item $b_i = \min\{b_i, b_{n+d} - s_i - \tau_{(i,n+d)}\},\,\forall i\in\mathcal{D}^+ \setminus \{n+d\}$
	\item $b_i = \min\{b_i, b_{n+i} - s_i - \tau_{(i,n+i)}\},\,\forall i\in\mathcal{O}^+ \setminus \{d\}$
	\item $a_i = \max\{a_i, a_d + s_d + \tau_{(d,i)}\},\,\forall i\in\mathcal{O}^+ \setminus \{d\}$
	\item $a_i = \max\{a_i, a_{i-n} + s_{i-n} + \tau_{(i-n,i)}\},\,\forall i\in\mathcal{D}^+ \setminus \{n+d\}$
\end{itemize}

\noindent
The following constraints and rules, derived by combining those
proposed by \cite{dumas1991} and \cite{cordeau2006}, are then applied
to identify and eliminate infeasible edges:

\begin{enumerate}[label={(\alph*)}]
	\item Driver: Edges $\{(d,n+i), (i,d), (i,n+d), (n+i,d), (n+d,i), (n+d,n+i)\,:\,i\in\mathcal{O}^+ \setminus \{d\}\}$. 
	\item Pairing and precedence: Edges $\{(n+i,i)\,:\,i\in\mathcal{O}^+\}$.
	\item Capacity: Edges $\{(i,j), (j,i), (i,n+j), (j,n+i), (n+i,n+j), (n+j,n+i)\,:\,i,j\in\mathcal{O}^+ \wedge i\neq j \wedge \kappa_i + \kappa_j > K\}$.
	\item Time windows: Edges $\{(i,j)\,:\,(i,j)\in\mathcal{A}_d^+ \wedge a_i + s_i + \tau_{(i,j)} > b_j\}$.
	\item Ride-duration limit: Edges $\{(i,j),(j,n+i)\,:\,i\in\mathcal{O}^+ \wedge j\in\mathcal{N}_d^+ \wedge i\neq j \wedge \tau_{(i,j)} + s_j + \tau_{(j,n+i)} > L_i\}$.
	\item Pairing, time windows, and ride-duration limit:
	\begin{itemize}
		\item Edges $\{(i,n+j)\,:\,i,j\in\mathcal{O}^+ \wedge i\neq j \wedge \neg feasible(j\rightarrow i\rightarrow n+j\rightarrow n+i)\}$.
		\item Edges $\{(n+i,j)\,:\,i,j\in\mathcal{O}^+ \wedge i\neq j \wedge \neg feasible(i\rightarrow n+i\rightarrow j\rightarrow n+j)\}$.
		\item Edges $\{(i,j)\,:\,i,j\in\mathcal{O}^+ \wedge i\neq j \wedge \neg feasible(i\rightarrow j\rightarrow n+i\rightarrow n+j) \wedge \neg feasible(i\rightarrow j\rightarrow n+j\rightarrow n+i)\}$.
		\item Edges $\{(n+i,n+j)\,:\,i,j\in\mathcal{O}^+ \wedge i\neq j \wedge \neg feasible(i\rightarrow j\rightarrow n+i\rightarrow n+j) \wedge \neg feasible(j\rightarrow i\rightarrow n+i\rightarrow n+j)\}$.
	\end{itemize}
\end{enumerate}

Note that the rules in (f) utilize the $feasible$ function introduced earlier to determine if a partial route satisfies time-window and ride-duration limit constraints. For instance, the first says edge $(i,n+j)$ is infeasible if route $j\rightarrow i\rightarrow n+j\rightarrow n+i$ is infeasible. Edge elimination rules for $\mathcal{G}_d^-$ are obtained by replacing $\mathcal{O}^+$, $\mathcal{D}^+$, and $\mathcal{A}_d^+$ in the above rules with $\mathcal{O}^-$, $\mathcal{D}^-$, and $\mathcal{A}_d^-$ respectively.

\subsection{The Resource-Constrained Shortest Path Algorithm}\label{sec:rcspa}

The PSP uses an RCSPA based on the label-setting dynamic program proposed
by \cite{desrochers1988b} to find the least-cost feasible path from any graph in $\mathcal{G}$, i.e., one that satisfies time-window, capacity, pairing, precedence, and ride-duration constraints. The path searched by this algorithm is identical to that sought in the PSP of the DARP by \cite{gschwind2015}. Their method incorporated novel dominance rules in the labeling procedure to directly enforce all constraints in the dynamic program. 

On the other hand, the RCSPA presented here first searches for the {\em minimum-cost, feasible route that ignores the wait times.} The routes that are infeasible with respect to the wait times are then pruned in a second step. This procedure is motivated by the fact that the optimal values for the wait times requires knowledge of the complete route, which is only known at the end of the search. By relaxing the wait times, the dynamic program first finds a candidate route which is later evaluated for feasibility with respect to the wait times once it is complete. Moreover, subsequent empirical evaluations revealed that for the problem instances considered, an overwhelming majority of the candidate routes are feasible with the inclusion of wait times, and the resources utilized in the algorithm are also capable of guaranteeing generation of elementary paths.

The RCSPA can therefore be seen as a middle ground approach between the method by \cite{ropke2006} which completely relaxes the ride-duration constraint in the PSP and prevents selection of paths that may violate the constraint through infeasible path elimination constraints in the RMP, and that of \cite{gschwind2015} which directly enforces all constraints in the dynamic program of the PSP. Without loss of generality, this section describes the
algorithm for $\mathcal{G}_d^+$.

\subsubsection{Label definition}

Let $\mathcal{P}_l^k$ denote the $k\textsuperscript{th}$ path from the
source $d$ to node $l$. A label $\mathcal{L}_l^k$ with five resources
$(c_l^k,T_l^k,\mathcal{C}_l^k,\mathcal{R}_l^k,\mathcal{W}_l^k)$ is
associated with each $\mathcal{P}_l^k$. $c_l^k$ represents the total
cost of edges in $\mathcal{P}_l^k$, i.e., $c_l^k =
\sum_{(i,j)\in\mathcal{P}_l^k} c_{(i,j)}$, whereas $T_l^k$ is the time
at which service at node $l$ begins for
$\mathcal{P}_l^k$. $\mathcal{C}_l^k$ denotes the set of riders on the
vehicle right after visiting node $l$ on $\mathcal{P}_l^k$. It is
equivalent to the set of pickup nodes visited on $\mathcal{P}_l^k$
whose corresponding drop-off nodes have yet to be visited. On the
other hand, $\mathcal{R}_l^k$ denotes the set of all riders that have
been picked up by $\mathcal{P}_l^k$ after visiting node $l$. It is
equivalent to the set of all pickup nodes visited by
$\mathcal{P}_l^k$. Finally, $\mathcal{W}_l^k$ is the set of trip
durations, \emph{excluding wait times}, for each rider in
$\mathcal{R}_l^k$. Letting $\mathcal{P}_l^k(m)$ denote the set of
edges from $\mathcal{P}_l^k$ on which rider $m$ is on the vehicle and
$w_l^k(m)$ be the trip duration of rider $m$ excluding wait times on
$\mathcal{P}_l^k$, i.e., $w_l^k(m) = \sum_{(i,j)\in \mathcal{P}_l^k(m)}
s_i + \tau_{(i,j)}$, then $\mathcal{W}_l^k =
\{w_l^k(m)\,:\,m\in\mathcal{R}_l^k\}$. The load $Y_l^k$ of a vehicle
after visiting node $l$ on path $\mathcal{P}_l^k$ can be easily
obtained from $Y_l^k = \sum_{i\in\mathcal{C}_l^k}
\kappa_i$. Therefore, $\mathcal{L}_l^k$ contains sufficient
information to ensure $\mathcal{P}_l^k$ satisfies pairing, precedence,
time-window, and capacity constraints. While resource
$\mathcal{W}_l^k$ is not sufficient for verifying compliance to the
ride-duration limit for each rider, it does provide a lower bound to
each ride duration which must necessarily satisfy the limit for
$\mathcal{P}_l^k$ to be feasible.

\subsubsection{Label extension}
$\mathcal{L}_l^k$ is maintained using a forward dynamic program. In the label-setting algorithm, an attempt is made to extend $\mathcal{L}_l^k$ along edge $(l,j)$ to produce label $\mathcal{L}_j^{k'}$ for path $\mathcal{P}_j^{k'}$. The resources in $\mathcal{L}_j^{k'}$ are calculated as follows:
\begin{equation}\label{eqn:cost_extend}
c_j^{k'} = c_l^k + c_{(l,j)}
\end{equation}
\begin{equation}\label{eqn:service_time_extend}
T_j^{k'} = 
\begin{cases}
\max\{a_j, T_l^k + s_l + \tau_{(l,j)}\}&\qquad\text{if }j\in\mathcal{O}^+\\
T_l^k + s_l + \tau_{(l,j)}&\qquad\text{otherwise}
\end{cases}
\end{equation}
\begin{equation}\label{eqn:commuters_in_vehicle_extend}
\mathcal{C}_j^{k'}=
\begin{cases}
\mathcal{C}_l^k\cup\{j\}&\qquad\text{if }j\in\mathcal{O}^+\\
\mathcal{C}_l^k\setminus\{j-n\}&\qquad\text{otherwise}
\end{cases}
\end{equation}
\begin{equation}
\mathcal{R}_j^{k'}=\mathcal{R}_l^k\cup\{j\}\qquad\text{if }j\in\mathcal{O}^+
\end{equation}
\begin{equation}
w_j^{k'}(j) = 0 \qquad \text{if } j\in\mathcal{O}^+ \wedge j\notin\mathcal{R}_l^k
\end{equation}
\begin{equation}\label{eqn:ride_duration_ignore_wait_extend}
w_j^{k'}(i) = w_l^k(i) + s_l + \tau_{(l,j)} \qquad \forall i\in\mathcal{C}_l^{k}
\end{equation}
The extension is performed if and only if:
\begin{equation}\label{eqn:service_time_extend_condition}
T_j^{k'}\leq b_j,
\end{equation}
\begin{equation}
j\notin\mathcal{C}_l^k\qquad\text{if }j\in\mathcal{O}^+,
\end{equation}
\begin{equation}
j-n\in\mathcal{C}_l^k\qquad\text{if }j\in\mathcal{D}^+,
\end{equation}
\begin{equation}
\sum_{i\in\mathcal{C}_j^{k'}} \kappa_i \leq K,\text{ and}
\end{equation}
\begin{equation}\label{eqn:ride_duration_ignore_wait_extend_condition}
w_j^{k'}(i)-s_i\leq L_i\qquad\forall i\in\mathcal{C}_l^{k}.
\end{equation}
Constraints \eqref{eqn:service_time_extend_condition}--\eqref{eqn:ride_duration_ignore_wait_extend_condition} list conditions that are necessary to ensure feasibility of $\mathcal{P}_j^{k'}$. Note that if $w_j^{k'}(i)-s_i$, which constitutes rider $i$'s ride duration excluding wait times and hence is the lower bound to her ride duration, is already exceeding $L_i$, then $L_i$ will certainly be exceeded if wait times were included. Therefore conditions in \eqref{eqn:ride_duration_ignore_wait_extend_condition} are necessary but not sufficient in enforcing the ride-duration limit constraint for each rider. 

The algorithm is initialized by path $\mathcal{P}_d^1$ whose label $\mathcal{L}_d^1=(0,a_d,\{d\},\{d\},\{0\})$, and a preliminary solution is given by path $\mathcal{P}_{n+d}^{k^*}$ whose cost $c_{n+d}^{k^*}$ is minimal and whose resource $\mathcal{C}_{n+d}^{k^*}=\text{\O}$. Note that a non-elementary path may result if the graph contains a negative-cost cycle. However, such paths may be eliminated by setting the ride-duration limit of each rider to be less than twice the ride duration of her direct trip, i.e., $L_i<2\tau_{(i,n+i)}+s_i$.
\begin{proposition}
	Non-elementary paths will not be generated by the RCSPA if $L_i<2\tau_{(i,n+i)}+s_i$ for each $i\in\mathcal{O}^+$.	
\end{proposition}

\proof{Proof.}
Suppose a non-elementary path is generated by the RCSPA. On the path, there must exist at least one rider $i$ who is served more than once. For such riders, both $i$ and $n+i$ must be visited more than once with $i$ preceding $n+i$ each time and $n+i$ being visited first before $i$ is visited again due to the pairing and precedence constraints. As a result, resource $w_{n+d}^{k^*}(i)\geq 2(s_i + \tau_{(i,n+i)})$ and therefore $w_{n+d}^{k^*}(i)-s_i\geq2\tau_{(i,n+i)}+s_i$. If $L_i<2\tau_{(i,n+i)}+s_i$, then $w_{n+d}^{k^*}(i)-s_i>L_i$. Condition \eqref{eqn:ride_duration_ignore_wait_extend_condition} is thus violated, causing the path to not be extended.\Halmos
\endproof

Also note that as the restrictions on $\mathcal{W}_l^k$ are not sufficient for ensuring satisfaction of the ride-duration constraints, $\mathcal{P}_{n+d}^{k^*}$ may be infeasible. Therefore, an additional step needs to be performed to verify the feasibility of $\mathcal{P}_{n+d}^{k^*}$.

\subsubsection{Forbidding paths violating the ride-duration limit}

Feasibility of the preliminary solution $\mathcal{P}_{n+d}^{k^*}$ with
the inclusion of wait times can be verified using the $feasible$
function once the path is complete. A feasible path
$\mathcal{P}_{n+d}^{k^*}$ represents the optimal solution to the
ESPPRC of the PSP. While
empirical evaluations revealed that the vast majority of preliminary
routes found ($>99\%$ of the paths found) are feasible, infeasible
paths are still discovered on rare occasions. In such cases, the
infeasible path is added to a set of forbidden paths associated with
the graph, after which the RCSPA is executed again repeatedly to
generate newer paths until a feasible one is found.

The shortest path problem with forbidden paths
\citep{villeneuve2005,pugliese2013a,pugliese2013b} is a method that
has been successfully applied for handling constraints which are hard
or impossible to model as resources. This work exploits this idea to
properly enforce the ride-duration limit constraints by preventing
infeasible preliminary routes from being discovered by the RCSPA
again. The dynamic-programming approach of \cite{pugliese2013b} is
employed for this purpose since it fits well into the label-setting
framework.

Firstly, let $\mathcal{F}_d^+$ denote the set of forbidden paths for $\mathcal{G}_d^+$. Also let $\dot{f}$ denote the first edge of a forbidden path $f$, $|f|$ denote the total number of edges on the path, and $h_l^k(f)$ denote the number of consecutive edges of $f$ starting from $\dot{f}$ that is present in $\mathcal{P}_l^k$. To forbid paths in $\mathcal{F}_d^+$ from being discovered by the RCSPA, an additional resource $\mathcal{H}_l^k = \{h_l^k(f)\,:\,f\in\mathcal{F}_d^+\}$ is introduced to the label so that $\mathcal{L}_l^k=(c_l^k,T_l^k,\mathcal{C}_l^k,\mathcal{R}_l^k,\mathcal{W}_l^k,\mathcal{H}_l^k)$. During label extension along edge $(l,j)$, $\mathcal{H}_j^{k'}$ is calculated as follows:
\begin{equation}
h_j^{k'}(f)=
\begin{cases}
1&\quad\text{if }(l,j)\in f \wedge h_l^k(f)=0 \wedge (l,j)=\dot{f}\\
0&\quad\text{if }(l,j)\in f \wedge h_l^k(f)=0 \wedge (l,j)\neq\dot{f}\\
h_l^k(f)+1&\quad\text{if }(l,j)\in f \wedge h_l^k(f)\geq 1 \wedge consecutive(\mathcal{P}_l^k,(l,j),f)\\
0&\quad\text{if }(l,j)\in f \wedge h_l^k(f)\geq 1 \wedge \neg consecutive(\mathcal{P}_l^k,(l,j),f)\\
0&\quad\text{if }(l,j)\notin f
\end{cases}
\quad\forall f\in\mathcal{F}_d^+
\end{equation}
$consecutive(\mathcal{P}_l^k,(l,j),f)$ is a function that returns true if there exists a set of consecutive edges in path $\{\mathcal{P}_l^k,(l,j)\}$ ending with $(l,j)$ that exactly matches a set of consecutive edges in path $f$ starting from $\dot{f}$, and returns false otherwise. The extended resource must then satisfy the following constraints:
\begin{equation}
h_j^{k'}(f)\leq |f|-1\qquad\forall f\in\mathcal{F}_d^+
\end{equation}
since $\mathcal{P}_j^{k'}$ would contain a forbidden path otherwise. The resource is initialized with $h_d^1(f)=0$ for each $f\in\mathcal{F}_d^+$. Resource $\mathcal{H}_l^{k}$ prevents the RCSPA from discovering infeasible preliminary routes stored in $\mathcal{F}_d^+$ again, thus ensuring the algorithm's solution is always feasible.

\subsubsection{Label elimination}
As efficiency of the label-setting algorithm increases with the amount of eliminated labels, a label and its associated path is eliminated if it is established that the label cannot belong to either an optimal or a feasible solution.  Firstly, dominance rules are applied to determine if a label does not belong to an optimal solution.

\begin{definition}[Label Domination]
$\mathcal{L}_l^k$ dominates $\mathcal{L}_l^{k'}$ if and only if:
\begin{equation}
c_l^k \leq c_l^{k'},
\end{equation}
\begin{equation}
T_l^k \leq T_l^{k'},
\end{equation}	
\begin{equation}
\mathcal{C}_l^k \subseteq \mathcal{C}_l^{k'},
\end{equation}	
\begin{equation}
w_l^k(i) \leq w_l^{k'}(i) \qquad \forall i\in\mathcal{C}_l^k,\text{ and}
\end{equation}	
\begin{equation}
h_l^k(f) \leq h_l^{k'}(f) \qquad \forall f\in\mathcal{F}_d^+.
\end{equation}		
\end{definition}

If $\mathcal{L}_l^{k'}$ is dominated by $\mathcal{L}_l^k$, then $\mathcal{L}_l^{k'}$ and its associated path $\mathcal{P}_l^{k'}$ cannot belong to an optimal solution to the ESPPRC as every feasible extension to $\mathcal{P}_l^{k'}$ is also applicable to $\mathcal{P}_l^{k}$ at an equal or lower cost. Therefore $\mathcal{L}_l^{k'}$ and $\mathcal{P}_l^{k'}$ may be eliminated.

Next, the following rules are applied to identify labels that cannot belong to a feasible solution:
\begin{enumerate}[label={(\alph*)}]
	\item $\mathcal{L}_l^k$ such that $\mathcal{C}_l^k\setminus\{d\}\neq\text{\O}$ is eliminated if there exists $i\in\mathcal{C}_l^k\setminus\{d\}$ where the path extension $l\rightarrow n+i\rightarrow n+d$ is infeasible.
	\item $\mathcal{L}_l^k$ such that $|\,\mathcal{C}_l^k\setminus\{d\}\,|\geq 2$ is eliminated if there exists $i,j\in\mathcal{C}_l^k\setminus\{d\} \wedge i\neq j$ where path extensions $l\rightarrow n+i\rightarrow n+j\rightarrow n+d$ and $l\rightarrow n+j\rightarrow n+i\rightarrow n+d$ are both infeasible.
\end{enumerate}
For the rules above, feasibility of path extensions are verified by
checking if they satisfy the time window and ride-duration
constraints excluding wait times, i.e., by checking if each node along
the extension satisfies conditions \eqref{eqn:service_time_extend_condition} and
\eqref{eqn:ride_duration_ignore_wait_extend_condition}. The rules are inspired
by the notion of non-post-feasible labels introduced by
\cite{dumas1991}. They are essentially heuristics which check if at
least one (for rule (a)) or two (for rule (b)) of the riders on the
vehicle, excluding the driver, can be delivered to their destinations
while respecting their time windows and ride-duration limits if wait
times are ignored. While not sufficient, these conditions are
necessary for the feasibility of any extension to $\mathcal{P}_l^{k}$,
and they result in the elimination of a large amount of infeasible
labels in practice.

\subsection{Obtaining an Integer Solution}

The unique structure of the MP lets us infer a few properties about
its solution. Firstly, since the total number of selected inbound
routes must match that of outbound routes in any solution, the total
number of selected routes in an integer solution must be even, i.e.,
$\sum_{r\in\Omega^{+}\cup\Omega^{-}} X_r \in
\{2a\,:\,a\in\mathbb{Z}_{\geq 0}\}$. Secondly, since only integral
distances are used in this work, all routes costs and consequently the
objective value of an integer solution must also be integral, i.e.,
$\sum_{r\in\Omega^{+}\cup\Omega^{-}} c_r X_r \in \mathbb{Z}_{\geq 0}$.

These two properties are leveraged to obtain an integer solution
should the optimal solution of the RMP not be integral. Let $\chi^*$
denote the total number of selected routes at convergence,
i.e., $\chi^*=\sum_{r\in\Omega^{+\prime}\cup\Omega^{-\prime}} X_r$, and
$z^*$ denote the objective value at convergence. If $\chi^*$ is not an
even integer, the following cut is introduced to the RMP to round up
the total number of selected routes to the nearest even integer.
\begin{equation}\label{eqn:trip_count_cut}
\sum_{r\in\Omega^{+\prime}\cup\Omega^{-\prime}} X_r \geq 2\ceil*{\frac{\chi^*}{2}}
\end{equation}
The dual of the cut is appropriately transfered to the PSP and the column-generation procedure is resumed until convergence again. If $z^*$ is not integral at this point, another cut is added to the RMP to round up its objective value to the nearest integer:
\begin{equation}\label{eqn:objective_value_cut}
\sum_{r\in\Omega^{+\prime}\cup\Omega^{-\prime}} c_r X_r \geq \ceil*{z^*}
\end{equation}
Once again, the dual of the cut is transferred to the PSP and the column-generation procedure is resumed until convergence. If the solution of the RMP is still not integral at this stage, then a branch-and-bound tree needs to be explored whereby additional columns may be generated at each tree node. 

A bi-level branching scheme is employed for the branch-and-bound tree, whereby integrality of driver selection is enforced in the first level and integrality of edge flow is enforced in the second. In the first level, let $V_i$ be a variable that indicates whether rider $i$ is selected as the driver in a solution. It is given by:
\begin{equation}
V_i = \sum_{r\in\Omega^{+\prime}} \beta_{r,i} X_r \qquad \forall i\in\mathcal{C}
\end{equation}
In an integral solution, all $V_i$s must be binary. Therefore if they are not, a fractional $V_i$ is selected and two branches are created; one fixing it to 0 and another fixing it to 1. The branch decision of $V_i = 0$ is enforced in the RMP by removing columns where rider $i$ is the driver, i.e., $\{r\in\Omega^{+\prime}\cup\Omega^{-\prime}\,:\,D_r = i\}$, while it is enforced in the PSP by not solving the ESPPRC on graphs where rider $i$ is the driver, i.e., $\mathcal{G}^+_i$ and $\mathcal{G}^-_i$. To enforce $V_i = 1$, the following cut is introduced to the RMP:
\begin{equation}\label{eqn:driver_selection_cut}
\sum_{r\in\Omega^{+\prime}\,:\,D_r=i} X_r = 1
\end{equation}
while ensuring its dual is properly incorporated into the PSP. No additional steps are needed to enforce the branch decision in the PSP since the ESPPRCs on graphs $\mathcal{G}^+_i$ and $\mathcal{G}^-_i$ are already being solved by default. 

If all $V_i$s are binary and the solution of the RMP is still fractional, then a second branching scheme based on that proposed by \cite{desrochers1992} is utilized. In the second level, let $\omega(i,j)$ denote the set of all routes utilizing edge $(i,j)$, i.e., $\omega(i,j) = \{r\in\Omega^{+\prime}\cup\Omega^{-\prime}\,:\,(i,j)\in r\}$, and let $F_{(i,j)}$ be the flow variable for edge $(i,j)$ that indicates if node $i$ should be served before node $j$ in a solution. It is given by:
\begin{equation}
F_{(i,j)} = \sum_{r\in\omega(i,j)} X_r \qquad \forall (i,j)\in \{\mathcal{A}_d^+\cup\mathcal{A}_d^-:d\in\mathcal{C}\}
\end{equation} 
Also let $\mathcal{A}^+$ and $\mathcal{A}^-$ denote the set of edges from all inbound and outbound graphs respectively, i.e., $\mathcal{A}^+ = \{\mathcal{A}_d^+\,:\,d\in\mathcal{C}\}$, $\mathcal{A}^- = \{\mathcal{A}_d^-\,:\,d\in\mathcal{C}\}$. In an integer solution, all $F_{(i,j)}$s must be binary. In a fractional solution however, one of the following cases may occur:
\begin{enumerate}[label={(\alph*)}]
	\item $F_{(i,j)}$ for all $(i,j)\in\mathcal{A}^+$ are binary, but there exists $(u,v)\in\mathcal{A}^-$ such that $F_{(u,v)}$ is fractional.
	\item $F_{(u,v)}$ for all $(u,v)\in\mathcal{A}^-$ are binary, but there exists $(i,j)\in\mathcal{A}^+$ such that $F_{(i,j)}$ is fractional.
	\item There exist $(i,j)\in\mathcal{A}^+$ and $(u,v)\in\mathcal{A}^-$ such that both $F_{(i,j)}$ and $F_{(u,v)}$ are fractional. 
\end{enumerate}
If either case (a) or (b) occurs, then an edge $(i,j)$ whose flow is fractional is selected (from either $\mathcal{A}^+$ or $\mathcal{A}^-$ depending on the case) and two branches are created; one setting $F_{(i,j)}=0$ and another setting $F_{(i,j)}=1$. Should case (c) occurs, then two edges whose flows are fractional are selected, $(i,j)\in\mathcal{A}^+$ and $(u,v)\in\mathcal{A}^-$, and four branches are created with the following decisions:
\begin{enumerate}
	\item $F_{(i,j)} = 0 \wedge F_{(u,v)} = 0$.
	\item $F_{(i,j)} = 0 \wedge F_{(u,v)} = 1$.
	\item $F_{(i,j)} = 1 \wedge F_{(u,v)} = 0$.
	\item $F_{(i,j)} = 1 \wedge F_{(u,v)} = 1$. 
\end{enumerate}

$F_{(i,j)}=0$ is enforced in the RMP by removing columns containing edge $(i,j)$, whereas in the PSP, edge $(i,j)$ is removed from all graphs to prevent columns containing it from being generated. To enforce $F_{(i,j)}=1$, edges in sets $\gamma^+(i)\setminus\{(i,j)\}$ and $\gamma^-(j)\setminus\{(i,j)\}$ are removed from all graphs in the PSP and columns containing the edges are correspondingly removed from the RMP.

In practice, cuts \eqref{eqn:trip_count_cut}, \eqref{eqn:objective_value_cut}, and \eqref{eqn:driver_selection_cut} are introduced into the RMP (one for every rider $i\in\mathcal{C}$ in the case of \eqref{eqn:driver_selection_cut}) from the very beginning with their right-hand sides initially set to $\geq 0$. The right-hand sides are then correspondingly updated to those shown in \eqref{eqn:trip_count_cut}, \eqref{eqn:objective_value_cut}, and \eqref{eqn:driver_selection_cut} as the algorithm progresses. Let $\mu$, $\nu$, and $\phi_i$ denote the duals of cuts \eqref{eqn:trip_count_cut}, \eqref{eqn:objective_value_cut}, and \eqref{eqn:driver_selection_cut} respectively. These duals are incorporated into the PSP by updating the costs of edges $(i,j)\in\mathcal{A}_d^+$ defined earlier in \eqref{eqn:inbound_edge_costs} to:
\begin{equation}\label{eqn:inbound_edge_costs_update}
c_{(i,j)}=
\begin{cases}
\bar{c}(1 - \nu) + \delta_{(i,j)}(1 - \nu) - \pi_i^+ - \sigma_d - \mu - \phi_i&\qquad\forall (i,j)\in\gamma^+(d)\\
\delta_{(i,j)}(1 - \nu) - \pi_i^+&\qquad\forall i\in\mathcal{O}^+\setminus\{d\},\, \forall (i,j)\in\gamma^+(i)\\
\delta_{(i,j)}(1 - \nu)&\qquad\forall i\in\mathcal{D}^+,\, \forall (i,j)\in\gamma^+(i)
\end{cases}
\end{equation}
and those of $(i,j)\in\mathcal{A}_d^-$ defined in \eqref{eqn:outbound_edge_costs} to:
\begin{equation}\label{eqn:outbound_edge_costs_update}
c_{(i,j)}=
\begin{cases}
\bar{c}(1 - \nu) + \delta_{(i,j)}(1 - \nu) - \pi_i^- + \sigma_d - \mu&\qquad\forall (i,j)\in\gamma^+(d)\\
\delta_{(i,j)}(1 - \nu) - \pi_i^-&\qquad\forall i\in\mathcal{O}^-\setminus\{d\},\, \forall (i,j)\in\gamma^+(i)\\
\delta_{(i,j)}(1 - \nu)&\qquad\forall i\in\mathcal{D}^-,\, \forall (i,j)\in\gamma^+(i)
\end{cases}
\end{equation}

\subsection{Implementation Strategies}
Several strategies are adopted in our implementation to reduce
execution time. Firstly, since the PSP involves solving at most $2n$
ESPPRCs which are independent, they are solved in parallel and
multiple columns are added to the RMP in each column-generation
iteration.

Secondly, to check the convergence of the column-generation phase, a
primal upper bound and a dual lower bound are maintained for the
optimal objective value, $z^*$. The objective value of the
RMP after each iteration, $z_\text{RMP}$, serves as the primal upper
bound while the lower bound proposed by \cite{lubbecke2005} is used as
the dual lower bound. It is given by $z_\text{LB} = z_\text{RMP} +
{rc}^* \lambda$, where ${rc}^*$ is the smallest reduced cost discovered
in the PSP and $\lambda$ is an upper bound to the number of selected
routes, i.e., $\lambda \geq \sum_{r\in\Omega^+\cup\Omega^-} X_r$. In this
case, it is easy to see that $\lambda$ can be chosen as $2n$.

Assume that $\chi_\text{RMP}$ and $\chi_\text{LB}$ are the upper and lower bounds to the total number of selected routes, obtained by considering only the fixed cost contributions to $z_\text{RMP}$ and $z_\text{LB}$ respectively. Since the number of selected routes must be even for an integer solution, the column generation is first suspended when $2 \ceil{\chi_\text{RMP} / 2} - \chi_\text{LB} < 2$. Cut \eqref{eqn:trip_count_cut} is then introduced to the MP to round the total to the nearest even integer after which the column generation is resumed. Since the optimal objective value of the MP must be integral, the column generation is terminated when $\ceil{z_\text{RMP}} - z_\text{LB} < 1$, after which cut \eqref{eqn:objective_value_cut} is introduced.

Finally, the branch-and-bound tree is explored depth-first to quickly obtain integer solutions. During tree exploration, a best integer solution may be obtained at any stage by solving the RMP as a MIP (in practice, this is only done for every 1,000 tree nodes explored beginning with the root node due to its potentially high expense). Let $z_\text{MIP}$ denote the objective value of the MIP solution, $z_\text{int}^*$ be that of the optimal integer solution sought, and $z_\text{min}^*$ be the smallest $z^*$ from all unexplored tree nodes. Since at any stage of tree exploration $z_\text{min}^* \leq z_\text{int}^* \leq z_\text{MIP}$, it is terminated when $z_\text{MIP} - z_\text{min}^* < 1$, at which point the optimal integer solution is given by the best integer solution.

\subsection{The Root-Node Heuristic}
To assess the algorithm's ability to produce high-quality solutions in an operational setting, a heuristic is conceived based on the BPA. It simply executes column generation at the root node of the branch-and-price tree within an allocated time budget $t_{\text{RMP}}$, and then finds an integer solution by solving the RMP as a MIP within another time budget $t_{\text{MIP}}$. The multi-objective function is simplified to only minimize the number of selected routes by setting route costs $\mathbf{c}_r \equiv \mathbf{1}$. The quality of the heuristic solution is assessed by calculating its optimality gap given by $(z_{\text{MIP}} - z_{\text{LB}}) / z_{\text{MIP}}$, where $z_{\text{MIP}}$ is the objective value of the MIP solution and $z_{\text{LB}}$ is its lower bound. $z_{\text{LB}}$ is given by the optimal objective value of the RMP at convergence, $z^{*}$. Should the RMP not converge within $t_{\text{RMP}}$, a lower bound to $z^{*}$ that is calculated using the method proposed by \cite{farley1990} is used instead. Farley's lower bound is given by:
\begin{equation}
z_\text{LB} = z_\text{RMP} \frac{c_{r^{\prime}}}{\boldsymbol{\pi}^{\intercal}\boldsymbol{a}_{r^{\prime}}}
\end{equation}
where $r^{\prime}=\argmin_{r\in\Omega^{+}\cup\Omega^{-}}\{c_r/\boldsymbol{\pi}^{\intercal}\boldsymbol{a}_r\,:\,\boldsymbol{\pi}^{\intercal}\boldsymbol{a}_r>0\}$, $\boldsymbol{\pi}$ is the dual optimal solution of the RMP, and $\boldsymbol{a}_r$ is the column of constraint coefficients of route $r$. The unit route costs simplify the lower bound to $z_\text{LB} = z_{\text{RMP}}/(1 - rc^{*})$.

An alternate variant of the heuristic which relaxes forbidden paths in the RCSPA is also considered. The consideration is made based on a couple of preliminary observations: (1) preliminary solutions to the RCSPA are very rarely infeasible, and (2) forbidding discovery of infeasible paths in the RCSPA is expensive. A consequence of this relaxation is that infeasible routes may be introduced into the RMP and therefore: (1) they will need to be filtered out before the RMP is solved as a MIP, and (2) the RMP may converge to a weaker lower bound, $z_\text{LB} \leq z^*$. Despite the potential loss in solution quality, the relaxation strategy may still be worthwhile as the loss may be very small and it may be outweighed by gains resulting from faster computation times.

\section{The Clustering Algorithm}\label{sec:clustering}

Commuters are grouped according to the neighborhoods they live in, using a clustering algorithm that groups no more than $N$ commuters together based on the spatial proximity of their home locations, before ride-sharing is optimized intra-cluster. This notion of breaking down the problem into smaller ones is not new; it is in line with the conclusion by \cite{agatz2012} that effective decomposition techniques will be necessary to make large-scale problems computationally feasible. Besides improving tractability by decomposing the problem into smaller, independent subproblems which could then be solved concurrently, the technique also limits the distance traveled by the driver when picking up and dropping off passengers and fosters intra-community interaction. In other words, it supports the notion of having riders commute with people living nearby, which was identified earlier as a desirable feature for the car-pooling platform. It must also be acknowledged that the tractability gained from this decomposition comes at a price: it precludes obtaining a global optimal solution. The trade-off is seen as a necessity however as empirical evaluations have shown that a global solution cannot be obtained for the dataset considered within a time frame that is reasonable for an operational setting, even with the faster root-node heuristic. The algorithm is therefore a mechanism by which problem instances of different sizes are generated for the computational experiments.

The algorithm treats commuters as points in $\mathbb{R}^2$ whose positions are specified by
the Cartesian coordinates of their homes. The algorithm is similar to
the $k$-means clustering algorithm \citep{lloyd1982} with the
exception of a small modification to its assignment step. The number
of clusters, $k=\ceil{|\mathcal{C}|/N}$, is first calculated where
$\mathcal{C}$ denotes the set of all commuters. Cluster centers are
then initialized using the $k$-means++ method by \cite{arthur2007},
whereby a center $\boldsymbol{u}_1$ is first selected uniformly at
random from $\mathcal{C}$. Let $S(\boldsymbol{x})$ denote the
Euclidean distance from point $\boldsymbol{x}$ to the nearest center
already selected. The $i\textsuperscript{th}$ center
$\boldsymbol{u}_i$ is then selected from $\mathcal{C}$ with
probability
$\frac{S(\boldsymbol{u}_i)^2}{\sum_{\boldsymbol{c}\in\mathcal{C}}S(\boldsymbol{c})^2}$
until $k$ centers are obtained.

An assignment step then assigns each point to its nearest center subject to a constraint that each center is assigned at most $N$ points. Let $\mathcal{U}$ denote the set of all cluster centers and $S(\boldsymbol{x},\boldsymbol{y})$ denote the Euclidean distance between points $\boldsymbol{x}$ and $\boldsymbol{y}$. The assignment step is performed by solving the generalized-assignment problem of \eqref{eqn:min_total_distance}--\eqref{eqn:assign_var}. It is defined in terms of a binary variable $x_{\boldsymbol{c},\boldsymbol{u}}$ which indicates whether commuter $\boldsymbol{c}$ is assigned to cluster center $\boldsymbol{u}$. Its objective function minimizes the total distance between commuters and their assigned cluster centers. Constraints \eqref{eqn:assignment} assigns each commuter to one cluster center, while constraints \eqref{eqn:size_limit} limit the number of commuters assigned to each center to $N$.

\noindent\begin{minipage}{\textwidth}
\begin{flalign}
\min & \quad \sum_{\boldsymbol{c}\in\mathcal{C}}\sum_{\boldsymbol{u}\in\mathcal{U}} S(\boldsymbol{c},\boldsymbol{u}) x_{\boldsymbol{c},\boldsymbol{u}} & \label{eqn:min_total_distance} \\
\text{s.t.}\nonumber & \\
& \quad \sum_{\boldsymbol{u}\in\mathcal{U}} x_{\boldsymbol{c},\boldsymbol{u}}=1\qquad\forall \boldsymbol{c}\in\mathcal{C} \label{eqn:assignment} \\
& \quad \sum_{\boldsymbol{c}\in\mathcal{C}} x_{\boldsymbol{c},\boldsymbol{u}}\leq N\qquad\forall \boldsymbol{u}\in\mathcal{U} \label{eqn:size_limit} \\
& \quad x_{\boldsymbol{c},\boldsymbol{u}}\in\{0,1\}\qquad\forall \boldsymbol{c}\in\mathcal{C},\forall \boldsymbol{u}\in\mathcal{U} \label{eqn:assign_var}
\end{flalign}
\vspace{\parskip}
\end{minipage}

After assignment, the coordinates of each cluster center is updated with the mean of the coordinates of all assigned commuters:
\begin{equation}
\boldsymbol{u} = \frac{\sum_{\boldsymbol{c}\in\mathcal{C}} x_{\boldsymbol{c},\boldsymbol{u}} \boldsymbol{c}}{\sum_{\boldsymbol{c}\in\mathcal{C}}x_{\boldsymbol{c},\boldsymbol{u}}} \qquad \forall \boldsymbol{u}\in\mathcal{U}
\end{equation} 
\noindent The assignment and update steps are repeated until the assignments stabilize, i.e., until the commuter-cluster center assignments stop changing, at which point the algorithm is terminated.

\section{Experimental Results}\label{sec:exp_results}

This section reports the computational results for the proposed
algorithms, as well as their effectiveness in reducing parking
pressure.

\subsection{Experimental Setting}

\begin{figure}[!t]
	\FIGURE
	{\centering
		\includegraphics[width=1.0\linewidth]{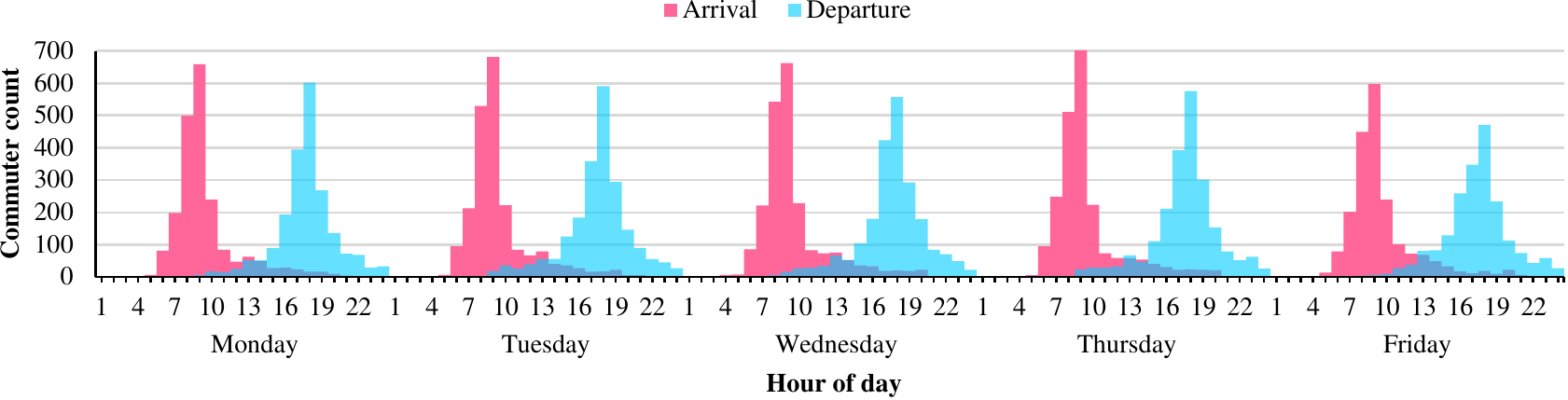}}
	{Commuting Patterns on Second Week of April 2017.\label{fig:AccessTimeDistribution}}
	{}
\end{figure}

The computational performance of the algorithms is evaluated using
problem instances derived from a real-world dataset consisting of
access information of 15 parking structures located in downtown Ann
Arbor, Michigan. The access information consists of the IDs, arrival
times, and departure times of customers of the parking structures
throughout the month of April 2017. This information is joined with
their home addresses to reconstruct their daily trips. It resulted in
the trip information of approximately 3,900 commuters living within
Ann Arbor's city limits (the region bounded by highways US-23, M-14,
and I-94), an area spanning 27 square miles, from which approximately
2,200 commute trips are made on a daily basis. Figure
\ref{fig:AccessTimeDistribution}, which shows the distribution of
arrival and departure times of this population over the busiest week
of the month (the second week), reveals the remarkable similarity and
consistency of their travel patterns over different
weekdays. Particularly notable are the peaks of the arrival and
departure time distributions which coincide with the typical 6--9 am
and 4--7 pm traffic peak hours.

Several assumptions are made regarding commuters using the trip
sharing platform. Firstly, it is assumed that, when requesting a
commute trip, rider $i$ would specify the desired arrival time at the
destination of her inbound trip $at_i^+$ and the desired departure
time $dt_i^-$ at the origin of her outbound trip.  This assumption is
consistent with that made in other DARP literature,
e.g. \cite{jaw1986}, \cite{cordeau2003a}, and \cite{cordeau2006}. It
is also assumed that the commuters are willing to tolerate a maximum
shift of $\pm\Delta$ to the desired times. Therefore, by treating the
arrival and departure times at the parking structures as the desired
times, time windows of $[a_{n+i}, b_{n+i}] = [at_i^+ - \Delta, at_i^+
  + \Delta]$ and $[a_i, b_i] = [dt_i^- - \Delta, dt_i^- + \Delta]$ are
associated with the destination of the inbound trip and the origin of
the outbound trip of rider $i$ respectively. Consequently, time
windows at the origin of the inbound trip and at the destination of
the outbound trip of rider $i$ are calculated using $[a_i, b_i] =
[a_{n+i} - s_i - L_i, b_{n+i} - s_i - \tau_{(i,n+i)}]$ and $[a_{n+i},
  b_{n+i}] = [a_i + s_i + \tau_{(i,n+i)}, b_i + s_i + L_i]$
respectively. It is also assumed that each commuter is willing to
tolerate at most an $R\%$ extension to her direct-ride duration,
i.e., $L_i = (1 + R)\cdot \tau_{(i,n+i)}$. This assumption is similar
to that made by \cite{hunsaker2002}. 

\subsection{Algorithmic Settings}

The clustering algorithm is used to construct problem instances of
varying sizes by varying $N$. Due to the non-deterministic nature of
its initialization step, the algorithm is executed 100 times for each
value of $N$, after which only the solution with the smallest assignment
objective value is selected. The shortest path, travel time estimate,
and travel distance estimate between any two locations are obtained
using GraphHopper's Directions API which uses data from
OpenStreetMap. All algorithms are implemented in C++ with
parallelization duties being handled by OpenMP. The
resource-constrained shortest path function from Boost 1.64.0's Graph
Library is used to implement the RCSPA, while Gurobi 7.5.1 is invoked
to solve all LPs and MIPs. The route fixed cost $\bar{c}$ is
obtained by making a very conservative overestimate of the longest
route length. The RMP of the BPA is initialized with the set of all
feasible single- and two-trip routes, which is generated using the REA
with $K = 2$. Each problem is solved on a high-performance computing
cluster with 12 cores of a 2.5 GHz Intel Xeon E5-2680v3 processor and
64 GB of RAM. Unless otherwise stated, a time limit of 12 hours is
applied to all problems and the best feasible solution is reported for
those that cannot be solved optimally within the time limit.

\subsection{Selecting Values for $\Delta$ and $R$}

Half of the time-window size, $\Delta$, and the ride-duration limit, $L_i = (1 + R)\cdot \tau_{(i,n+i)}$, directly influence the quality of service for rider $i$; the former represents the maximum amount of time by which the rider needs to shift (up or down) her desired arrival or departure time at a parking lot, whereas the latter represents the maximum amount of time the rider has to spend on the vehicle. Therefore, it is ideal for any rider to have the values of $\Delta$ and $R$ be as small as possible. However, doing so will also limit the potential for trip sharing. Indeed, selecting values for either parameter involves a trade-off between user convenience and trip shareability. A sensitivity analysis was therefore performed to study the impact of these two parameters on the vehicle reduction for the CTSP, by first applying the clustering algorithm with $N=100$ on commuters traveling on each of the first four weekdays of week 2, and then optimizing trip sharing within each cluster with the REA using $K = 4$ (the capacity of a car). The following values were used in the analysis: $\Delta = \{5,10,15\}$ mins and $R = \{25,50,75\}\%$. 

\begin{figure}[!t]
	\centering
	\begin{minipage}{.48\textwidth}
		\FIGURE
		{\centering
			\includegraphics[width=1.0\linewidth]{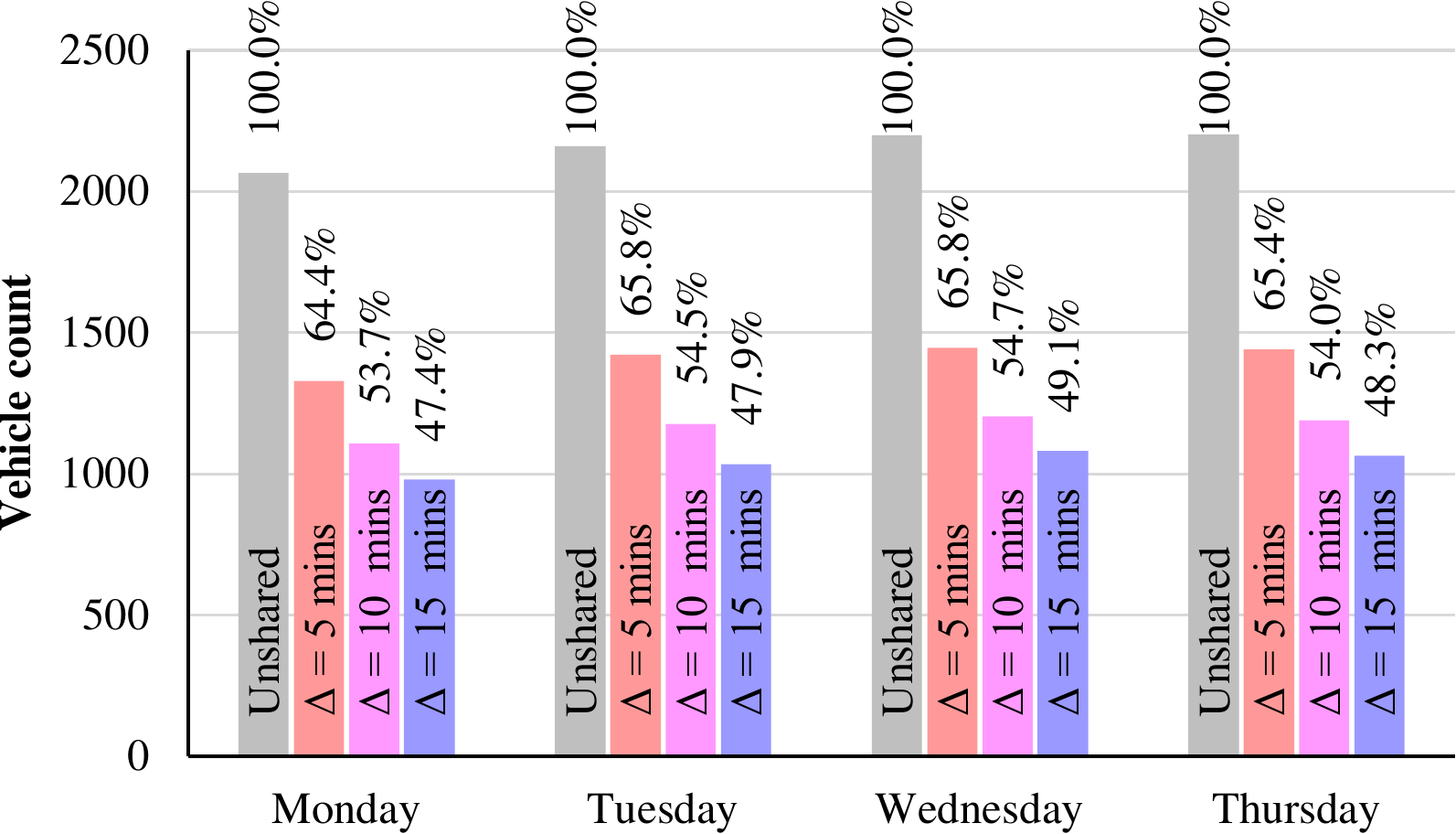}}
		{Effect of Increasing $\Delta$ on Total Vehicle Count ($R=50\%, K=4, N=100$).\label{fig:DeltaSensitivity}}
		{}
	\end{minipage}
	\hspace{0.02\textwidth}
	\begin{minipage}{.48\textwidth}
		\FIGURE
		{\centering
			\includegraphics[width=1.0\linewidth]{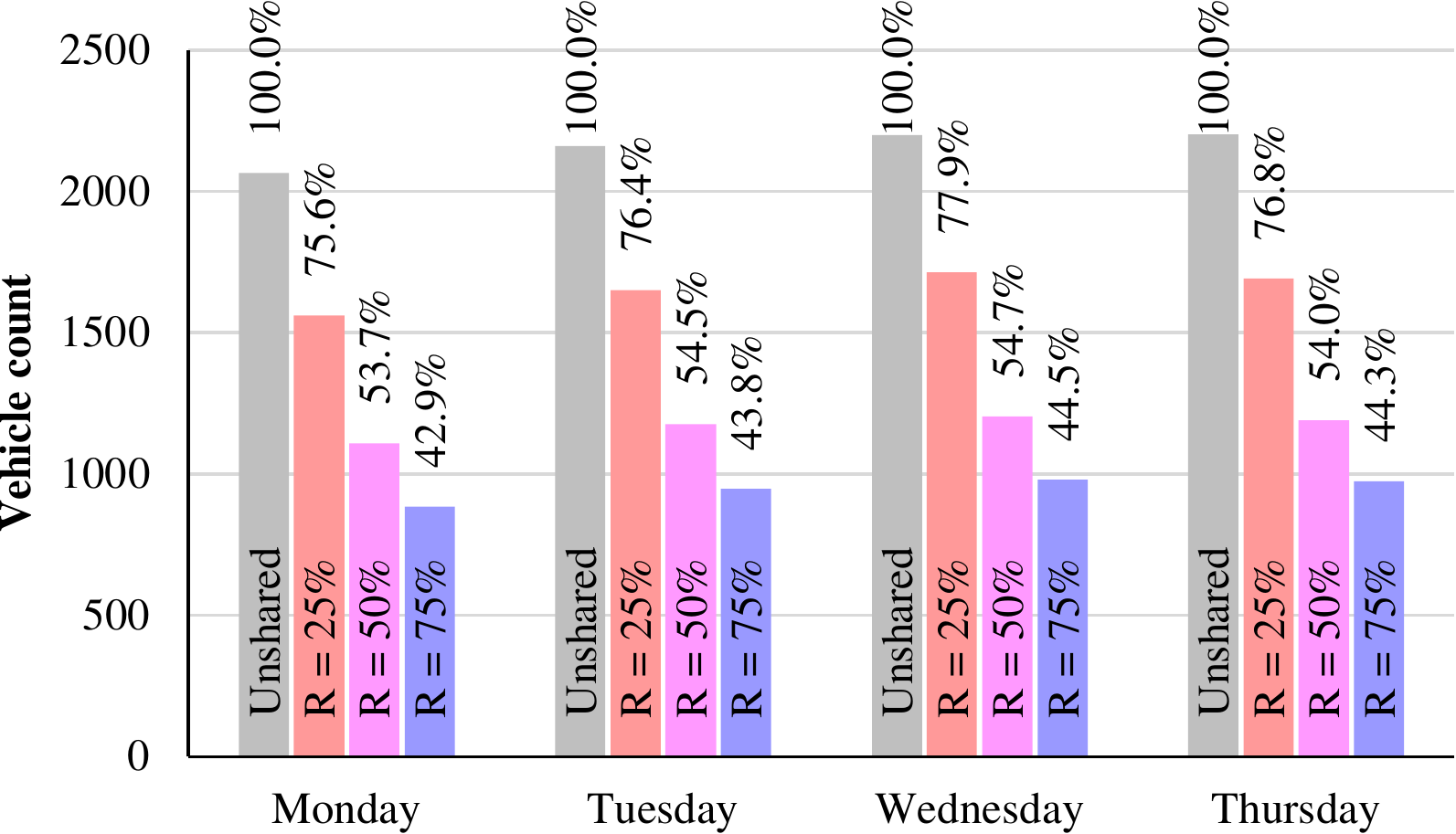}}
		{Effect of Increasing $R$ on Total Vehicle Count ($\Delta=10$ mins, $K=4,N=100$).\label{fig:RSensitivity}}
		{}
	\end{minipage}
\end{figure}

Figures \ref{fig:DeltaSensitivity} and \ref{fig:RSensitivity} summarize the analysis results for $\Delta$ and $R$ respectively. They include the required number of vehicles under no sharing conditions as well as the percentage of each vehicle count as a fraction of the no sharing count for additional perspective. The figures quantify the trade-off mentioned earlier; while the smaller values of $\Delta = 5$ mins and $R = 25\%$ may be convenient for the riders, the results indicate that the vehicle reduction potential is significantly hampered by these values. Vehicle reduction increases, albeit at a decreasing rate, as both $\Delta$ and $R$ are increased. While the largest values of both parameters produce the best vehicle reduction potential, they also demand the highest amount of tolerance to inconvenience from the riders. Therefore, $\Delta = 10$ mins and $R = 50\%$ were deemed to be the best compromise, as they still produced a sizeable amount of vehicle reduction while not inducing a significant amount of inconvenience to the riders. These values are therefore used in the rest of the experiments. 

\subsection{Vehicle Capacity Scaling}
\label{sec:veh_cap_scaling}

\begin{figure}[!t]
	\FIGURE
	{\centering
		\includegraphics[width=0.90\linewidth]{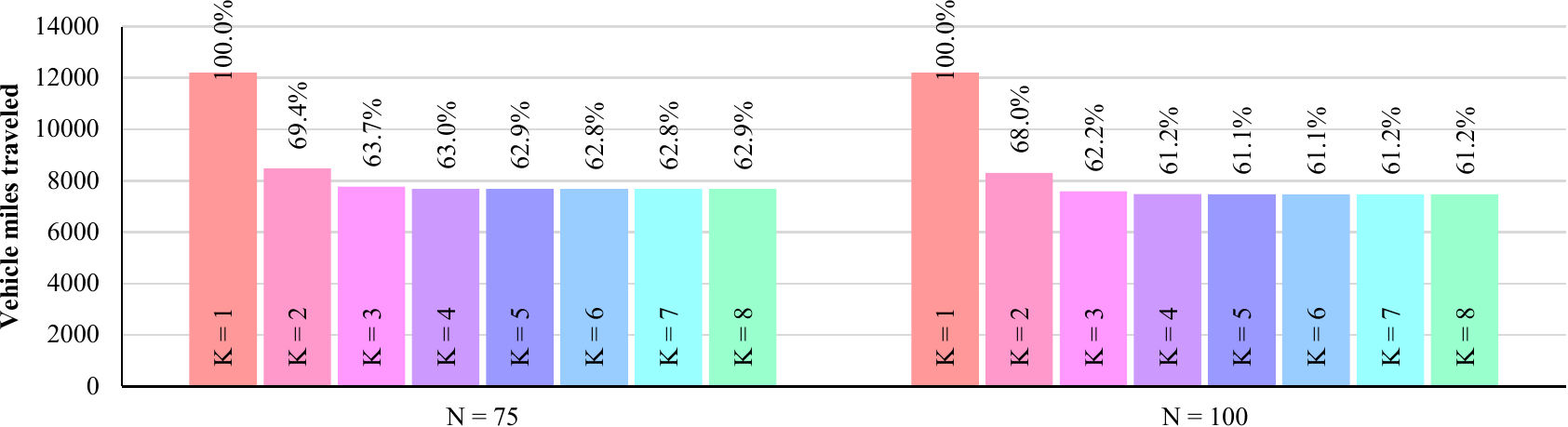}}
	{Effect of Increasing Vehicle Capacity on Total Route Distance.\label{fig:VehCapScaling_VehicleMilesTraveled}}
	{}
\end{figure}

\begin{figure}[!t]
	\FIGURE
	{\centering
		\includegraphics[width=0.90\linewidth]{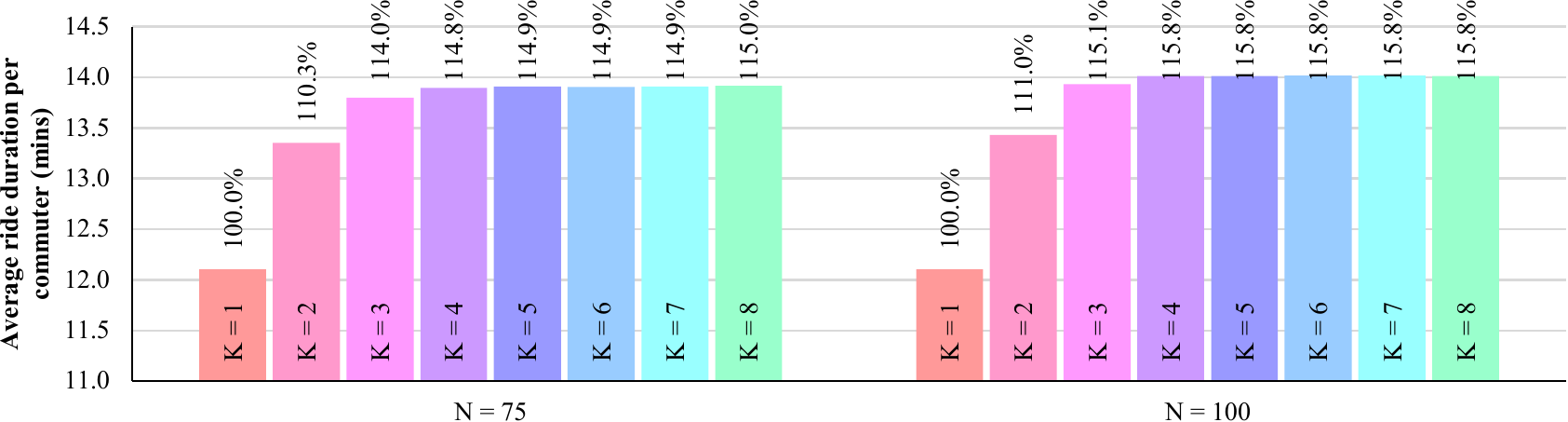}}
	{Effect of Increasing Vehicle Capacity on Average Ride Duration.\label{fig:VehCapScaling_AvgRideDuration}}
	{}
\end{figure}

The next set of computational experiments explores the scalability of
both algorithms with increasing vehicle capacity. A variety of
car-pooling programs provide small vans to commuters: These vans can
typically carry about 8 people and it is important to evaluate the
benefits of using such vehicles. Problem instances are created by applying the clustering algorithm with $N = \{75,100\}$ on commuters traveling on a selected day (Wednesday of week 2, which had 2,200 commute trips) and setting $K = \{4,5,6,7,8\}$. Let $n$ denote the size of a cluster. Since $N$ only controls the upper bound for the size of clusters produced by the algorithm, residual clusters with $n<N$ are also generated when the total number of commuters, $|\mathcal{C}|$, is not an exact multiple of $N$. For the experiments, only clusters with sizes of exactly 75 and 100 are selected as the main problem instances (24 clusters with $n = 75$ and 22 with $n = 100$) for in-depth study. Detailed results of the REA and BPA on these selected clusters are presented in Tables \ref{tbl:rea_vehicle_capacity_scaling} and \ref{tbl:bpa_vehicle_capacity_scaling} respectively in the appendix. However, Figures \ref{fig:VehCapScaling_VehicleCount}--\ref{fig:VehCapScaling_AvgRideDuration} which aggregates vehicle count, route distance, and average ride duration for the day \emph{utilizes results from all clusters}. 

{\em The REA is unable to complete route enumeration
  within the time limit when $K>6$, therefore results of the algorithm
  for $K = \{7, 8\}$ are not available}. The time limit for clusters
C9-100 and C20-100 when $K=6$ also had to be extended to obtain a
solution. While the BPA is able to handle larger vehicle capacities
for the most part, it could not find a root-node solution within the
time limit for cluster C10-75 when $K = 8$, so the time limit for this
case had to be extended too. As expected, when problems are solved to
optimality, identical objective values are produced by both algorithms
as shown by the same vehicle counts and total route distances in their
results.

\begin{figure}[!t]
	\FIGURE
	{\centering
		\includegraphics[width=0.90\linewidth]{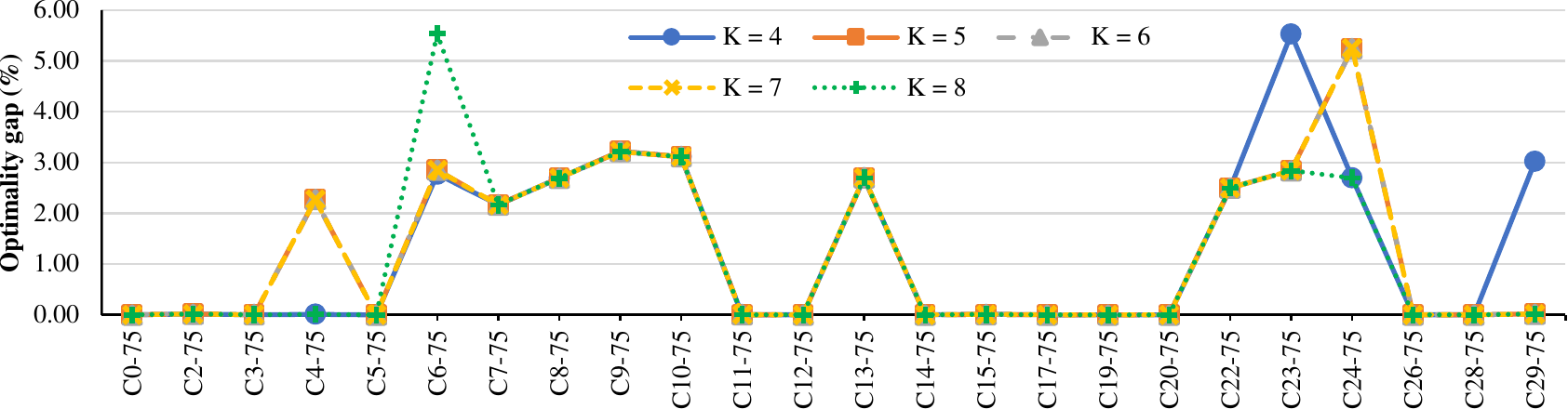}}
	{Optimality Gap of MIP Solution at Root Node of BPA for Problem Instances with $n=75$.\label{fig:RootNodeOptimalityGap_N75}}
	{}
\end{figure}

\begin{figure}[!t]
	\FIGURE
	{\centering
		\includegraphics[width=0.90\linewidth]{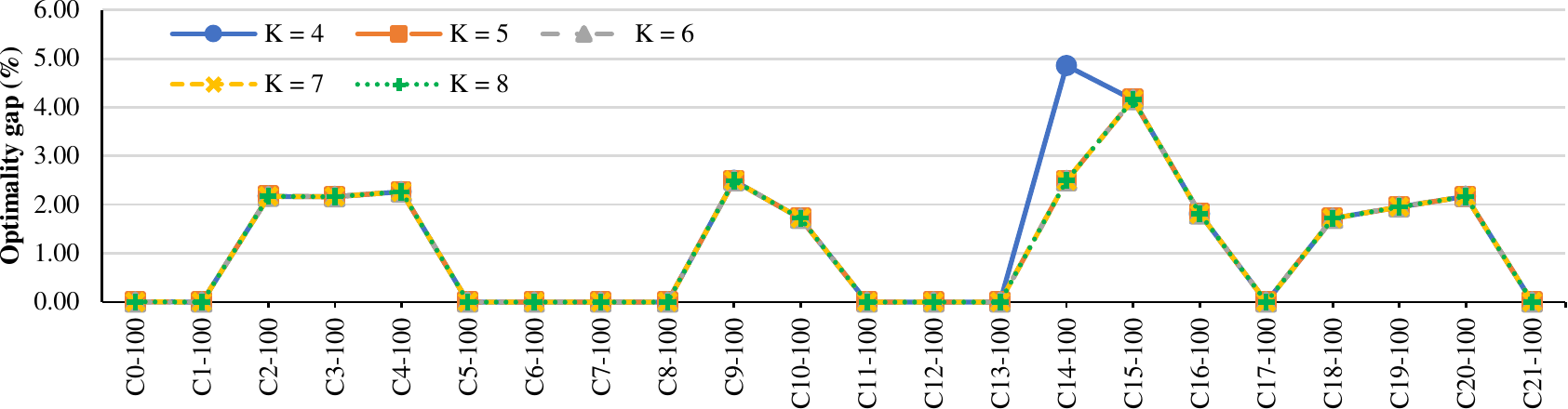}}
	{Optimality Gap of MIP Solution at Root Node of BPA for Problem Instances with $n=100$.\label{fig:RootNodeOptimalityGap_N100}}
	{}
\end{figure}

The REA produces optimal results in all instances when $K\leq6$, while
the BPA does so for all but 14 instances. Unsurprisingly, these 14
instances are typically characterized by large vehicle capacities
($K\geq5$) as well as relatively large edge counts. {\em For these
  instances, the optimality gap of the best feasible solution is
  consistently $<5\%$, and comparison of their vehicle count results
  against those of the REA that are available reveal that they are in
  fact optimal}. Also notable is the number of columns generated by
the BPA being consistently less than the REA, and in some cases
significantly so. Another notable result is the excellent quality of
the BPA's root-node solution which is summarized in Figures
\ref{fig:RootNodeOptimalityGap_N75} and
\ref{fig:RootNodeOptimalityGap_N100} for problem instances with $n=75$
and $n=100$ respectively. Its optimality gap is $<6\%$ in all
instances and is optimal in some cases, making it a viable option when
optimality is not crucial. The integrality gap, also being $<6\%$ for
all instances, emphasizes the strength of the primal lower bound
provided by the RMP's optimal objective value.

Lastly, another notable observation is the disparity in the total
number of feasible inbound and outbound edges of the graphs of the
BPA. Recall that feasible edges are those that satisfy the a priori
feasibility constraints outlined in Section
\ref{sec:edge_elimination}. The edge counts can be seen as a rough
measure of the shareability potential of the set of trips being
considered, and the number of outbound edges being less than inbound
edges in all problem instances indicates fewer sharing opportunities
for outbound trips. This can be attributed to the wider distribution
of their departure times as shown in Figure
\ref{fig:AccessTimeDistribution} which further complicates ride
sharing. It also highlights another unique challenge to solving the
CTSP, as maximal sharing is sought over two sets of trips (inbound and
outbound) with different shareability potential.

\begin{figure}[!t]
	\FIGURE
	{\centering
		\includegraphics[width=0.90\linewidth]{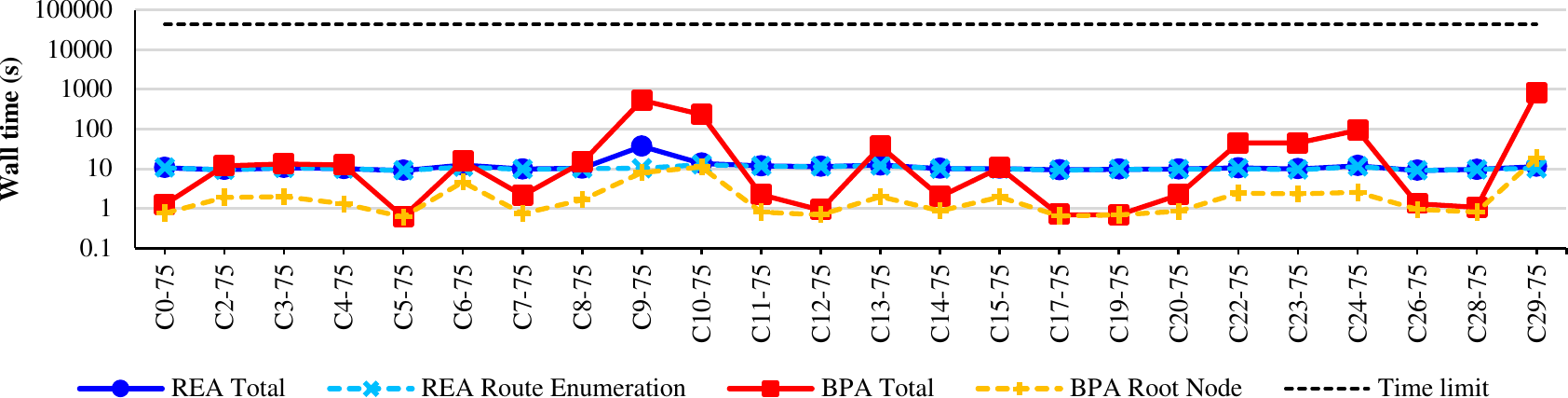}}
	{Computation Times for Problem Instances with $n = 75$ and $K = 4$.\label{fig:CPUTime_N75_VC4}}
	{}
\end{figure}

\begin{figure}[!t]
	\FIGURE
	{\centering
		\includegraphics[width=0.90\linewidth]{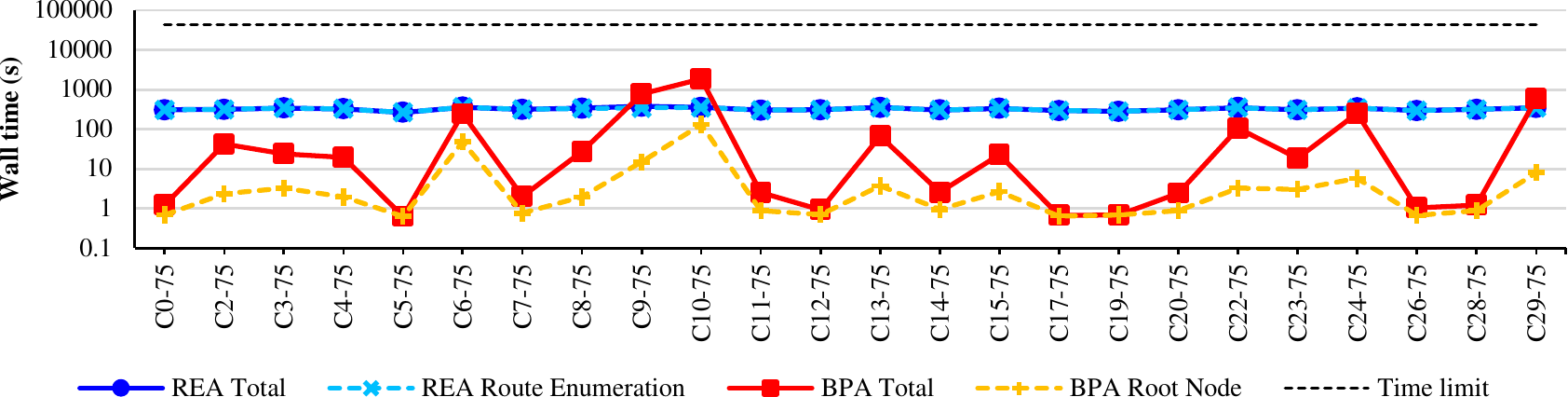}}
	{Computation Times for Problem Instances with $n = 75$ and $K = 5$.\label{fig:CPUTime_N75_VC5}}
	{}
\end{figure}

\begin{figure}[!t]
	\FIGURE
	{\centering
		\includegraphics[width=0.90\linewidth]{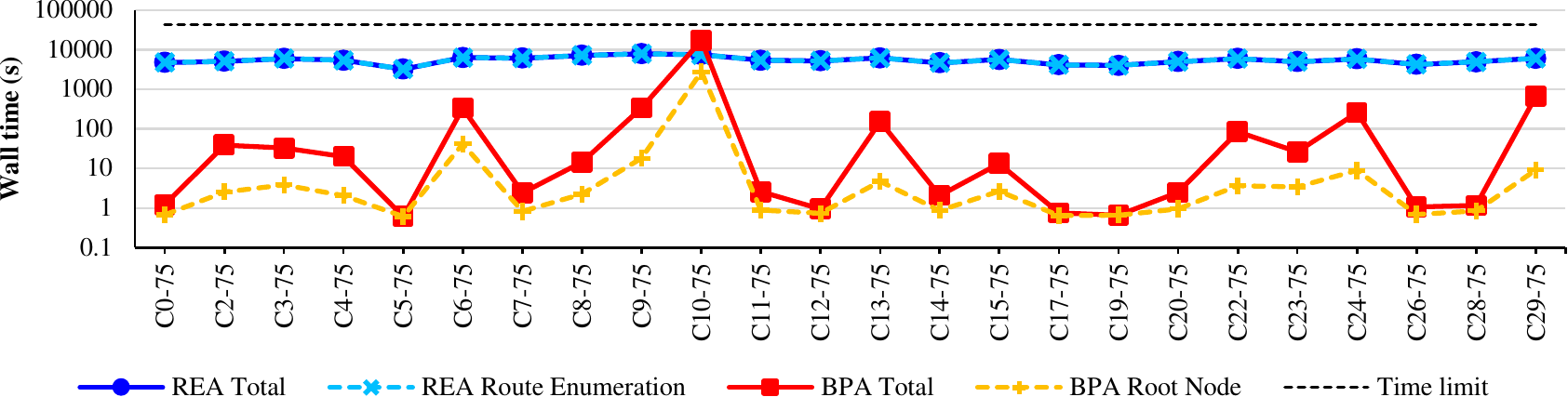}}
	{Computation Times for Problem Instances with $n = 75$ and $K = 6$.\label{fig:CPUTime_N75_VC6}}
	{}
\end{figure}

\begin{figure}[!t]
	\FIGURE
	{\centering
		\includegraphics[width=0.90\linewidth]{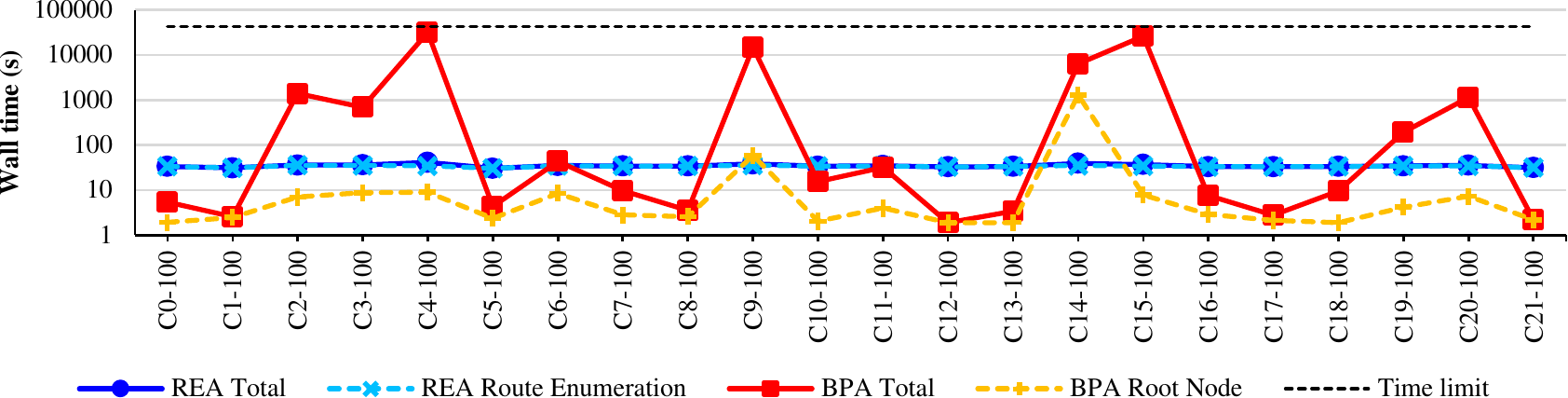}}
	{Computation Times for Problem Instances with $n = 100$ and $K = 4$.\label{fig:CPUTime_N100_VC4}}
	{}
\end{figure}

\begin{figure}[!t]
	\FIGURE
	{\centering
		\includegraphics[width=0.90\linewidth]{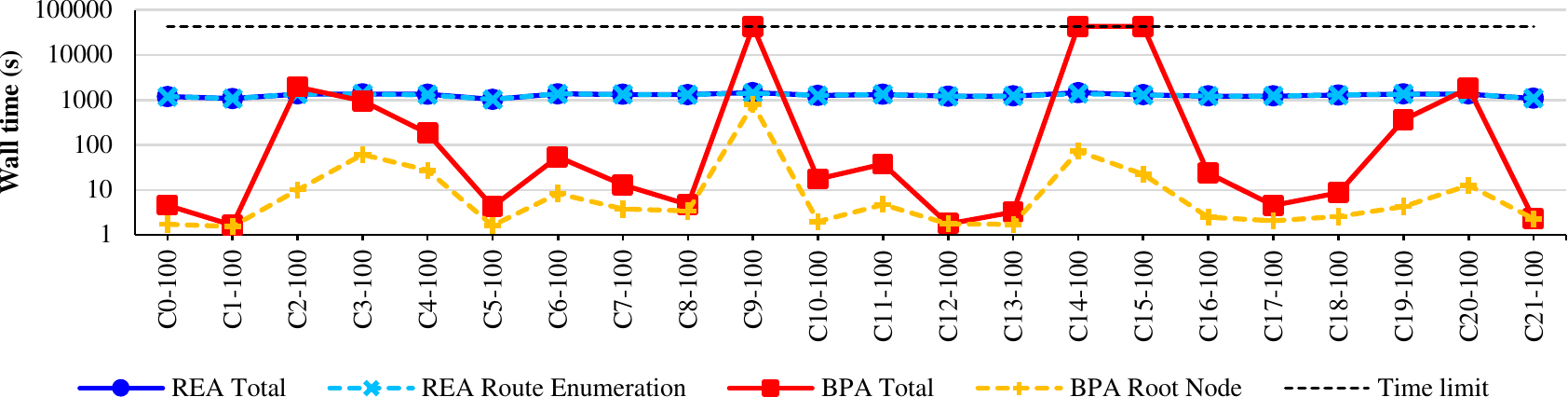}}
	{Computation Times for Problem Instances with $n = 100$ and $K = 5$.\label{fig:CPUTime_N100_VC5}}
	{}
\end{figure}

\begin{figure}[!t]
	\FIGURE
	{\centering
		\includegraphics[width=0.90\linewidth]{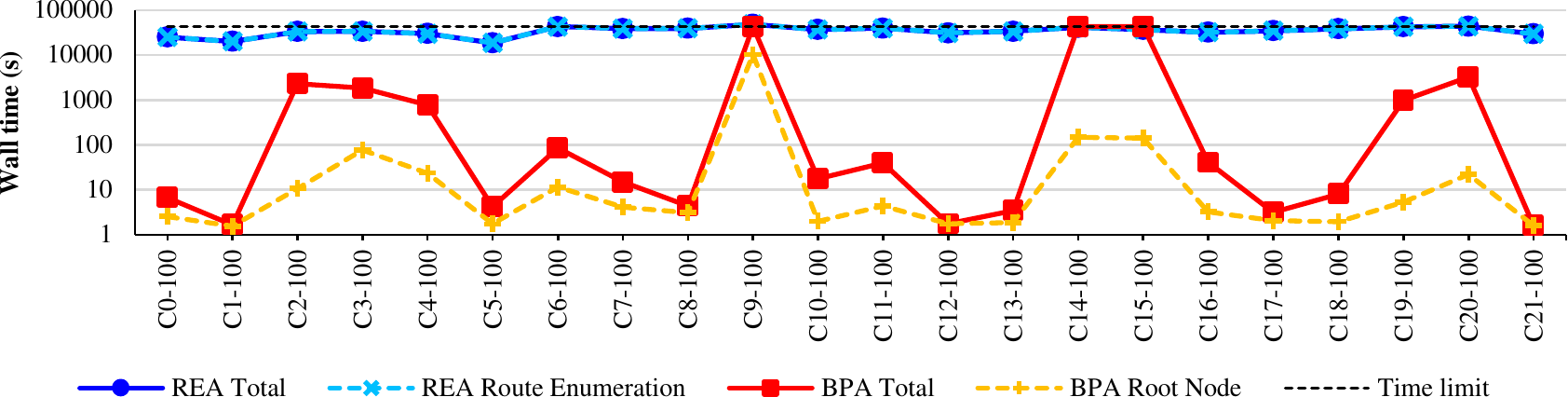}}
	{Computation Times for Problem Instances with $n = 100$ and $K = 6$.\label{fig:CPUTime_N100_VC6}}
	{}
\end{figure}

Figures \ref{fig:CPUTime_N75_VC4}--\ref{fig:CPUTime_N75_VC6} summarize
computation times of both algorithms for the problem instances with $n
= 75$ and $K = \{4,5,6\}$, while Figures
\ref{fig:CPUTime_N100_VC4}--\ref{fig:CPUTime_N100_VC6} do the same for
$n = 100$ and $K = \{4,5,6\}$. The figures reveal that computation
times of the REA are more consistent across problem instances with the
same $n$ and $K$ values. They also appear to be dominated by the route-enumeration phase for these instances. The figures also show that
computation times of the REA are more sensitive to $K$; they appear to
increase more rapidly with increasing $K$ than those of the BPA. The
BPA is slower in 13 out of 24 instances when $K=4$ and $n=75$, and it
is slower in nine out of 22 instances when $K=4$ and $n=100$. These
fractions decrease however as $K$ becomes larger to the point where
the BPA is faster in all but one instance when $K=6$ and $n=75$ and in
all but three instances when $K=6$ and $n=100$. These observations,
combined with results showing the BPA's ability to obtain solutions
when $K\geq6$, indicate that the BPA scales better with increasing
vehicle capacity. Also noteworthy is the time taken to produce the
root-node solution for the BPA; it is faster than the REA in all but
one instance when $n = 75$, and in all but two instances when $n =
100$. This further strengthens the case for it being a viable option
when an optimal solution is not sought.

Finally, Figures \ref{fig:VehCapScaling_VehicleCount},
\ref{fig:VehCapScaling_VehicleMilesTraveled}, and
\ref{fig:VehCapScaling_AvgRideDuration} summarize the effect of
increasing $K$ on total vehicle count, total distance of selected
routes, and average ride time per commuter respectively. They show
aggregated results from all clusters for all trips from a single
weekday (Wednesday of week 2) as $N$ is kept constant. The percentage
of each quantity as a fraction of its value when $K=1$ as well as
results for $K = \{1,2,3\}$ are included for additional
perspective. Firstly, Figures \ref{fig:VehCapScaling_VehicleCount} and
\ref{fig:VehCapScaling_VehicleMilesTraveled} reveal diminishing
marginal decreases in total vehicle count and total travel distance as
$K$ is increased. Furthermore, the benefit of increasing vehicle
capacity almost diminishes completely beyond $K=3$. This can be
attributed to the nature of the routes in the CTSP being very
short. Each needs to start and end at the origin and destination of
its driver respectively, and ride-duration limits are imposed on the
driver in addition to all passengers. The length of each route is
therefore constrained by the ride-duration limit of its driver. As
longer routes are needed to fully utilize the capacities of larger
vehicles (to pick up and drop off more riders), \emph{the routes of this
problem do not benefit from larger vehicle capacities}. For this
reason, subsequent computational experiments only consider the use of cars by limiting $K$ to 4. Figure
\ref{fig:VehCapScaling_AvgRideDuration} provides a first glimpse of
the trade-off in ride sharing: Increased average ride durations. As
vehicle capacity is increased, so does ride-sharing
opportunities. Consequently, an increase in ride duration should also
be expected as a ride is shared with more and more people.  There
appears to be an inverse relationship between average ride duration
and vehicle count or total travel distance, as the former increases as
either of the latter decreases. A similar diminishing effect in the
marginal increase of average ride duration is seen as vehicle capacity
is increased.

\begin{figure}[!t]
	\FIGURE
	{\centering
		\includegraphics[width=0.90\linewidth]{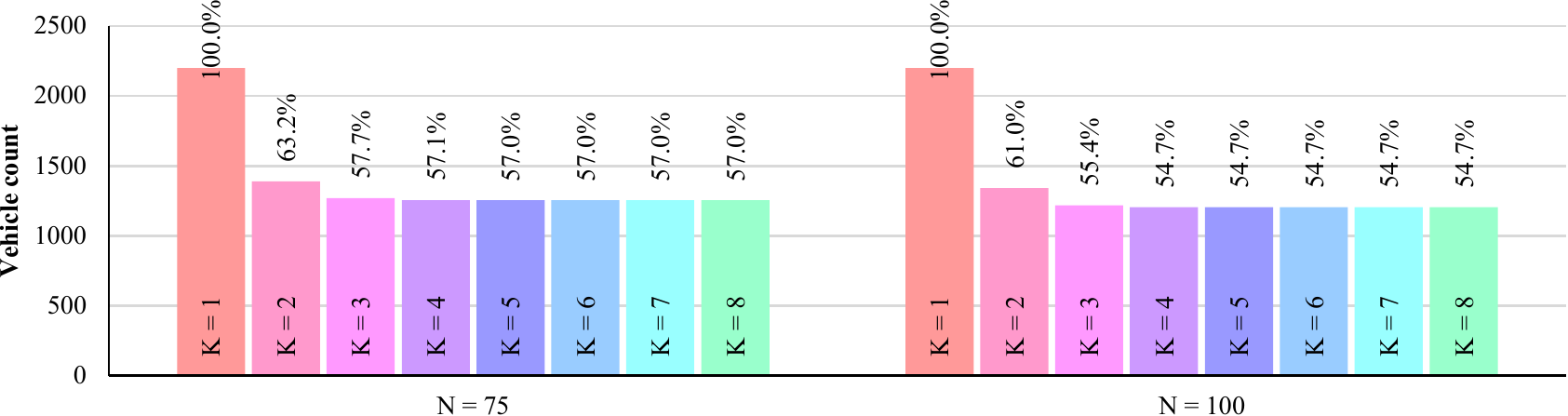}}
	{Effect of Increasing Vehicle Capacity on Total Vehicle Count.\label{fig:VehCapScaling_VehicleCount}}
	{}
\end{figure}

\subsection{Cluster Size Scaling}
\label{sec:clust_size_scaling}

\begin{figure}[!t]
	\FIGURE
	{\centering
		\includegraphics[width=0.90\linewidth]{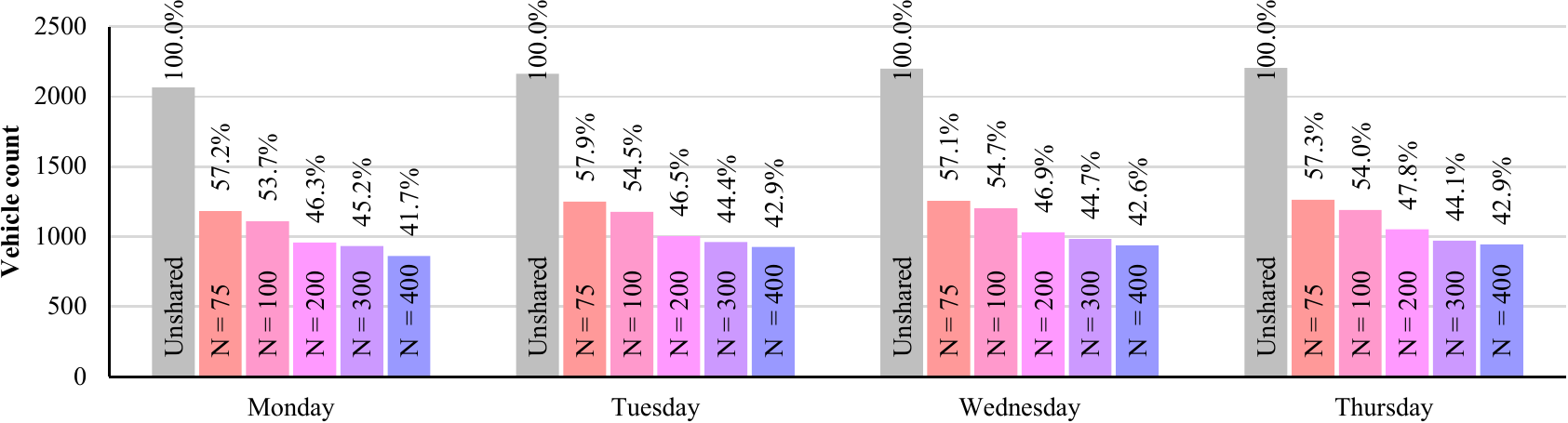}}
	{Effect of Increasing Cluster Size on Total Vehicle Count.\label{fig:ClustSizeScaling_VehicleCount}}
	{}
\end{figure}

\begin{figure}[!t]
	\FIGURE
	{\centering
		\includegraphics[width=0.90\linewidth]{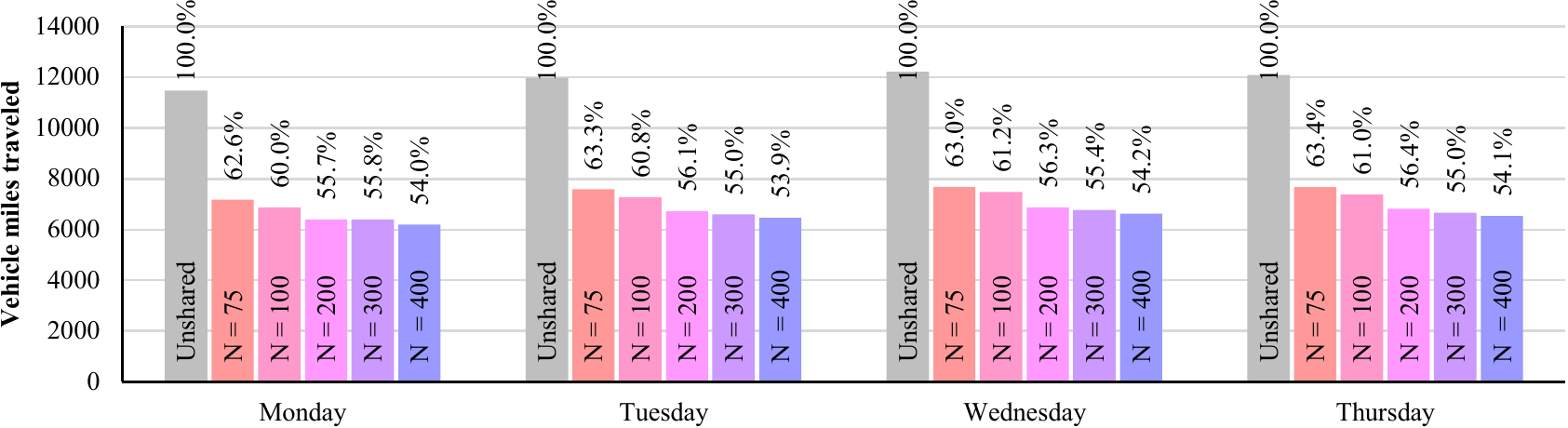}}
	{Effect of Increasing Cluster Size on Total Route Distance.\label{fig:ClustSizeScaling_VehicleMilesTraveled}}
	{}
\end{figure}

\begin{figure}[!t]
	\FIGURE
	{\centering
		\includegraphics[width=0.90\linewidth]{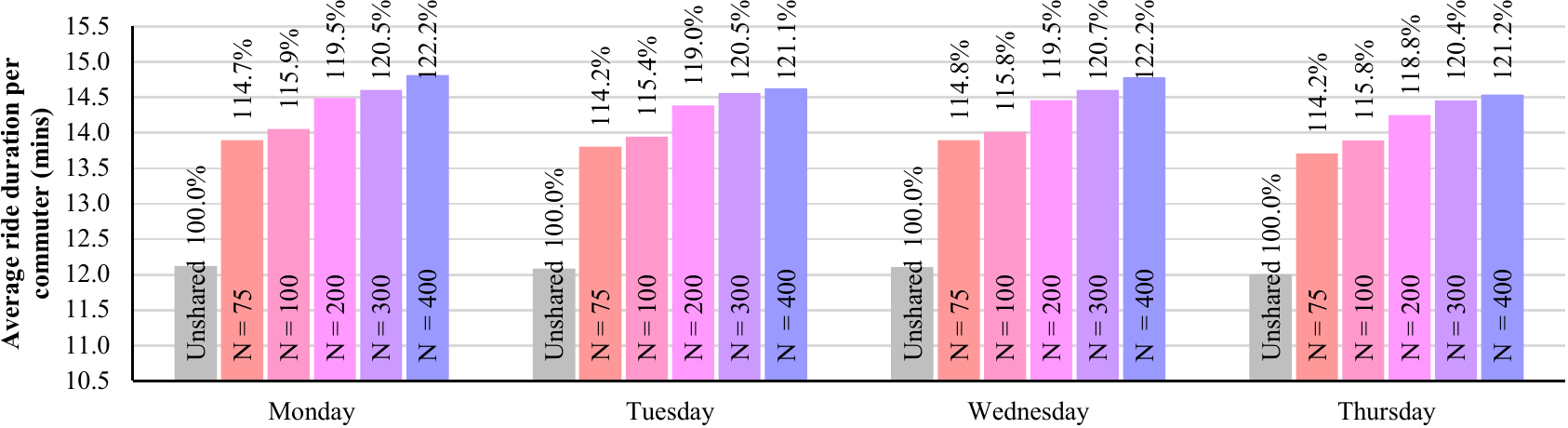}}
	{Effect of Increasing Cluster Size on Average Ride Duration.\label{fig:ClustSizeScaling_AvgRideDuration}}
	{}
\end{figure}

The next set of experiments explores the scalability of both algorithms with increasing cluster size. To this end, the clustering algorithm is applied to commuters traveling on each of the first four weekdays of week 2 (which had 2065, 2161, 2200, and 2203 commute trips respectively) with $N$ set to $\{200,300,400\}$ (results for $N = \{75,100\}$ are already available from Section \ref{sec:veh_cap_scaling})). $K$ is fixed to 4 in all experiments. Detailed results of the REA and BPA, shown in Tables \ref{tbl:rea_cluster_size_scaling} and \ref{tbl:bpa_cluster_size_scaling} respectively in the appendix, and the rest of the discussions in this section focus on selected clusters with sizes of exactly 200, 300, and 400 (11 with $n = 200$, 13 with $n = 300$, and 8 with $n = 400$) as the results are intended to show the effect of progressively increasing the cluster size by increments of 100. However, results from all clusters, including residuals with $n<N$, are used in Figures \ref{fig:ClustSizeScaling_VehicleCount}, \ref{fig:ClustSizeScaling_VehicleMilesTraveled}, and \ref{fig:ClustSizeScaling_AvgRideDuration} which show aggregate vehicle count, route distance, and average ride duration respectively for each day.

\begin{figure}[!t]
	\FIGURE
	{\centering
		\includegraphics[width=0.45\linewidth]{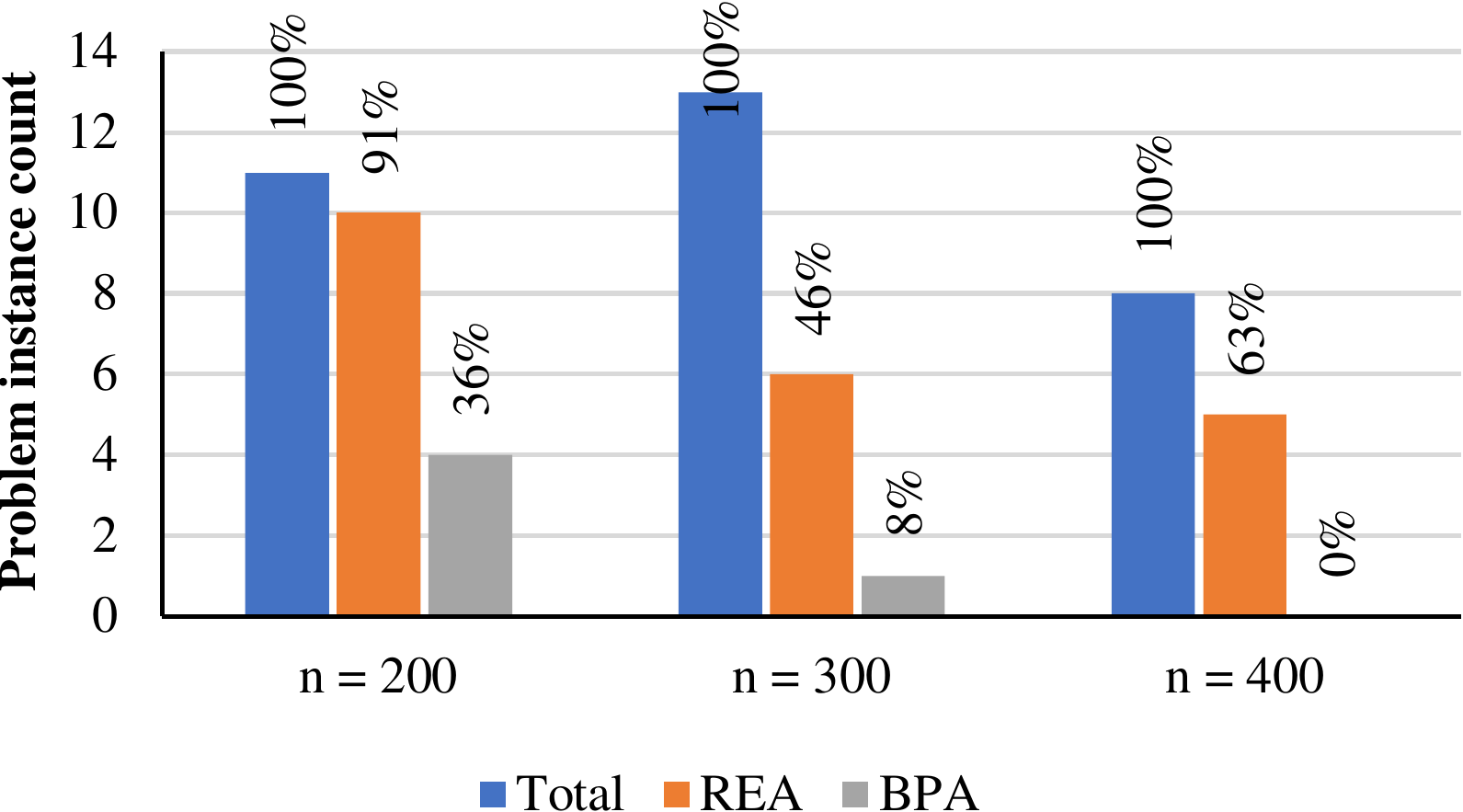}}
	{Number of Problem Instances Solved Optimally when $n=\{200,300,400\}$.\label{fig:SummaryOptimalCount}}
	{}
\end{figure}

Figure \ref{fig:SummaryOptimalCount} summarizes the number of instances that can be solved optimally by both algorithms together with the total number of instances for each $n$ value and the percentage of each quantity as a fraction of the total. When $n = 200$, the REA is able to obtain the optimal solution for all but one instance. Conversely, the BPA is only able to produce optimality for four instances. As $n$ is increased, so does the size of the problem instances as evident from the number of columns generated and edge counts, which leads to fewer instances being solved to optimality by either algorithm. When $n = 400$, none of the instances could be solved to optimality by the BPA, and only five out of eight instances could be solved optimally by the REA.

\begin{figure}[!t]
	\FIGURE
	{\centering
		\includegraphics[width=0.90\linewidth]{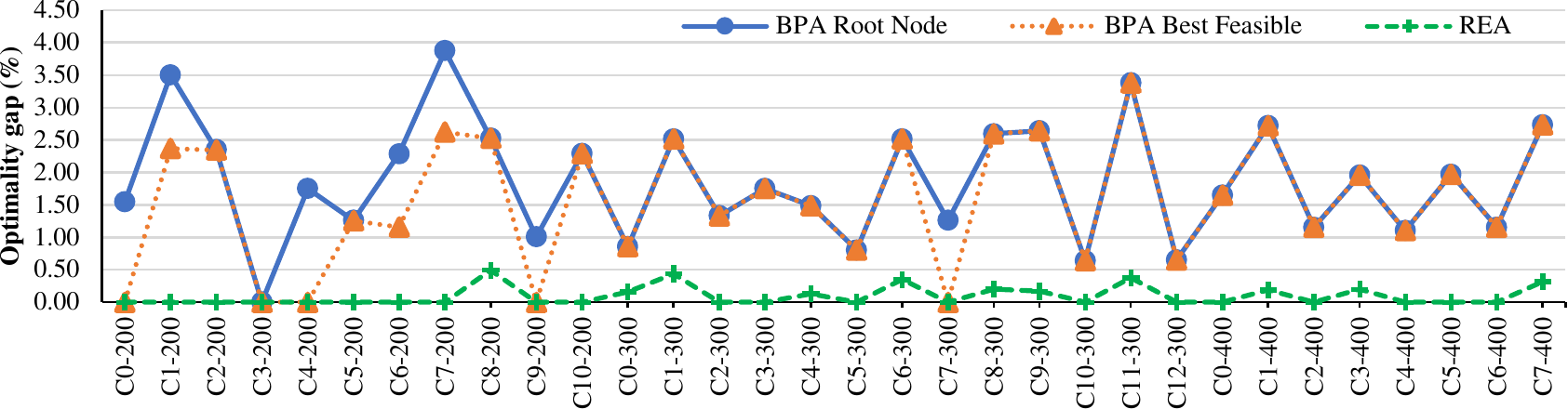}}
	{Optimality Gaps of the REA and the BPA for Problem Instances with $n=\{200,300,400\}$.\label{fig:OptimalityGap_N200_300_400}}
	{}
\end{figure}

Figure \ref{fig:OptimalityGap_N200_300_400} displays the optimality gaps produced by both algorithms in these experiments. The optimality gaps of suboptimal solutions of the REA are excellent, being $<0.5\%$ for all instances. Moreover, the optimal vehicle count is obtained in all these instances. The optimality gaps of suboptimal BPA solutions are competitive, being $<4\%$ in all instances, however there are a few instances where the optimal vehicle count is not obtained, e.g. clusters C8-200, C3-300, and C0-400. For these instances, the vehicle count is typically off by one when compared to the optimal counts of the REA. Root-node solutions of the BPA remain excellent as their optimality gaps are also $<4\%$ in all of the instances tested, with it being optimal for one instance (cluster C3-200). The total number of columns generated by the BPA also remain fewer than the REA in all instances. 

\begin{figure}[!t]
	\FIGURE
	{\centering
		\includegraphics[width=0.90\linewidth]{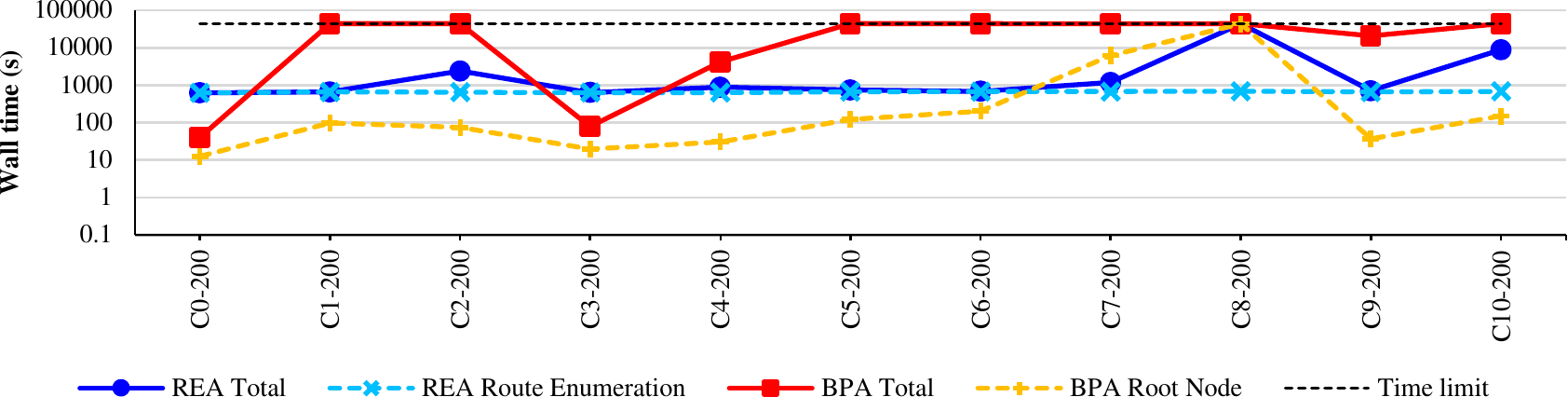}}
	{Computation Times for Problem Instances with $n = 200$.\label{fig:CPUTime_N200_VC4}}
	{}
\end{figure}

\begin{figure}[!t]
	\FIGURE
	{\centering
		\includegraphics[width=0.90\linewidth]{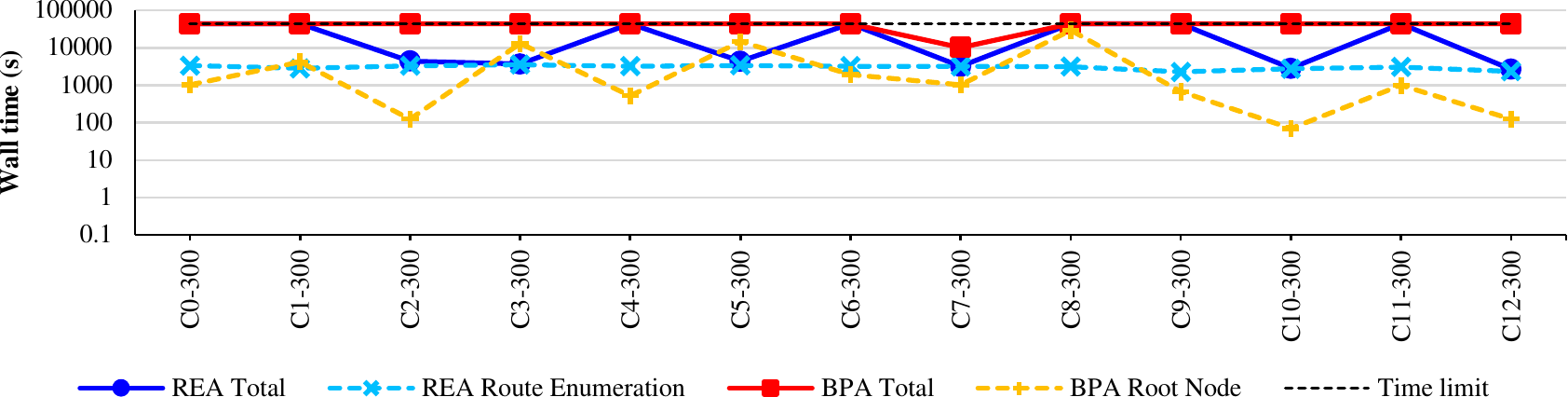}}
	{Computation Times for Problem Instances with $n = 300$.\label{fig:CPUTime_N300_VC4}}
	{}
\end{figure}

\begin{figure}[!t]
	\FIGURE
	{\centering
		\includegraphics[width=0.90\linewidth]{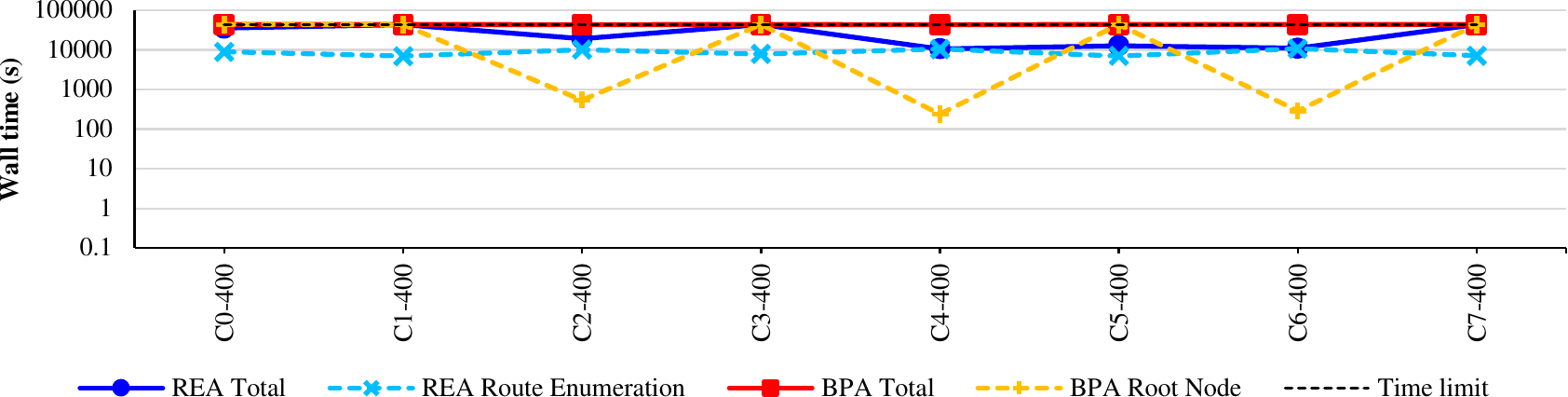}}
	{Computation Times for Problem Instances with $n = 400$.\label{fig:CPUTime_N400_VC4}}
	{}
\end{figure}

Figures \ref{fig:CPUTime_N200_VC4}--\ref{fig:CPUTime_N400_VC4} compare computation times of both algorithms for all problem instances with $n = \{200,300,400\}$. When combined with Figures \ref{fig:CPUTime_N75_VC4} and \ref{fig:CPUTime_N100_VC4}, they provide a clear picture of how both algorithms scale with increasing cluster size as $K$ is kept constant at 4. Figure \ref{fig:CPUTime_N200_VC4} shows the REA being faster then the BPA in eight out of the 11 instances tested when $n = 200$. The allocated time quickly gets saturated in most problem instances by either algorithm as $n$ is increased to the the point where the BPA reaches the time limit in all problem instances when $n = 400$, while the REA also does so for three out of the eight instances tested. The computation times of the REA are still dominated by its route-enumeration phase in most instances, however there also exist a few instances where solving the MP consumes a bigger portion of the total computation times, e.g. cluster C8-200 in Figure \ref{fig:CPUTime_N200_VC4}. This is due to the increased complexity of solving the MIP from the higher column counts. More fluctuations are also observed in the computation times of the REA across problems with the same $n$ value, and this too can be attributed to the complexity of solving the MIP as the time for the route-enumeration phase remains consistent for all instances. The root-node solution of the BPA remains a viable option if a quick, high-quality solution is required, as it is faster than the REA in all but nine instances. The results indicate that the REA is a better choice for large cluster sizes if an optimal solution, or one that is as close to optimal as possible, is sought. However, the root-node solution of the BPA is the best option if only a high-quality solution is needed, as it is faster than the REA in most of the instances tested.

Finally, Figures \ref{fig:ClustSizeScaling_VehicleCount}, \ref{fig:ClustSizeScaling_VehicleMilesTraveled}, and \ref{fig:ClustSizeScaling_AvgRideDuration} show the effect of increasing $N$ on the overall results in terms of total vehicle count, total distance of selected routes, and the average ride time per commuter. The figures show aggregated results from all clusters for the first four weekdays of week 2. Results for unshared trips as well as the value of each quantity as a percentage of their corresponding unshared values are included for added context. As expected, total vehicle count and total route distance results improve as $N$ is increased, as larger $N$ values produce larger `neighborhoods' which provide increased ride-sharing opportunities. Average ride duration also increases with $N$ as expected. Diminishing marginal decreases (respectively increases) in total vehicle count and total route distance (respectively average ride duration) are also observed with increasing $N$, however the effect is not as pronounced as that from increasing $K$, signaling that increasing maximum cluster size is more effective in improving ride-sharing results than increasing vehicle capacity for the CTSP. 

\subsection{Efficiency of the Root-Node Heuristic}
Finally, the root-node heuristic is applied on the selected clusters with $n=\{100,200,300,400\}$ from Sections \ref{sec:veh_cap_scaling} and \ref{sec:clust_size_scaling} with time budgets of $t_{\text{RMP}} = 8$ minutes and $t_{\text{MIP}} = 2$ minutes, resulting in a total time budget of 10 minutes per instance which is deemed reasonable for obtaining solutions in an operational setting. Detailed results are presented in Table \ref{tbl:root_node_heuristic} in the appendix. 

\begin{figure}[!t]
	\FIGURE
	{\centering
		\includegraphics[width=0.9\linewidth]{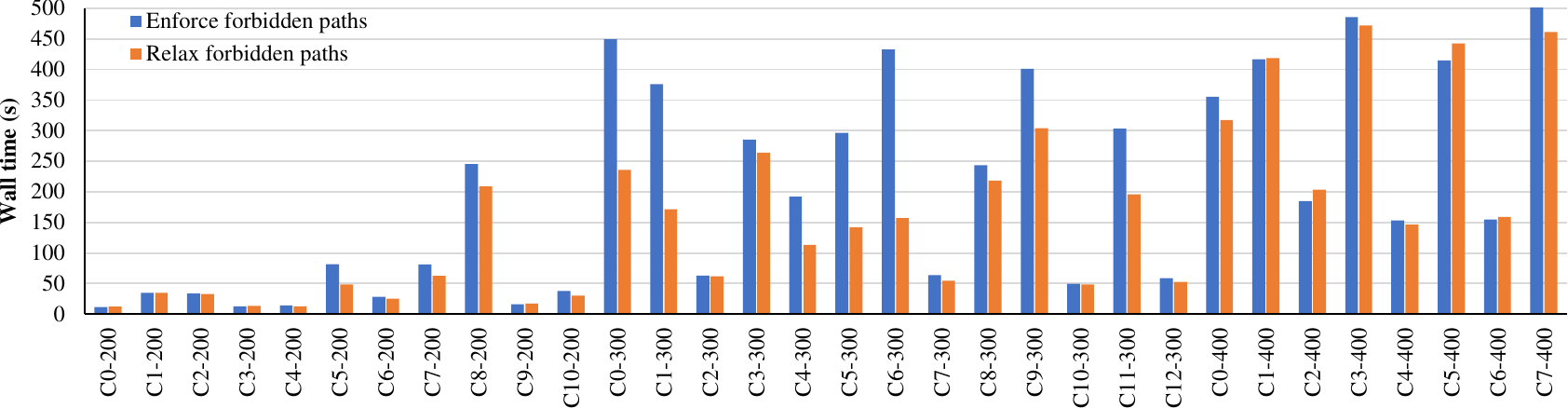}}
	{RMP Convergence Times of Root-Node Heuristics for Problem Instances with $n=\{200,300,400\}$\label{fig:RootNodeHeuristicWallTime}}
	{}
\end{figure}

Figure \ref{fig:RootNodeHeuristicOptimalityGap} compares the
optimality gaps of both heuristics for problem instances with $n \geq
200$, while Figure \ref{fig:RootNodeHeuristicWallTime} compares the
time spent for their RMPs to converge. It can be seen from the first
figure that the optimality gaps of both variants are typically $<5\%$,
with the relaxation heuristic having an optimality gap that is only
$0.1\%$ larger on average. This minimal loss can be attributed to its
small fraction of infeasible routes ($<0.5\%$ of all routes in the
instances tested). In some instances, the relaxation heuristic
produces optimality gaps that are smaller, and this can be attributed
to it being able to solve significantly more column-generation
iterations within its time budget compared to the other heuristic
which has to execute the expensive forbidden path algorithm. The
second figure shows the relaxation heuristic completing its
column-generation phase faster in the majority of the problem
instances, with it being on average $26\%$ faster than the first
heuristic. Nevertheless, regardless of the variant used, the results
further reinforce initial claims that the root-node heuristic is
indeed able to generate provably high-quality solutions in
time-constrained scenarios for problem instances of various sizes.

\begin{figure}[!t]
	\FIGURE
	{\centering
		\includegraphics[width=0.9\linewidth]{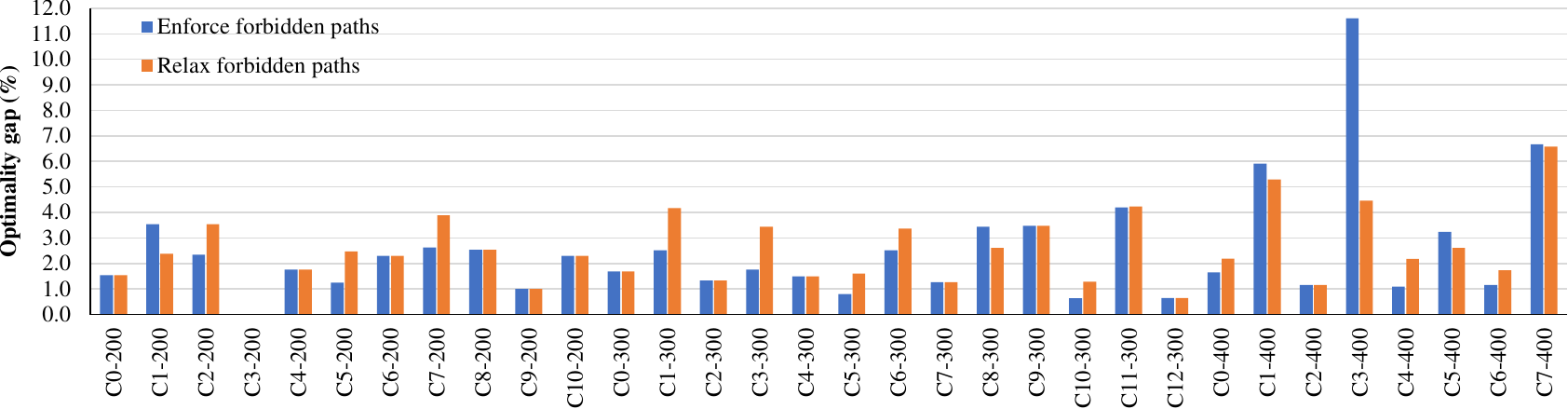}}
	{Optimality Gaps of Root-Node Heuristics for Problem Instances with $n=\{200,300,400\}$.\label{fig:RootNodeHeuristicOptimalityGap}}
	{}
\end{figure}

\section{Conclusion}
\label{sec:conclusion}

To relieve parking pressure which has been steadily increasing in
cities and in university and corporate campuses, this paper explored a
car-pooling platform that would match riders and drivers, while
guaranteeing a ride back and exploiting spatial and temporal
locality. It formalized the Commute Trip Sharing Problem (CTSP) to
find a routing plan that maximizes ride sharing for a set of commute
trips and proposed two exact algorithms for the CTSP: A
Route-Enumeration Algorithm (REA) and a Branch-and-Price
Algorithm (BPA). The former exhaustively searches for feasible routes from
all possible trip combinations, which are then supplied to a MIP which
solves a set-partitioning problem. The latter uses column generation
which applies a dynamic-programming algorithm to search for feasible
routes with negative marginal costs on demand. A clustering algorithm
is also proposed to group trips based on commuter home locations to
maintain problem tractability. The REA and the BPA are then used to
optimally match commute trips from a real-word dataset for the city of
Ann Arbor, Michigan.
 
Results of the computational experiments revealed that the BPA is
better suited for problems with larger vehicle capacities, although
they also revealed that there is very little benefit to utilizing
vehicles with capacity greater than 4 in the CTSP as their
effectiveness is mitigated by ride-duration constraints which
restrict route length. On the other hand, the REA is found to be
better suited for problems with large commuter counts, as it
consistently produced results that are optimal or have optimality gaps
that are often smaller than the BPA. The root-node solution of the BPA
is also found to be a good heuristic for producing high-quality
solutions in time-constrained scenarios, as it typically produces
solutions with optimality gaps of $<5$\%, even for larger problem
instances, within a 10-minute time span.
 
When it is assumed that commuters are willing to shift their desired
arrival and departure times by $\pm10$ minutes and tolerate a 50\%
increase to their ride durations, the algorithms can produce optimal
solutions for problems with up to 200 commuters and achieve
high-quality results for problems with up to 400 commuters. When the
maximum cluster size is set to 400, the results show the CTSP plans
potentially reducing vehicle utilization by 57\% and decreasing
vehicle miles traveled by 46\% at the cost of a 22\% increase in
average ride duration. The results thus highlight the significant
potential and effectiveness of the CTSP in easing traffic and parking
pressure on otherwise congested areas. Future work will be dedicated
to further increasing the efficiency of the algorithms while making
them more robust to changes in commuter schedules and additions of new
customers to the commuting pool. Similarly, behavioral studies to
determine adoption of such a car-pooling will be performed to
understand how much of the theoretical potential can be achieved in
various practical settings.

\ACKNOWLEDGMENT{%

We would like to thank Steve Dolen for his help with the case
study. This work was partially supported by a grant from the Michigan
Institute of Data Science.
  
}

%
\newpage
\begin{APPENDIX}{Computational Results}

Results of the vehicle capacity scaling experiments for the REA are presented in Table
  \ref{tbl:rea_vehicle_capacity_scaling}. The first three columns show
  the cluster size, vehicle capacity, and cluster ID which
  characterize each problem instance. The next column lists the number
  of columns generated by the algorithm, while the following three
  show the results of solving the MP with a MIP solver. They show the vehicle
  count, the total distance of selected routes, and the optimality gap
  of the MIP solution. Finally, the remaining two columns display
computation times of the route-enumeration phase and of the entire
  algorithm including MIP solve times.

Table \ref{tbl:bpa_vehicle_capacity_scaling} shows results of the same set of experiments 
for the BPA. Similar to Table \ref{tbl:rea_vehicle_capacity_scaling}, the
first three columns list the cluster size, vehicle capacity, and
cluster ID for each problem instance. The next two columns present the
total number of unique feasible edges from all inbound and outbound
graphs to further characterize the size of the problem instances,
while the following two show the total number of tree nodes explored
and columns generated by the algorithm. The next two display the
results in terms of the vehicle count and the total distance of
selected routes. The following two columns list optimality gaps of the
MIP solution at the root node and of the best feasible solution, while
the next lists the integrality gap of the best feasible solution,
which is the relative gap between its objective value and
$z^*$. Finally, the remaining four columns show computation times for
the RMP to converge, for finding the MIP solution at the root node,
for arriving at the best feasible solution, and for executing the
entire algorithm.

Tables \ref{tbl:rea_cluster_size_scaling} and \ref{tbl:bpa_cluster_size_scaling} summarizes results of the cluster size scaling experiments for the REA and the BPA respectively. They list the same quantities as those listed in Tables \ref{tbl:rea_vehicle_capacity_scaling} and \ref{tbl:bpa_vehicle_capacity_scaling} respectively.

Finally, Table \ref{tbl:root_node_heuristic} shows results of the experiments which investigate the efficiency of the root-node heuristic. The first two columns show the cluster size and ID of each problem instance. The next set of five columns list the results of the heuristic which enforces forbidden paths. They show the number of columns generated, the resulting vehicle count, its optimality gap, and the times spent on solving the RMP and its MIP. The next set of six columns display the same information for the heuristic which relaxes forbidden paths, with an additional column showing the number of infeasible columns it generated.

\begin{table}[!htbp]
	\TABLE
	{Results of REA Scalability with Increasing Vehicle Capacity ($\Delta = 10$ mins, $R = 50\%$).\label{tbl:rea_vehicle_capacity_scaling}}
	{\resizebox*{!}{\dimexpr\textheight-5\baselineskip\relax}{\begin{tabular}{ccrrrrrrr}
		\hline \up\down
		\multirow{2}{*}{\begin{tabular}[c]{@{}c@{}}Cluster\\ size\end{tabular}} & \multirow{2}{*}{\begin{tabular}[c]{@{}c@{}}Vehicle\\ capacity\end{tabular}} & \multicolumn{1}{c}{\multirow{2}{*}{\begin{tabular}[c]{@{}c@{}}Cluster\\ ID\end{tabular}}} & \multicolumn{1}{c}{\multirow{2}{*}{\begin{tabular}[c]{@{}c@{}}Column\\ \#\end{tabular}}} & \multicolumn{1}{c}{\multirow{2}{*}{\begin{tabular}[c]{@{}c@{}}Vehicle\\ \#\end{tabular}}} & \multicolumn{1}{c}{\multirow{2}{*}{\begin{tabular}[c]{@{}c@{}}Total\\ distance\\ (m)\end{tabular}}} & \multicolumn{1}{c}{\multirow{2}{*}{\begin{tabular}[c]{@{}c@{}}Optimality\\ gap (\%)\end{tabular}}} & \multicolumn{2}{c}{Wall time (s)} \\ \cline{8-9} 
		&  & \multicolumn{1}{c}{} & \multicolumn{1}{c}{} & \multicolumn{1}{c}{} & \multicolumn{1}{c}{} & \multicolumn{1}{c}{} & \multicolumn{1}{c}{\begin{tabular}[c]{@{}c@{}}Route\\ enumeration\end{tabular}} & \multicolumn{1}{c}{Total} \\ \hline \up
		\multirow{72}{*}{75} & \multirow{24}{*}{4} & C0-75 & 508 & 47 & 356296 & 0.00 & 11 & 11 \\  
		&  & C2-75 & 1946 & 40 & 643541 & 0.00 & 9 & 9 \\  
		&  & C3-75 & 2068 & 38 & 487061 & 0.00 & 11 & 11 \\  
		&  & C4-75 & 1483 & 43 & 398265 & 0.00 & 10 & 10 \\  
		&  & C5-75 & 274 & 61 & 230460 & 0.00 & 9 & 9 \\  
		&  & C6-75 & 3071 & 36 & 451262 & 0.00 & 11 & 12 \\  
		&  & C7-75 & 690 & 46 & 437743 & 0.00 & 10 & 10 \\  
		&  & C8-75 & 1350 & 37 & 487561 & 0.00 & 10 & 10 \\  
		&  & C9-75 & 3926 & 31 & 328525 & 0.00 & 10 & 37 \\  
		&  & C10-75 & 4137 & 32 & 527544 & 0.00 & 13 & 14 \\  
		&  & C11-75 & 866 & 47 & 443770 & 0.00 & 12 & 12 \\  
		&  & C12-75 & 592 & 46 & 298162 & 0.00 & 11 & 11 \\  
		&  & C13-75 & 2143 & 37 & 481869 & 0.00 & 12 & 12 \\  
		&  & C14-75 & 863 & 47 & 301108 & 0.00 & 10 & 10 \\  
		&  & C15-75 & 1512 & 36 & 460087 & 0.00 & 10 & 10 \\  
		&  & C17-75 & 457 & 50 & 348581 & 0.00 & 9 & 9 \\  
		&  & C19-75 & 496 & 57 & 245827 & 0.00 & 10 & 10 \\  
		&  & C20-75 & 752 & 46 & 420825 & 0.00 & 10 & 10 \\  
		&  & C22-75 & 2288 & 40 & 472247 & 0.00 & 10 & 11 \\  
		&  & C23-75 & 1574 & 36 & 385434 & 0.00 & 9 & 10 \\  
		&  & C24-75 & 1926 & 37 & 425807 & 0.00 & 11 & 12 \\  
		&  & C26-75 & 468 & 52 & 333692 & 0.00 & 9 & 9 \\  
		&  & C28-75 & 750 & 44 & 293838 & 0.00 & 10 & 10 \\  
		&  & C29-75 & 2541 & 32 & 548779 & 0.00 & 10 & 11 \\ \cline{2-9} 
		& \multirow{24}{*}{5} & C0-75 & 509 & 47 & 356296 & 0.00 & 306 & 306 \\  
		&  & C2-75 & 2014 & 40 & 643540 & 0.00 & 316 & 316 \\  
		&  & C3-75 & 2117 & 38 & 487061 & 0.00 & 341 & 342 \\  
		&  & C4-75 & 1559 & 43 & 398216 & 0.00 & 326 & 326 \\  
		&  & C5-75 & 274 & 61 & 230460 & 0.00 & 264 & 264 \\  
		&  & C6-75 & 3331 & 35 & 450414 & 0.00 & 355 & 356 \\  
		&  & C7-75 & 690 & 46 & 437743 & 0.00 & 316 & 316 \\  
		&  & C8-75 & 1352 & 37 & 487245 & 0.00 & 335 & 335 \\  
		&  & C9-75 & 4110 & 31 & 327993 & 0.00 & 346 & 378 \\  
		&  & C10-75 & 5589 & 32 & 526188 & 0.00 & 352 & 354 \\  
		&  & C11-75 & 870 & 47 & 443770 & 0.00 & 300 & 300 \\  
		&  & C12-75 & 592 & 46 & 298162 & 0.00 & 308 & 308 \\  
		&  & C13-75 & 2207 & 37 & 479166 & 0.00 & 351 & 351 \\  
		&  & C14-75 & 867 & 47 & 301108 & 0.00 & 304 & 304 \\  
		&  & C15-75 & 1522 & 36 & 460087 & 0.00 & 335 & 335 \\  
		&  & C17-75 & 457 & 50 & 348581 & 0.00 & 293 & 293 \\  
		&  & C19-75 & 496 & 57 & 245827 & 0.00 & 279 & 279 \\  
		&  & C20-75 & 752 & 46 & 420825 & 0.00 & 309 & 309 \\  
		&  & C22-75 & 2344 & 40 & 472040 & 0.00 & 347 & 348 \\  
		&  & C23-75 & 1628 & 35 & 379234 & 0.00 & 306 & 307 \\  
		&  & C24-75 & 1996 & 37 & 425806 & 0.00 & 339 & 339 \\  
		&  & C26-75 & 468 & 52 & 333692 & 0.00 & 292 & 292 \\  
		&  & C28-75 & 751 & 44 & 293838 & 0.00 & 315 & 315 \\  
		&  & C29-75 & 2594 & 32 & 543850 & 0.00 & 350 & 351 \\ \cline{2-9}
		& \multirow{24}{*}{6} & C0-75 & 509 & 47 & 356296 & 0.00 & 4761 & 4761 \\  
		&  & C2-75 & 2016 & 40 & 643540 & 0.00 & 5109 & 5110 \\  
		&  & C3-75 & 2120 & 38 & 487061 & 0.00 & 6009 & 6009 \\  
		&  & C4-75 & 1562 & 43 & 398216 & 0.00 & 5375 & 5375 \\  
		&  & C5-75 & 274 & 61 & 230460 & 0.00 & 3220 & 3220 \\  
		&  & C6-75 & 3406 & 35 & 450414 & 0.00 & 6341 & 6341 \\  
		&  & C7-75 & 690 & 46 & 437743 & 0.00 & 6135 & 6135 \\  
		&  & C8-75 & 1352 & 37 & 487245 & 0.00 & 7172 & 7172 \\  
		&  & C9-75 & 4116 & 31 & 327993 & 0.00 & 7893 & 7916 \\  
		&  & C10-75 & 6571 & 32 & 525115 & 0.00 & 7480 & 7482 \\  
		&  & C11-75 & 870 & 47 & 443770 & 0.00 & 5413 & 5413 \\  
		&  & C12-75 & 592 & 46 & 298162 & 0.00 & 5206 & 5206 \\  
		&  & C13-75 & 2209 & 37 & 479166 & 0.00 & 6164 & 6164 \\  
		&  & C14-75 & 867 & 47 & 301108 & 0.00 & 4658 & 4658 \\  
		&  & C15-75 & 1522 & 36 & 460087 & 0.00 & 5668 & 5668 \\  
		&  & C17-75 & 457 & 50 & 348581 & 0.00 & 4147 & 4147 \\  
		&  & C19-75 & 496 & 57 & 245827 & 0.00 & 4000 & 4000 \\  
		&  & C20-75 & 752 & 46 & 420825 & 0.00 & 4983 & 4983 \\  
		&  & C22-75 & 2345 & 40 & 472040 & 0.00 & 5961 & 5961 \\  
		&  & C23-75 & 1629 & 35 & 379234 & 0.00 & 5008 & 5008 \\  
		&  & C24-75 & 2001 & 37 & 425806 & 0.00 & 5843 & 5844 \\  
		&  & C26-75 & 468 & 52 & 333692 & 0.00 & 4254 & 4254 \\  
		&  & C28-75 & 751 & 44 & 293838 & 0.00 & 4935 & 4935 \\ \down 
		&  & C29-75 & 2594 & 32 & 543850 & 0.00 & 6071 & 6071 \\ \hline
	\end{tabular}}}
{}
\end{table}

\begin{table}[]
	\centering
	\resizebox*{!}{\textheight}
	{\begin{tabular}{ccrrrrrrr}
		\hline
		\up\down \multirow{2}{*}{\begin{tabular}[c]{@{}c@{}}Cluster\\ size\end{tabular}} & \multirow{2}{*}{\begin{tabular}[c]{@{}c@{}}Vehicle\\ capacity\end{tabular}} & \multicolumn{1}{c}{\multirow{2}{*}{\begin{tabular}[c]{@{}c@{}}Cluster\\ ID\end{tabular}}} & \multicolumn{1}{c}{\multirow{2}{*}{\begin{tabular}[c]{@{}c@{}}Column\\ \#\end{tabular}}} & \multicolumn{1}{c}{\multirow{2}{*}{\begin{tabular}[c]{@{}c@{}}Vehicle\\ \#\end{tabular}}} & \multicolumn{1}{c}{\multirow{2}{*}{\begin{tabular}[c]{@{}c@{}}Total\\ distance\\ (m)\end{tabular}}} & \multicolumn{1}{c}{\multirow{2}{*}{\begin{tabular}[c]{@{}c@{}}Optimality\\ gap (\%)\end{tabular}}} & \multicolumn{2}{c}{Wall time (s)} \\ \cline{8-9} 
		&  & \multicolumn{1}{c}{} & \multicolumn{1}{c}{} & \multicolumn{1}{c}{} & \multicolumn{1}{c}{} & \multicolumn{1}{c}{} & \multicolumn{1}{c}{\begin{tabular}[c]{@{}c@{}}Route\\ enumeration\end{tabular}} & \multicolumn{1}{c}{Total} \\ \hline \up
		\multirow{66}{*}{100} & \multirow{22}{*}{4} & C0-100 & 854 & 63 & 488119 & 0.00 & 34 & 34 \\  
		&  & C1-100 & 456 & 75 & 331497 & 0.00 & 31 & 31 \\  
		&  & C2-100 & 3802 & 46 & 596824 & 0.00 & 36 & 36 \\  
		&  & C3-100 & 4510 & 46 & 586434 & 0.00 & 35 & 36 \\  
		&  & C4-100 & 4232 & 44 & 558323 & 0.00 & 35 & 41 \\  
		&  & C5-100 & 732 & 70 & 366565 & 0.00 & 31 & 31 \\  
		&  & C6-100 & 2579 & 47 & 724844 & 0.00 & 35 & 35 \\  
		&  & C7-100 & 2020 & 53 & 661648 & 0.00 & 34 & 34 \\  
		&  & C8-100 & 1803 & 51 & 609120 & 0.00 & 34 & 34 \\  
		&  & C9-100 & 12034 & 40 & 527217 & 0.00 & 36 & 38 \\  
		&  & C10-100 & 1095 & 58 & 426386 & 0.00 & 34 & 34 \\  
		&  & C11-100 & 2374 & 51 & 427829 & 0.00 & 34 & 35 \\  
		&  & C12-100 & 989 & 62 & 413012 & 0.00 & 33 & 33 \\  
		&  & C13-100 & 909 & 61 & 483087 & 0.00 & 34 & 34 \\  
		&  & C14-100 & 4306 & 40 & 693825 & 0.00 & 35 & 40 \\  
		&  & C15-100 & 4605 & 48 & 790005 & 0.00 & 34 & 37 \\  
		&  & C16-100 & 1578 & 55 & 627423 & 0.00 & 33 & 33 \\  
		&  & C17-100 & 1151 & 60 & 640945 & 0.00 & 33 & 33 \\  
		&  & C18-100 & 952 & 58 & 681240 & 0.00 & 34 & 34 \\  
		&  & C19-100 & 2226 & 51 & 503252 & 0.00 & 34 & 35 \\  
		&  & C20-100 & 4254 & 46 & 569724 & 0.00 & 35 & 36 \\  
		&  & C21-100 & 667 & 78 & 328216 & 0.00 & 32 & 32 \\ \cline{2-9} 
		& \multirow{22}{*}{5} & C0-100 & 854 & 63 & 488119 & 0.00 & 1185 & 1186 \\  
		&  & C1-100 & 456 & 75 & 331497 & 0.00 & 1072 & 1072 \\  
		&  & C2-100 & 3910 & 46 & 596658 & 0.00 & 1356 & 1357 \\  
		&  & C3-100 & 4814 & 46 & 584077 & 0.00 & 1357 & 1359 \\  
		&  & C4-100 & 4487 & 44 & 552862 & 0.00 & 1298 & 1342 \\  
		&  & C5-100 & 735 & 70 & 366565 & 0.00 & 1034 & 1034 \\  
		&  & C6-100 & 2596 & 47 & 724844 & 0.00 & 1362 & 1362 \\  
		&  & C7-100 & 2043 & 53 & 661528 & 0.00 & 1326 & 1326 \\  
		&  & C8-100 & 1842 & 51 & 609038 & 0.00 & 1317 & 1317 \\  
		&  & C9-100 & 16749 & 40 & 519564 & 0.00 & 1419 & 1451 \\  
		&  & C10-100 & 1096 & 58 & 426386 & 0.00 & 1254 & 1254 \\  
		&  & C11-100 & 2418 & 51 & 427829 & 0.00 & 1320 & 1321 \\  
		&  & C12-100 & 990 & 62 & 413012 & 0.00 & 1208 & 1208 \\  
		&  & C13-100 & 912 & 61 & 483087 & 0.00 & 1221 & 1221 \\  
		&  & C14-100 & 4448 & 40 & 689319 & 0.00 & 1383 & 1457 \\  
		&  & C15-100 & 5985 & 48 & 788340 & 0.00 & 1287 & 1294 \\  
		&  & C16-100 & 1599 & 55 & 626787 & 0.00 & 1226 & 1227 \\  
		&  & C17-100 & 1155 & 60 & 640945 & 0.00 & 1218 & 1218 \\  
		&  & C18-100 & 953 & 58 & 681240 & 0.00 & 1292 & 1293 \\  
		&  & C19-100 & 2291 & 51 & 503252 & 0.00 & 1352 & 1352 \\  
		&  & C20-100 & 4761 & 46 & 568797 & 0.00 & 1352 & 1354 \\  
		&  & C21-100 & 667 & 78 & 328216 & 0.00 & 1080 & 1080 \\ \cline{2-9} 
		& \multirow{22}{*}{6} & C0-100 & 854 & 63 & 488119 & 0.00 & 25088 & 25088 \\  
		&  & C1-100 & 456 & 75 & 331497 & 0.00 & 20117 & 20117 \\  
		&  & C2-100 & 3917 & 46 & 596658 & 0.00 & 33455 & 33456 \\  
		&  & C3-100 & 4872 & 46 & 584048 & 0.00 & 33801 & 33802 \\  
		&  & C4-100 & 4510 & 44 & 552862 & 0.00 & 29993 & 30180 \\  
		&  & C5-100 & 735 & 70 & 366565 & 0.00 & 18613 & 18613 \\  
		&  & C6-100 & 2598 & 47 & 724844 & 0.00 & 43135 & 43136 \\  
		&  & C7-100 & 2043 & 53 & 661528 & 0.00 & 38791 & 38792 \\  
		&  & C8-100 & 1848 & 51 & 609038 & 0.00 & 39437 & 39437 \\  
		&  & C9-100 & 18758 & 40 & 517624 & 0.00 & 48359 & 48403 \\  
		&  & C10-100 & 1096 & 58 & 426386 & 0.00 & 36941 & 36941 \\  
		&  & C11-100 & 2419 & 51 & 427829 & 0.00 & 39761 & 39763 \\  
		&  & C12-100 & 990 & 62 & 413012 & 0.00 & 31038 & 31038 \\  
		&  & C13-100 & 912 & 61 & 483087 & 0.00 & 33881 & 33881 \\  
		&  & C14-100 & 4456 & 40 & 689319 & 0.00 & 41455 & 41492 \\  
		&  & C15-100 & 6939 & 48 & 788340 & 0.00 & 36980 & 36990 \\  
		&  & C16-100 & 1600 & 55 & 626787 & 0.00 & 32031 & 32031 \\  
		&  & C17-100 & 1155 & 60 & 640945 & 0.00 & 34810 & 34810 \\  
		&  & C18-100 & 953 & 58 & 681240 & 0.00 & 38501 & 38501 \\  
		&  & C19-100 & 2294 & 51 & 503252 & 0.00 & 42573 & 42574 \\  
		&  & C20-100 & 4898 & 46 & 568797 & 0.00 & 43900 & 43903 \\ \down
		&  & C21-100 & 667 & 78 & 328216 & 0.00 & 29507 & 29507 \\ \hline
	\end{tabular}}
\end{table}

\begin{table}[]
	\TABLE
	{Results of BPA Scalability with Increasing Vehicle Capacity ($\Delta = 10$ mins, $R = 50\%$).\label{tbl:bpa_vehicle_capacity_scaling}}
	{\resizebox*{!}{\dimexpr\textheight-1\baselineskip\relax}{
			\begin{tabular}{ccrrrrrrrrrrrrrr}
				\hline \up\down
				\multirow{2}{*}{\begin{tabular}[c]{@{}c@{}}Cluster\\ size\end{tabular}} & \multirow{2}{*}{\begin{tabular}[c]{@{}c@{}}Vehicle\\ capacity\end{tabular}} & \multicolumn{1}{c}{\multirow{2}{*}{\begin{tabular}[c]{@{}c@{}}Cluster\\ ID\end{tabular}}} & \multicolumn{1}{c}{\multirow{2}{*}{\begin{tabular}[c]{@{}c@{}}Inbound\\ edge \#\end{tabular}}} & \multicolumn{1}{c}{\multirow{2}{*}{\begin{tabular}[c]{@{}c@{}}Outbound\\ edge \#\end{tabular}}} & \multicolumn{1}{c}{\multirow{2}{*}{\begin{tabular}[c]{@{}c@{}}Tree\\ node\\ \#\end{tabular}}} & \multicolumn{1}{c}{\multirow{2}{*}{\begin{tabular}[c]{@{}c@{}}Column\\ \#\end{tabular}}} & \multicolumn{1}{c}{\multirow{2}{*}{\begin{tabular}[c]{@{}c@{}}Vehicle\\ \#\end{tabular}}} & \multicolumn{1}{c}{\multirow{2}{*}{\begin{tabular}[c]{@{}c@{}}Total\\ distance\\ (m)\end{tabular}}} & \multicolumn{2}{c}{Optimality gap (\%)} & \multicolumn{1}{c}{\multirow{2}{*}{\begin{tabular}[c]{@{}c@{}}Integrality\\ gap (\%)\end{tabular}}} & \multicolumn{4}{c}{Wall time (s)} \\ \cline{10-11} \cline{13-16} 
				&  & \multicolumn{1}{c}{} & \multicolumn{1}{c}{} & \multicolumn{1}{c}{} & \multicolumn{1}{c}{} & \multicolumn{1}{c}{} & \multicolumn{1}{c}{} & \multicolumn{1}{c}{} & \multicolumn{1}{c}{\begin{tabular}[c]{@{}c@{}}Root\\ node soln.\end{tabular}} & \multicolumn{1}{c}{\begin{tabular}[c]{@{}c@{}}Best\\ feasible\\ soln.\end{tabular}} & \multicolumn{1}{c}{} & \multicolumn{1}{c}{\begin{tabular}[c]{@{}c@{}}RMP\\ conv.\end{tabular}} & \multicolumn{1}{c}{\begin{tabular}[c]{@{}c@{}}Root\\ node\\ soln.\end{tabular}} & \multicolumn{1}{c}{\begin{tabular}[c]{@{}c@{}}Best\\ feasible\\ soln.\end{tabular}} & \multicolumn{1}{c}{Total} \\ \hline \up
				\multirow{72}{*}{75} & \multirow{24}{*}{4} & C0-75 & 1661 & 1565 & 15 & 462 & 47 & 356296 & 0.00 & 0.00 & 0.00 & 1 & 1 & 1 & 1 \\ 
				&  & C2-75 & 1906 & 1405 & 57 & 1043 & 40 & 643541 & 0.01 & 0.00 & 0.01 & 2 & 2 & 7 & 12 \\ 
				&  & C3-75 & 3099 & 1809 & 31 & 1247 & 38 & 487061 & 0.00 & 0.00 & 0.00 & 2 & 2 & 12 & 13 \\ 
				&  & C4-75 & 2052 & 1851 & 79 & 915 & 43 & 398265 & 0.01 & 0.00 & 0.01 & 1 & 1 & 8 & 12 \\ 
				&  & C5-75 & 823 & 810 & 1 & 267 & 61 & 230460 & 0.00 & 0.00 & 0.00 & 1 & 1 & 1 & 1 \\ 
				&  & C6-75 & 3427 & 2012 & 17 & 1420 & 36 & 451262 & 2.77 & 0.00 & 2.77 & 4 & 5 & 5 & 16 \\ 
				&  & C7-75 & 2458 & 1649 & 25 & 615 & 46 & 437743 & 2.17 & 0.00 & 2.17 & 1 & 1 & 1 & 2 \\ 
				&  & C8-75 & 3002 & 2183 & 27 & 1028 & 37 & 487561 & 2.69 & 0.00 & 2.69 & 2 & 2 & 11 & 15 \\ 
				&  & C9-75 & 3550 & 2341 & 428 & 2613 & 31 & 328525 & 3.22 & 0.00 & 3.21 & 7 & 8 & 493 & 528 \\ 
				&  & C10-75 & 2856 & 1580 & 151 & 2085 & 32 & 527544 & 3.11 & 0.00 & 3.11 & 11 & 11 & 168 & 232 \\ 
				&  & C11-75 & 1955 & 1193 & 23 & 639 & 47 & 443770 & 0.00 & 0.00 & 0.00 & 1 & 1 & 1 & 2 \\ 
				&  & C12-75 & 1801 & 1347 & 5 & 501 & 46 & 298162 & 0.00 & 0.00 & 0.00 & 1 & 1 & 1 & 1 \\ 
				&  & C13-75 & 2699 & 1302 & 141 & 1392 & 37 & 481869 & 2.70 & 0.00 & 2.70 & 2 & 2 & 17 & 37 \\ 
				&  & C14-75 & 2045 & 1332 & 17 & 626 & 47 & 301108 & 0.00 & 0.00 & 0.00 & 1 & 1 & 1 & 2 \\ 
				&  & C15-75 & 3007 & 1443 & 23 & 1098 & 36 & 460087 & 0.01 & 0.00 & 0.01 & 2 & 2 & 2 & 11 \\ 
				&  & C17-75 & 1492 & 835 & 3 & 404 & 50 & 348581 & 0.00 & 0.00 & 0.00 & 1 & 1 & 1 & 1 \\ 
				&  & C19-75 & 1497 & 664 & 1 & 390 & 57 & 245827 & 0.00 & 0.00 & 0.00 & 1 & 1 & 1 & 1 \\ 
				&  & C20-75 & 2323 & 1254 & 17 & 683 & 46 & 420825 & 0.00 & 0.00 & 0.00 & 1 & 1 & 2 & 2 \\ 
				&  & C22-75 & 3113 & 1524 & 119 & 1373 & 40 & 472247 & 2.49 & 0.00 & 2.49 & 2 & 2 & 17 & 45 \\ 
				&  & C23-75 & 2069 & 1548 & 169 & 1105 & 36 & 385434 & 5.53 & 0.00 & 5.53 & 2 & 2 & 28 & 45 \\ 
				&  & C24-75 & 2446 & 1694 & 305 & 1479 & 37 & 425807 & 2.70 & 0.00 & 2.69 & 2 & 3 & 72 & 94 \\ 
				&  & C26-75 & 1957 & 1039 & 17 & 432 & 52 & 333692 & 0.00 & 0.00 & 0.00 & 1 & 1 & 1 & 1 \\ 
				&  & C28-75 & 1781 & 1452 & 5 & 611 & 44 & 293838 & 0.00 & 0.00 & 0.00 & 1 & 1 & 1 & 1 \\ 
				&  & C29-75 & 2495 & 2238 & 1242 & 2217 & 32 & 548779 & 3.03 & 0.00 & 0.02 & 4 & 19 & 623 & 805 \\ \cline{2-16} 
				& \multirow{24}{*}{5} & C0-75 & 1661 & 1565 & 15 & 462 & 47 & 356296 & 0.00 & 0.00 & 0.00 & 1 & 1 & 1 & 1 \\ 
				&  & C2-75 & 1906 & 1405 & 115 & 1085 & 40 & 643540 & 0.02 & 0.00 & 0.01 & 2 & 2 & 33 & 42 \\ 
				&  & C3-75 & 3099 & 1809 & 37 & 1274 & 38 & 487061 & 0.00 & 0.00 & 0.00 & 3 & 3 & 22 & 24 \\ 
				&  & C4-75 & 2052 & 1851 & 75 & 946 & 43 & 398216 & 2.27 & 0.00 & 0.01 & 1 & 2 & 14 & 19 \\ 
				&  & C5-75 & 823 & 810 & 1 & 267 & 61 & 230460 & 0.00 & 0.00 & 0.00 & 1 & 1 & 1 & 1 \\ 
				&  & C6-75 & 3427 & 2012 & 39 & 1519 & 35 & 450414 & 2.85 & 0.00 & 2.85 & 48 & 48 & 107 & 243 \\ 
				&  & C7-75 & 2458 & 1649 & 21 & 619 & 46 & 437743 & 2.17 & 0.00 & 2.17 & 1 & 1 & 1 & 2 \\ 
				&  & C8-75 & 3002 & 2183 & 41 & 1049 & 37 & 487245 & 2.69 & 0.00 & 2.69 & 2 & 2 & 26 & 27 \\ 
				&  & C9-75 & 3550 & 2341 & 399 & 2597 & 31 & 327993 & 3.22 & 0.00 & 3.22 & 14 & 15 & 721 & 781 \\ 
				&  & C10-75 & 2856 & 1580 & 197 & 2350 & 32 & 526188 & 3.12 & 0.00 & 3.12 & 130 & 131 & 382 & 1854 \\ 
				&  & C11-75 & 1955 & 1193 & 23 & 645 & 47 & 443770 & 0.00 & 0.00 & 0.00 & 1 & 1 & 1 & 3 \\ 
				&  & C12-75 & 1801 & 1347 & 5 & 506 & 46 & 298162 & 0.00 & 0.00 & 0.00 & 1 & 1 & 1 & 1 \\ 
				&  & C13-75 & 2699 & 1302 & 161 & 1367 & 37 & 479166 & 2.70 & 0.00 & 2.70 & 3 & 4 & 59 & 67 \\ 
				&  & C14-75 & 2045 & 1333 & 21 & 626 & 47 & 301108 & 0.00 & 0.00 & 0.00 & 1 & 1 & 2 & 3 \\ 
				&  & C15-75 & 3003 & 1443 & 35 & 1117 & 36 & 460087 & 0.01 & 0.00 & 0.01 & 3 & 3 & 3 & 23 \\ 
				&  & C17-75 & 1492 & 835 & 3 & 404 & 50 & 348581 & 0.00 & 0.00 & 0.00 & 1 & 1 & 1 & 1 \\ 
				&  & C19-75 & 1497 & 664 & 1 & 390 & 57 & 245827 & 0.00 & 0.00 & 0.00 & 1 & 1 & 1 & 1 \\ 
				&  & C20-75 & 2323 & 1254 & 17 & 683 & 46 & 420825 & 0.00 & 0.00 & 0.00 & 1 & 1 & 2 & 2 \\ 
				&  & C22-75 & 3113 & 1524 & 199 & 1399 & 40 & 472040 & 2.50 & 0.00 & 2.49 & 3 & 3 & 19 & 105 \\ 
				&  & C23-75 & 2069 & 1548 & 25 & 983 & 35 & 379234 & 2.84 & 0.00 & 2.84 & 3 & 3 & 13 & 18 \\ 
				&  & C24-75 & 2446 & 1694 & 705 & 1620 & 37 & 425806 & 5.25 & 0.00 & 2.70 & 3 & 6 & 109 & 252 \\ 
				&  & C26-75 & 1957 & 1039 & 17 & 430 & 52 & 333692 & 0.00 & 0.00 & 0.00 & 1 & 1 & 1 & 1 \\ 
				&  & C28-75 & 1781 & 1452 & 5 & 628 & 44 & 293838 & 0.00 & 0.00 & 0.00 & 1 & 1 & 1 & 1 \\ 
				&  & C29-75 & 2495 & 2238 & 619 & 2097 & 32 & 543850 & 0.02 & 0.00 & 0.02 & 8 & 8 & 454 & 578 \\ \cline{2-16} 
				& \multirow{24}{*}{6} & C0-75 & 1661 & 1565 & 15 & 462 & 47 & 356296 & 0.00 & 0.00 & 0.00 & 1 & 1 & 1 & 1 \\ 
				&  & C2-75 & 1906 & 1405 & 87 & 1072 & 40 & 643540 & 0.02 & 0.00 & 0.01 & 2 & 3 & 27 & 39 \\ 
				&  & C3-75 & 3099 & 1809 & 35 & 1281 & 38 & 487061 & 0.00 & 0.00 & 0.00 & 4 & 4 & 30 & 33 \\ 
				&  & C4-75 & 2052 & 1851 & 75 & 929 & 43 & 398216 & 2.27 & 0.00 & 0.01 & 2 & 2 & 16 & 20 \\ 
				&  & C5-75 & 823 & 810 & 1 & 267 & 61 & 230460 & 0.00 & 0.00 & 0.00 & 1 & 1 & 1 & 1 \\ 
				&  & C6-75 & 3427 & 2012 & 53 & 1557 & 35 & 450414 & 2.85 & 0.00 & 2.85 & 42 & 42 & 224 & 336 \\ 
				&  & C7-75 & 2458 & 1649 & 29 & 617 & 46 & 437743 & 2.17 & 0.00 & 2.17 & 1 & 1 & 1 & 2 \\ 
				&  & C8-75 & 3002 & 2183 & 19 & 1019 & 37 & 487245 & 2.69 & 0.00 & 2.69 & 2 & 2 & 10 & 14 \\ 
				&  & C9-75 & 3550 & 2341 & 103 & 2127 & 31 & 327993 & 3.22 & 0.00 & 3.22 & 17 & 18 & 246 & 339 \\ 
				&  & C10-75 & 2856 & 1580 & 73 & 2029 & 32 & 525115 & 3.11 & 0.00 & 3.11 & 2763 & 2764 & 2764 & 17162 \\ 
				&  & C11-75 & 1955 & 1193 & 23 & 648 & 47 & 443770 & 0.00 & 0.00 & 0.00 & 1 & 1 & 1 & 3 \\ 
				&  & C12-75 & 1801 & 1347 & 5 & 505 & 46 & 298162 & 0.00 & 0.00 & 0.00 & 1 & 1 & 1 & 1 \\ 
				&  & C13-75 & 2699 & 1302 & 369 & 1484 & 37 & 479166 & 2.70 & 0.00 & 2.70 & 5 & 5 & 88 & 151 \\ 
				&  & C14-75 & 2045 & 1333 & 15 & 620 & 47 & 301108 & 0.00 & 0.00 & 0.00 & 1 & 1 & 1 & 2 \\ 
				&  & C15-75 & 3007 & 1443 & 19 & 1061 & 36 & 460087 & 0.01 & 0.00 & 0.01 & 3 & 3 & 3 & 13 \\ 
				&  & C17-75 & 1492 & 835 & 5 & 411 & 50 & 348581 & 0.00 & 0.00 & 0.00 & 1 & 1 & 1 & 1 \\ 
				&  & C19-75 & 1497 & 664 & 1 & 390 & 57 & 245827 & 0.00 & 0.00 & 0.00 & 1 & 1 & 1 & 1 \\ 
				&  & C20-75 & 2323 & 1254 & 17 & 671 & 46 & 420825 & 0.00 & 0.00 & 0.00 & 1 & 1 & 2 & 2 \\ 
				&  & C22-75 & 3113 & 1524 & 141 & 1358 & 40 & 472040 & 2.50 & 0.00 & 2.49 & 3 & 4 & 53 & 84 \\ 
				&  & C23-75 & 2069 & 1548 & 35 & 981 & 35 & 379234 & 2.84 & 0.00 & 2.84 & 3 & 3 & 16 & 26 \\ 
				&  & C24-75 & 2446 & 1694 & 613 & 1572 & 37 & 425806 & 5.25 & 0.00 & 2.70 & 4 & 9 & 94 & 256 \\ 
				&  & C26-75 & 1957 & 1039 & 17 & 430 & 52 & 333692 & 0.00 & 0.00 & 0.00 & 1 & 1 & 1 & 1 \\ 
				&  & C28-75 & 1781 & 1452 & 5 & 622 & 44 & 293838 & 0.00 & 0.00 & 0.00 & 1 & 1 & 1 & 1 \\ \down
				&  & C29-75 & 2495 & 2238 & 555 & 2050 & 32 & 543850 & 0.02 & 0.00 & 0.02 & 9 & 9 & 9 & 667 \\ \hline
			\end{tabular}
	}}
	{}
\end{table}

\begin{table}[]
	\resizebox{\textwidth}{!}
	{
		\begin{tabular}{ccrrrrrrrrrrrrrr}
			\hline \up\down
			\multirow{2}{*}{\begin{tabular}[c]{@{}c@{}}Cluster\\ size\end{tabular}} & \multirow{2}{*}{\begin{tabular}[c]{@{}c@{}}Vehicle\\ capacity\end{tabular}} & \multicolumn{1}{c}{\multirow{2}{*}{\begin{tabular}[c]{@{}c@{}}Cluster\\ ID\end{tabular}}} & \multicolumn{1}{c}{\multirow{2}{*}{\begin{tabular}[c]{@{}c@{}}Inbound\\ edge \#\end{tabular}}} & \multicolumn{1}{c}{\multirow{2}{*}{\begin{tabular}[c]{@{}c@{}}Outbound\\ edge \#\end{tabular}}} & \multicolumn{1}{c}{\multirow{2}{*}{\begin{tabular}[c]{@{}c@{}}Tree\\ node\\ \#\end{tabular}}} & \multicolumn{1}{c}{\multirow{2}{*}{\begin{tabular}[c]{@{}c@{}}Column\\ \#\end{tabular}}} & \multicolumn{1}{c}{\multirow{2}{*}{\begin{tabular}[c]{@{}c@{}}Vehicle\\ \#\end{tabular}}} & \multicolumn{1}{c}{\multirow{2}{*}{\begin{tabular}[c]{@{}c@{}}Total\\ distance\\ (m)\end{tabular}}} & \multicolumn{2}{c}{Optimality gap (\%)} & \multicolumn{1}{c}{\multirow{2}{*}{\begin{tabular}[c]{@{}c@{}}Integrality\\ gap (\%)\end{tabular}}} & \multicolumn{4}{c}{Wall time (s)} \\ \cline{10-11} \cline{13-16} 
			&  & \multicolumn{1}{c}{} & \multicolumn{1}{c}{} & \multicolumn{1}{c}{} & \multicolumn{1}{c}{} & \multicolumn{1}{c}{} & \multicolumn{1}{c}{} & \multicolumn{1}{c}{} & \multicolumn{1}{c}{\begin{tabular}[c]{@{}c@{}}Root\\ node soln.\end{tabular}} & \multicolumn{1}{c}{\begin{tabular}[c]{@{}c@{}}Best\\ feasible\\ soln.\end{tabular}} & \multicolumn{1}{c}{} & \multicolumn{1}{c}{\begin{tabular}[c]{@{}c@{}}RMP\\ conv.\end{tabular}} & \multicolumn{1}{c}{\begin{tabular}[c]{@{}c@{}}Root\\ node\\ soln.\end{tabular}} & \multicolumn{1}{c}{\begin{tabular}[c]{@{}c@{}}Best\\ feasible\\ soln.\end{tabular}} & \multicolumn{1}{c}{Total} \\ \hline \up
			\multirow{48}{*}{75} & \multirow{24}{*}{7} & C0-75 & 1661 & 1565 & 15 & 463 & 47 & 356296 & 0.00 & 0.00 & 0.00 & 1 & 1 & 1 & 1 \\ 
			&  & C2-75 & 1906 & 1405 & 79 & 1072 & 40 & 643540 & 0.02 & 0.00 & 0.01 & 3 & 3 & 30 & 39 \\ 
			&  & C3-75 & 3099 & 1809 & 21 & 1238 & 38 & 487061 & 0.00 & 0.00 & 0.00 & 4 & 4 & 4 & 20 \\ 
			&  & C4-75 & 2052 & 1851 & 139 & 962 & 43 & 398216 & 2.27 & 0.00 & 0.01 & 1 & 2 & 33 & 33 \\ 
			&  & C5-75 & 823 & 810 & 1 & 267 & 61 & 230460 & 0.00 & 0.00 & 0.00 & 1 & 1 & 1 & 1 \\ 
			&  & C6-75 & 3427 & 2012 & 33 & 1515 & 35 & 450414 & 2.85 & 0.00 & 2.85 & 48 & 48 & 186 & 247 \\ 
			&  & C7-75 & 2458 & 1649 & 31 & 616 & 46 & 437743 & 2.17 & 0.00 & 2.17 & 1 & 1 & 1 & 3 \\ 
			&  & C8-75 & 3002 & 2183 & 23 & 1025 & 37 & 487245 & 2.69 & 0.00 & 2.69 & 2 & 2 & 9 & 18 \\ 
			&  & C9-75 & 3550 & 2341 & 189 & 2347 & 31 & 327993 & 3.22 & 0.00 & 3.22 & 20 & 20 & 482 & 594 \\ 
			&  & C10-75 & 2856 & 1580 & 2 & 1459 & 32 & 525115 & 3.12 & 3.12 & 3.12 & 39153 & 39154 & 39154 & 43200 \\ 
			&  & C11-75 & 1955 & 1193 & 24 & 650 & 47 & 443770 & 0.00 & 0.00 & 0.00 & 1 & 1 & 1 & 3 \\ 
			&  & C12-75 & 1801 & 1347 & 5 & 502 & 46 & 298162 & 0.00 & 0.00 & 0.00 & 1 & 1 & 1 & 1 \\ 
			&  & C13-75 & 2699 & 1302 & 337 & 1470 & 37 & 479166 & 2.70 & 0.00 & 2.70 & 2 & 3 & 78 & 105 \\ 
			&  & C14-75 & 2045 & 1332 & 19 & 619 & 47 & 301108 & 0.00 & 0.00 & 0.00 & 1 & 1 & 2 & 2 \\ 
			&  & C15-75 & 3007 & 1443 & 25 & 1089 & 36 & 460087 & 0.01 & 0.00 & 0.01 & 2 & 3 & 3 & 17 \\ 
			&  & C17-75 & 1492 & 835 & 7 & 411 & 50 & 348581 & 0.00 & 0.00 & 0.00 & 1 & 1 & 1 & 1 \\ 
			&  & C19-75 & 1497 & 664 & 1 & 390 & 57 & 245827 & 0.00 & 0.00 & 0.00 & 1 & 1 & 1 & 1 \\ 
			&  & C20-75 & 2323 & 1254 & 11 & 672 & 46 & 420825 & 0.00 & 0.00 & 0.00 & 1 & 1 & 2 & 2 \\ 
			&  & C22-75 & 3113 & 1524 & 161 & 1382 & 40 & 472040 & 2.50 & 0.00 & 2.49 & 3 & 3 & 52 & 79 \\ 
			&  & C23-75 & 2069 & 1548 & 23 & 940 & 35 & 379234 & 2.84 & 0.00 & 2.84 & 3 & 3 & 12 & 19 \\ 
			&  & C24-75 & 2446 & 1694 & 745 & 1666 & 37 & 425806 & 5.25 & 0.00 & 2.70 & 3 & 6 & 148 & 329 \\ 
			&  & C26-75 & 1957 & 1039 & 17 & 430 & 52 & 333692 & 0.00 & 0.00 & 0.00 & 1 & 1 & 1 & 1 \\ 
			&  & C28-75 & 1781 & 1452 & 5 & 621 & 44 & 293838 & 0.00 & 0.00 & 0.00 & 1 & 1 & 1 & 1 \\ 
			&  & C29-75 & 2495 & 2238 & 523 & 2052 & 32 & 543850 & 0.02 & 0.00 & 0.02 & 7 & 8 & 8 & 619 \\ \cline{2-16} 
			& \multirow{24}{*}{8} & C0-75 & 1661 & 1565 & 15 & 461 & 47 & 356296 & 0.00 & 0.00 & 0.00 & 1 & 1 & 1 & 1 \\ 
			&  & C2-75 & 1906 & 1405 & 75 & 1047 & 40 & 643540 & 0.02 & 0.00 & 0.01 & 2 & 2 & 26 & 36 \\ 
			&  & C3-75 & 3099 & 1809 & 35 & 1268 & 38 & 487061 & 0.00 & 0.00 & 0.00 & 4 & 4 & 30 & 33 \\ 
			&  & C4-75 & 2052 & 1851 & 101 & 928 & 43 & 398216 & 0.01 & 0.00 & 0.01 & 2 & 2 & 27 & 28 \\ 
			&  & C5-75 & 823 & 810 & 1 & 267 & 61 & 230460 & 0.00 & 0.00 & 0.00 & 1 & 1 & 1 & 1 \\ 
			&  & C6-75 & 3427 & 2012 & 62 & 1584 & 35 & 450414 & 5.54 & 0.00 & 2.85 & 39 & 40 & 266 & 362 \\ 
			&  & C7-75 & 2458 & 1649 & 31 & 618 & 46 & 437743 & 2.17 & 0.00 & 2.17 & 1 & 1 & 1 & 2 \\ 
			&  & C8-75 & 3002 & 2183 & 49 & 1064 & 37 & 487245 & 2.69 & 0.00 & 2.69 & 3 & 4 & 9 & 81 \\ 
			&  & C9-75 & 3550 & 2341 & 105 & 2199 & 31 & 327993 & 3.22 & 0.00 & 3.22 & 26 & 27 & 300 & 398 \\ 
			&  & C10-75 & 2856 & 1580 & 1 & 1425 & 32 & 525207 & 3.12 & 3.12 & 3.12 & 84928 & 84928 & 84928 & 84928 \\ 
			&  & C11-75 & 1955 & 1193 & 24 & 641 & 47 & 443770 & 0.00 & 0.00 & 0.00 & 1 & 1 & 1 & 3 \\ 
			&  & C12-75 & 1801 & 1347 & 5 & 508 & 46 & 298162 & 0.00 & 0.00 & 0.00 & 1 & 1 & 1 & 1 \\ 
			&  & C13-75 & 2699 & 1302 & 427 & 1513 & 37 & 479166 & 2.70 & 0.00 & 2.70 & 4 & 4 & 77 & 144 \\ 
			&  & C14-75 & 2045 & 1332 & 17 & 624 & 47 & 301108 & 0.00 & 0.00 & 0.00 & 1 & 1 & 1 & 2 \\ 
			&  & C15-75 & 3005 & 1443 & 21 & 1093 & 36 & 460087 & 0.01 & 0.00 & 0.01 & 3 & 3 & 3 & 16 \\ 
			&  & C17-75 & 1492 & 835 & 7 & 410 & 50 & 348581 & 0.00 & 0.00 & 0.00 & 1 & 1 & 1 & 1 \\ 
			&  & C19-75 & 1497 & 664 & 1 & 390 & 57 & 245827 & 0.00 & 0.00 & 0.00 & 1 & 1 & 1 & 1 \\ 
			&  & C20-75 & 2323 & 1254 & 11 & 671 & 46 & 420825 & 0.00 & 0.00 & 0.00 & 1 & 1 & 2 & 2 \\ 
			&  & C22-75 & 3113 & 1524 & 161 & 1357 & 40 & 472040 & 2.50 & 0.00 & 2.49 & 3 & 4 & 65 & 103 \\ 
			&  & C23-75 & 2069 & 1548 & 25 & 959 & 35 & 379234 & 2.84 & 0.00 & 2.84 & 4 & 4 & 16 & 22 \\ 
			&  & C24-75 & 2446 & 1694 & 677 & 1567 & 37 & 425806 & 2.70 & 0.00 & 2.70 & 3 & 3 & 83 & 265 \\ 
			&  & C26-75 & 1957 & 1039 & 15 & 430 & 52 & 333692 & 0.00 & 0.00 & 0.00 & 1 & 1 & 1 & 1 \\ 
			&  & C28-75 & 1781 & 1452 & 5 & 620 & 44 & 293838 & 0.00 & 0.00 & 0.00 & 1 & 1 & 1 & 1 \\ \down
			&  & C29-75 & 2495 & 2238 & 573 & 2163 & 32 & 543850 & 0.02 & 0.00 & 0.02 & 9 & 9 & 353 & 625 \\ \hline \up
			\multirow{22}{*}{100} & \multirow{22}{*}{4} & C0-100 & 3488 & 1930 & 45 & 777 & 63 & 488119 & 0.00 & 0.00 & 0.00 & 2 & 2 & 2 & 5 \\ 
			&  & C1-100 & 1825 & 1642 & 5 & 429 & 75 & 331497 & 0.00 & 0.00 & 0.00 & 2 & 2 & 2 & 3 \\ 
			&  & C2-100 & 5540 & 3106 & 1445 & 2652 & 46 & 596824 & 2.17 & 0.00 & 2.17 & 6 & 7 & 995 & 1397 \\ 
			&  & C3-100 & 5383 & 3127 & 823 & 2794 & 46 & 586434 & 2.17 & 0.00 & 2.16 & 8 & 9 & 558 & 695 \\ 
			&  & C4-100 & 4211 & 2787 & 50653 & 6120 & 44 & 558323 & 2.27 & 0.00 & 2.26 & 8 & 9 & 730 & 31745 \\ 
			&  & C5-100 & 1960 & 1452 & 39 & 575 & 70 & 366565 & 0.00 & 0.00 & 0.00 & 2 & 2 & 2 & 4 \\ 
			&  & C6-100 & 4708 & 3393 & 33 & 1805 & 47 & 724844 & 0.00 & 0.00 & 0.00 & 8 & 8 & 39 & 44 \\ 
			&  & C7-100 & 4120 & 2205 & 23 & 1227 & 53 & 661648 & 0.00 & 0.00 & 0.00 & 3 & 3 & 3 & 10 \\ 
			&  & C8-100 & 4952 & 2656 & 3 & 1285 & 51 & 609120 & 0.00 & 0.00 & 0.00 & 3 & 3 & 3 & 4 \\ 
			&  & C9-100 & 5664 & 3572 & 2316 & 6235 & 40 & 527217 & 2.50 & 0.00 & 2.49 & 57 & 59 & 7410 & 14795 \\ 
			&  & C10-100 & 2995 & 2611 & 147 & 980 & 58 & 426386 & 1.72 & 0.00 & 1.72 & 2 & 2 & 14 & 16 \\ 
			&  & C11-100 & 4606 & 2964 & 61 & 1506 & 51 & 427829 & 0.00 & 0.00 & 0.00 & 3 & 4 & 14 & 32 \\ 
			&  & C12-100 & 2863 & 2459 & 1 & 794 & 62 & 413012 & 0.00 & 0.00 & 0.00 & 2 & 2 & 2 & 2 \\ 
			&  & C13-100 & 3232 & 2441 & 17 & 755 & 61 & 483087 & 0.00 & 0.00 & 0.00 & 2 & 2 & 3 & 3 \\ 
			&  & C14-100 & 4335 & 3403 & 2695 & 3757 & 40 & 693825 & 4.85 & 0.00 & 2.49 & 15 & 1297 & 3093 & 6348 \\ 
			&  & C15-100 & 3711 & 2642 & 26493 & 3217 & 48 & 790005 & 4.16 & 0.00 & 4.16 & 7 & 8 & 2068 & 26100 \\ 
			&  & C16-100 & 3278 & 2302 & 21 & 1085 & 55 & 627423 & 1.81 & 0.00 & 1.81 & 3 & 3 & 5 & 8 \\ 
			&  & C17-100 & 3413 & 1603 & 7 & 881 & 60 & 640945 & 0.00 & 0.00 & 0.00 & 2 & 2 & 2 & 3 \\ 
			&  & C18-100 & 3778 & 2518 & 81 & 872 & 58 & 681240 & 1.72 & 0.00 & 1.72 & 2 & 2 & 2 & 10 \\ 
			&  & C19-100 & 3722 & 3405 & 353 & 1656 & 51 & 503252 & 1.96 & 0.00 & 1.96 & 4 & 4 & 64 & 192 \\ 
			&  & C20-100 & 4377 & 2653 & 1863 & 2708 & 46 & 569724 & 2.17 & 0.00 & 2.17 & 7 & 7 & 1135 & 1144 \\ \down
			&  & C21-100 & 2597 & 1223 & 1 & 524 & 78 & 328216 & 0.00 & 0.00 & 0.00 & 2 & 2 & 2 & 2 \\ \hline
		\end{tabular}
	}
\end{table}

\begin{table}[]
	\centering
	\resizebox*{!}{\textheight}
	{
\begin{tabular}{ccrrrrrrrrrrrrrr}
	\hline \up\down
	\multirow{2}{*}{\begin{tabular}[c]{@{}c@{}}Cluster\\ size\end{tabular}} & \multirow{2}{*}{\begin{tabular}[c]{@{}c@{}}Vehicle\\ capacity\end{tabular}} & \multicolumn{1}{c}{\multirow{2}{*}{\begin{tabular}[c]{@{}c@{}}Cluster\\ ID\end{tabular}}} & \multicolumn{1}{c}{\multirow{2}{*}{\begin{tabular}[c]{@{}c@{}}Inbound\\ edge \#\end{tabular}}} & \multicolumn{1}{c}{\multirow{2}{*}{\begin{tabular}[c]{@{}c@{}}Outbound\\ edge \#\end{tabular}}} & \multicolumn{1}{c}{\multirow{2}{*}{\begin{tabular}[c]{@{}c@{}}Tree\\ node\\ \#\end{tabular}}} & \multicolumn{1}{c}{\multirow{2}{*}{\begin{tabular}[c]{@{}c@{}}Column\\ \#\end{tabular}}} & \multicolumn{1}{c}{\multirow{2}{*}{\begin{tabular}[c]{@{}c@{}}Vehicle\\ \#\end{tabular}}} & \multicolumn{1}{c}{\multirow{2}{*}{\begin{tabular}[c]{@{}c@{}}Total\\ distance\\ (m)\end{tabular}}} & \multicolumn{2}{c}{Optimality gap (\%)} & \multicolumn{1}{c}{\multirow{2}{*}{\begin{tabular}[c]{@{}c@{}}Integrality\\ gap (\%)\end{tabular}}} & \multicolumn{4}{c}{Wall time (s)} \\ \cline{10-11} \cline{13-16} 
	&  & \multicolumn{1}{c}{} & \multicolumn{1}{c}{} & \multicolumn{1}{c}{} & \multicolumn{1}{c}{} & \multicolumn{1}{c}{} & \multicolumn{1}{c}{} & \multicolumn{1}{c}{} & \multicolumn{1}{c}{\begin{tabular}[c]{@{}c@{}}Root\\ node soln.\end{tabular}} & \multicolumn{1}{c}{\begin{tabular}[c]{@{}c@{}}Best\\ feasible\\ soln.\end{tabular}} & \multicolumn{1}{c}{} & \multicolumn{1}{c}{\begin{tabular}[c]{@{}c@{}}RMP\\ conv.\end{tabular}} & \multicolumn{1}{c}{\begin{tabular}[c]{@{}c@{}}Root\\ node\\ soln.\end{tabular}} & \multicolumn{1}{c}{\begin{tabular}[c]{@{}c@{}}Best\\ feasible\\ soln.\end{tabular}} & \multicolumn{1}{c}{Total} \\ \hline \up
	\multirow{88}{*}{100} & \multirow{22}{*}{5} & C0-100 & 3488 & 1930 & 43 & 772 & 63 & 488119 & 0.00 & 0.00 & 0.00 & 2 & 2 & 2 & 5 \\  
	&  & C1-100 & 1825 & 1642 & 5 & 429 & 75 & 331497 & 0.00 & 0.00 & 0.00 & 2 & 2 & 2 & 2 \\  
	&  & C2-100 & 5540 & 3106 & 1513 & 2637 & 46 & 596658 & 2.17 & 0.00 & 2.17 & 9 & 10 & 1269 & 1942 \\  
	&  & C3-100 & 5383 & 3127 & 191 & 2449 & 46 & 584077 & 2.17 & 0.00 & 2.16 & 60 & 61 & 860 & 950 \\  
	&  & C4-100 & 4209 & 2787 & 67 & 2111 & 44 & 552862 & 2.26 & 0.00 & 2.26 & 18 & 27 & 141 & 183 \\  
	&  & C5-100 & 1964 & 1452 & 47 & 586 & 70 & 366565 & 0.00 & 0.00 & 0.00 & 2 & 2 & 2 & 4 \\  
	&  & C6-100 & 4708 & 3393 & 29 & 1807 & 47 & 724844 & 0.00 & 0.00 & 0.00 & 8 & 8 & 26 & 54 \\  
	&  & C7-100 & 4120 & 2205 & 27 & 1234 & 53 & 661528 & 0.00 & 0.00 & 0.00 & 4 & 4 & 4 & 13 \\  
	&  & C8-100 & 4952 & 2656 & 3 & 1287 & 51 & 609038 & 0.00 & 0.00 & 0.00 & 3 & 3 & 3 & 5 \\  
	&  & C9-100 & 5559 & 3570 & 459 & 4860 & 40 & 521285 & 2.49 & 2.49 & 2.49 & 798 & 799 & 29422 & 43304 \\  
	&  & C10-100 & 2995 & 2611 & 169 & 971 & 58 & 426386 & 1.72 & 0.00 & 1.72 & 2 & 2 & 2 & 17 \\  
	&  & C11-100 & 4606 & 2964 & 53 & 1462 & 51 & 427829 & 0.00 & 0.00 & 0.00 & 4 & 5 & 18 & 38 \\  
	&  & C12-100 & 2863 & 2459 & 1 & 787 & 62 & 413012 & 0.00 & 0.00 & 0.00 & 2 & 2 & 2 & 2 \\  
	&  & C13-100 & 3232 & 2441 & 17 & 758 & 61 & 483087 & 0.00 & 0.00 & 0.00 & 2 & 2 & 2 & 3 \\  
	&  & C14-100 & 4335 & 3403 & 1067 & 3293 & 40 & 690964 & 2.49 & 2.49 & 2.49 & 70 & 72 & 39918 & 43218 \\  
	&  & C15-100 & 3711 & 2642 & 11238 & 4146 & 48 & 788340 & 4.16 & 4.16 & 4.16 & 22 & 22 & 4167 & 43205 \\  
	&  & C16-100 & 3278 & 2302 & 65 & 1151 & 55 & 626787 & 1.81 & 0.00 & 1.81 & 2 & 2 & 14 & 24 \\  
	&  & C17-100 & 3413 & 1603 & 19 & 902 & 60 & 640945 & 0.00 & 0.00 & 0.00 & 2 & 2 & 2 & 5 \\  
	&  & C18-100 & 3778 & 2518 & 65 & 868 & 58 & 681240 & 1.72 & 0.00 & 1.72 & 2 & 3 & 3 & 9 \\  
	&  & C19-100 & 3722 & 3405 & 633 & 1695 & 51 & 503252 & 1.96 & 0.00 & 1.96 & 4 & 4 & 263 & 359 \\  
	&  & C20-100 & 4377 & 2653 & 1785 & 2948 & 46 & 568797 & 2.17 & 0.00 & 2.17 & 12 & 13 & 1765 & 1815 \\  
	&  & C21-100 & 2597 & 1223 & 1 & 525 & 78 & 328216 & 0.00 & 0.00 & 0.00 & 2 & 2 & 2 & 2 \\ \cline{2-16} 
	& \multirow{22}{*}{6} & C0-100 & 3488 & 1931 & 59 & 782 & 63 & 488119 & 0.00 & 0.00 & 0.00 & 2 & 3 & 3 & 7 \\  
	&  & C1-100 & 1825 & 1642 & 5 & 429 & 75 & 331497 & 0.00 & 0.00 & 0.00 & 2 & 2 & 2 & 2 \\  
	&  & C2-100 & 5540 & 3106 & 1439 & 2614 & 46 & 596658 & 2.17 & 0.00 & 2.17 & 10 & 11 & 1618 & 2309 \\  
	&  & C3-100 & 5383 & 3127 & 415 & 2690 & 46 & 584048 & 2.17 & 0.00 & 2.16 & 76 & 77 & 1791 & 1841 \\  
	&  & C4-100 & 4211 & 2787 & 309 & 2656 & 44 & 552862 & 2.26 & 0.00 & 2.26 & 21 & 24 & 178 & 773 \\  
	&  & C5-100 & 1959 & 1452 & 43 & 583 & 70 & 366565 & 0.00 & 0.00 & 0.00 & 2 & 2 & 2 & 4 \\  
	&  & C6-100 & 4708 & 3393 & 27 & 1821 & 47 & 724844 & 0.00 & 0.00 & 0.00 & 11 & 11 & 30 & 85 \\  
	&  & C7-100 & 4120 & 2205 & 29 & 1243 & 53 & 661528 & 0.00 & 0.00 & 0.00 & 4 & 4 & 4 & 15 \\  
	&  & C8-100 & 4952 & 2656 & 3 & 1277 & 51 & 609038 & 0.00 & 0.00 & 0.00 & 3 & 3 & 3 & 4 \\  
	&  & C9-100 & 5664 & 3572 & 35 & 2997 & 40 & 518571 & 2.49 & 2.49 & 2.49 & 10273 & 10275 & 10275 & 43211 \\  
	&  & C10-100 & 2995 & 2611 & 149 & 982 & 58 & 426386 & 1.72 & 0.00 & 1.72 & 2 & 2 & 16 & 17 \\  
	&  & C11-100 & 4606 & 2964 & 45 & 1453 & 51 & 427829 & 0.00 & 0.00 & 0.00 & 4 & 4 & 18 & 40 \\  
	&  & C12-100 & 2863 & 2459 & 1 & 790 & 62 & 413012 & 0.00 & 0.00 & 0.00 & 2 & 2 & 2 & 2 \\  
	&  & C13-100 & 3232 & 2441 & 17 & 760 & 61 & 483087 & 0.00 & 0.00 & 0.00 & 2 & 2 & 2 & 3 \\  
	&  & C14-100 & 4335 & 3403 & 532 & 3161 & 40 & 693512 & 2.50 & 2.50 & 2.50 & 144 & 148 & 21052 & 43203 \\  
	&  & C15-100 & 3711 & 2642 & 3810 & 3876 & 48 & 788340 & 4.16 & 4.16 & 4.16 & 142 & 142 & 15775 & 43217 \\  
	&  & C16-100 & 3278 & 2302 & 69 & 1151 & 55 & 626787 & 1.81 & 0.00 & 1.81 & 3 & 3 & 3 & 41 \\  
	&  & C17-100 & 3411 & 1603 & 9 & 881 & 60 & 640945 & 0.00 & 0.00 & 0.00 & 2 & 2 & 2 & 3 \\  
	&  & C18-100 & 3778 & 2518 & 59 & 867 & 58 & 681240 & 1.72 & 0.00 & 1.72 & 2 & 2 & 2 & 8 \\  
	&  & C19-100 & 3722 & 3405 & 1364 & 1761 & 51 & 503252 & 1.96 & 0.00 & 1.96 & 5 & 5 & 696 & 977 \\  
	&  & C20-100 & 4377 & 2653 & 2993 & 3190 & 46 & 568797 & 2.17 & 0.00 & 2.17 & 21 & 22 & 3242 & 3252 \\  
	&  & C21-100 & 2597 & 1223 & 1 & 527 & 78 & 328216 & 0.00 & 0.00 & 0.00 & 2 & 2 & 2 & 2 \\ \cline{2-16} 
	& \multirow{22}{*}{7} & C0-100 & 3488 & 1931 & 65 & 783 & 63 & 488119 & 0.00 & 0.00 & 0.00 & 2 & 2 & 2 & 6 \\  
	&  & C1-100 & 1825 & 1642 & 5 & 429 & 75 & 331497 & 0.00 & 0.00 & 0.00 & 2 & 2 & 2 & 2 \\  
	&  & C2-100 & 5540 & 3106 & 1613 & 2710 & 46 & 596658 & 2.17 & 0.00 & 2.17 & 11 & 12 & 1811 & 2822 \\  
	&  & C3-100 & 5383 & 3127 & 247 & 2572 & 46 & 584048 & 2.17 & 0.00 & 2.16 & 68 & 69 & 1250 & 1297 \\  
	&  & C4-100 & 4209 & 2787 & 87 & 2212 & 44 & 552334 & 2.26 & 0.00 & 2.26 & 13 & 15 & 220 & 240 \\  
	&  & C5-100 & 1959 & 1452 & 35 & 583 & 70 & 366565 & 0.00 & 0.00 & 0.00 & 2 & 2 & 2 & 4 \\  
	&  & C6-100 & 4708 & 3393 & 27 & 1807 & 47 & 724844 & 0.00 & 0.00 & 0.00 & 11 & 11 & 29 & 82 \\  
	&  & C7-100 & 4120 & 2205 & 27 & 1209 & 53 & 661528 & 0.00 & 0.00 & 0.00 & 3 & 4 & 4 & 13 \\  
	&  & C8-100 & 4952 & 2656 & 3 & 1270 & 51 & 609038 & 0.00 & 0.00 & 0.00 & 3 & 3 & 3 & 4 \\  
	&  & C9-100 & 5664 & 3572 & 1 & 2482 & 40 & 520597 & 2.50 & 2.50 & 2.50 & 43846 & 43847 & 43847 & 43847 \\  
	&  & C10-100 & 2995 & 2611 & 155 & 981 & 58 & 426386 & 1.72 & 0.00 & 1.72 & 2 & 2 & 16 & 17 \\  
	&  & C11-100 & 4606 & 2964 & 57 & 1456 & 51 & 427829 & 0.00 & 0.00 & 0.00 & 5 & 5 & 29 & 49 \\  
	&  & C12-100 & 2863 & 2459 & 1 & 793 & 62 & 413012 & 0.00 & 0.00 & 0.00 & 2 & 2 & 2 & 2 \\  
	&  & C13-100 & 3232 & 2441 & 17 & 757 & 61 & 483087 & 0.00 & 0.00 & 0.00 & 2 & 2 & 2 & 3 \\  
	&  & C14-100 & 4335 & 3403 & 369 & 2927 & 40 & 692956 & 2.49 & 2.49 & 2.49 & 240 & 242 & 242 & 43244 \\  
	&  & C15-100 & 3711 & 2642 & 993 & 2445 & 48 & 789091 & 4.16 & 4.16 & 4.16 & 942 & 942 & 942 & 43215 \\  
	&  & C16-100 & 3278 & 2302 & 63 & 1131 & 55 & 626787 & 1.81 & 0.00 & 1.81 & 4 & 4 & 28 & 53 \\  
	&  & C17-100 & 3413 & 1603 & 9 & 875 & 60 & 640945 & 0.00 & 0.00 & 0.00 & 2 & 2 & 2 & 3 \\  
	&  & C18-100 & 3778 & 2518 & 81 & 872 & 58 & 681240 & 1.72 & 0.00 & 1.72 & 2 & 2 & 2 & 10 \\  
	&  & C19-100 & 3722 & 3405 & 639 & 1714 & 51 & 503252 & 1.96 & 0.00 & 1.96 & 5 & 6 & 306 & 407 \\  
	&  & C20-100 & 4377 & 2653 & 1669 & 2890 & 46 & 568797 & 2.17 & 0.00 & 2.17 & 16 & 17 & 1180 & 2009 \\  
	&  & C21-100 & 2597 & 1223 & 1 & 525 & 78 & 328216 & 0.00 & 0.00 & 0.00 & 2 & 2 & 2 & 2 \\ \cline{2-16} 
	& \multirow{22}{*}{8} & C0-100 & 3488 & 1931 & 65 & 784 & 63 & 488119 & 0.00 & 0.00 & 0.00 & 2 & 2 & 2 & 7 \\  
	&  & C1-100 & 1825 & 1642 & 5 & 429 & 75 & 331497 & 0.00 & 0.00 & 0.00 & 1 & 1 & 1 & 2 \\  
	&  & C2-100 & 5540 & 3106 & 1477 & 2619 & 46 & 596658 & 2.17 & 0.00 & 2.17 & 13 & 13 & 1790 & 2605 \\  
	&  & C3-100 & 5383 & 3127 & 71 & 2209 & 46 & 584048 & 2.17 & 0.00 & 2.16 & 90 & 91 & 389 & 518 \\  
	&  & C4-100 & 4210 & 2787 & 71 & 2176 & 44 & 552334 & 2.26 & 0.00 & 2.26 & 13 & 15 & 253 & 271 \\  
	&  & C5-100 & 1964 & 1452 & 53 & 575 & 70 & 366565 & 0.00 & 0.00 & 0.00 & 2 & 2 & 2 & 5 \\  
	&  & C6-100 & 4708 & 3393 & 27 & 1822 & 47 & 724844 & 0.00 & 0.00 & 0.00 & 11 & 12 & 29 & 87 \\  
	&  & C7-100 & 4120 & 2205 & 29 & 1250 & 53 & 661528 & 0.00 & 0.00 & 0.00 & 3 & 4 & 4 & 15 \\  
	&  & C8-100 & 4952 & 2656 & 3 & 1273 & 51 & 609038 & 0.00 & 0.00 & 0.00 & 3 & 3 & 3 & 5 \\  
	&  & C9-100 & 5664 & 3572 & 4 & 2577 & 40 & 517627 & 2.49 & 2.49 & 2.49 & 34816 & 34818 & 34818 & 43200 \\  
	&  & C10-100 & 2995 & 2611 & 159 & 976 & 58 & 426386 & 1.72 & 0.00 & 1.72 & 2 & 2 & 16 & 18 \\  
	&  & C11-100 & 4606 & 2964 & 53 & 1457 & 51 & 427829 & 0.00 & 0.00 & 0.00 & 4 & 5 & 19 & 45 \\  
	&  & C12-100 & 2863 & 2459 & 1 & 789 & 62 & 413012 & 0.00 & 0.00 & 0.00 & 2 & 2 & 2 & 2 \\  
	&  & C13-100 & 3232 & 2441 & 13 & 760 & 61 & 483087 & 0.00 & 0.00 & 0.00 & 2 & 2 & 2 & 3 \\  
	&  & C14-100 & 4060 & 3403 & 460 & 3150 & 40 & 693850 & 2.50 & 2.50 & 2.50 & 196 & 198 & 33298 & 43354 \\  
	&  & C15-100 & 3711 & 2642 & 321 & 2353 & 48 & 789091 & 4.16 & 4.16 & 4.16 & 2597 & 2598 & 2598 & 44157 \\  
	&  & C16-100 & 3278 & 2302 & 29 & 1116 & 55 & 626787 & 1.81 & 0.00 & 1.81 & 4 & 4 & 13 & 32 \\  
	&  & C17-100 & 3413 & 1603 & 9 & 873 & 60 & 640945 & 0.00 & 0.00 & 0.00 & 2 & 2 & 2 & 3 \\  
	&  & C18-100 & 3778 & 2518 & 71 & 871 & 58 & 681240 & 1.72 & 0.00 & 1.72 & 2 & 2 & 2 & 9 \\  
	&  & C19-100 & 3722 & 3405 & 719 & 1709 & 51 & 503252 & 1.96 & 0.00 & 1.96 & 5 & 5 & 394 & 487 \\  
	&  & C20-100 & 4377 & 2653 & 1760 & 2899 & 46 & 568797 & 2.17 & 0.00 & 2.17 & 21 & 22 & 1737 & 2192 \\ \down
	&  & C21-100 & 2597 & 1223 & 1 & 524 & 78 & 328216 & 0.00 & 0.00 & 0.00 & 1 & 2 & 2 & 2 \\ \hline
\end{tabular}	
	}
\end{table}

\begin{table}[]
	\TABLE
	{Results of REA Scalability with Increasing Cluster Size ($K = 4$, $\Delta = 10$ mins, $R = 50\%$).\label{tbl:rea_cluster_size_scaling}}
	{
		\begin{tabular}{crrrrrrr}
			\hline \up\down
			\multirow{2}{*}{\begin{tabular}[c]{@{}c@{}}Cluster\\ size\end{tabular}} & \multicolumn{1}{c}{\multirow{2}{*}{\begin{tabular}[c]{@{}c@{}}Cluster\\ ID\end{tabular}}} & \multicolumn{1}{c}{\multirow{2}{*}{\begin{tabular}[c]{@{}c@{}}Column\\ \#\end{tabular}}} & \multicolumn{1}{c}{\multirow{2}{*}{\begin{tabular}[c]{@{}c@{}}Vehicle\\ \#\end{tabular}}} & \multicolumn{1}{c}{\multirow{2}{*}{\begin{tabular}[c]{@{}c@{}}Total\\ distance\\ (m)\end{tabular}}} & \multicolumn{1}{c}{\multirow{2}{*}{\begin{tabular}[c]{@{}c@{}}Optimality\\ gap (\%)\end{tabular}}} & \multicolumn{2}{c}{Wall time (s)} \\ \cline{7-8}
			& \multicolumn{1}{c}{} & \multicolumn{1}{c}{} & \multicolumn{1}{c}{} & \multicolumn{1}{c}{} & \multicolumn{1}{c}{} & \multicolumn{1}{c}{\begin{tabular}[c]{@{}c@{}}Route\\ enumeration\end{tabular}} & \multicolumn{1}{c}{Total} \\ \hline \up
			\multirow{11}{*}{200} & C0-200 & 2286 & 129 & 702670 & 0.00 & 618 & 618 \\  
			& C1-200 & 10955 & 84 & 1122469 & 0.00 & 653 & 654 \\  
			& C2-200 & 15149 & 85 & 1067503 & 0.00 & 646 & 2361 \\  
			& C3-200 & 3101 & 113 & 824730 & 0.00 & 629 & 630 \\  
			& C4-200 & 2627 & 114 & 684981 & 0.00 & 632 & 878 \\  
			& C5-200 & 22061 & 80 & 937514 & 0.00 & 660 & 744 \\  
			& C6-200 & 11844 & 86 & 1136870 & 0.00 & 673 & 680 \\  
			& C7-200 & 16235 & 76 & 1365475 & 0.00 & 668 & 1165 \\  
			& C8-200 & 27913 & 78 & 923290 & 0.49 & 684 & 43211 \\  
			& C9-200 & 5289 & 99 & 1055393 & 0.00 & 655 & 702 \\  \down
			& C10-200 & 7790 & 87 & 1249784 & 0.00 & 676 & 8731 \\ \hline \up
			\multirow{13}{*}{300} & C0-300 & 33262 & 118 & 1652972 & 0.16 & 3301 & 43206 \\  
			& C1-300 & 32536 & 119 & 1928887 & 0.44 & 2796 & 43203 \\  
			& C2-300 & 7994 & 150 & 1096341 & 0.00 & 3249 & 4367 \\  
			& C3-300 & 36568 & 113 & 1511806 & 0.00 & 3499 & 3650 \\  
			& C4-300 & 15394 & 134 & 1477823 & 0.13 & 3198 & 43205 \\  
			& C5-300 & 22072 & 124 & 1773852 & 0.00 & 3342 & 4325 \\  
			& C6-300 & 33541 & 119 & 1927882 & 0.35 & 3191 & 43205 \\  
			& C7-300 & 6554 & 157 & 1105570 & 0.00 & 3152 & 3158 \\  
			& C8-300 & 53120 & 114 & 1418494 & 0.20 & 3129 & 43204 \\  
			& C9-300 & 30568 & 114 & 1911085 & 0.17 & 2258 & 43203 \\  
			& C10-300 & 6370 & 156 & 1085481 & 0.00 & 2736 & 2776 \\  
			& C11-300 & 34630 & 118 & 1903034 & 0.37 & 3037 & 43204 \\  \down
			& C12-300 & 7137 & 153 & 1067479 & 0.00 & 2339 & 2632 \\ \hline \up
			\multirow{8}{*}{400} & C0-400 & 28968 & 180 & 1617748 & 0.00 & 8897 & 35048 \\  
			& C1-400 & 52194 & 145 & 1901974 & 0.19 & 6972 & 43203 \\  
			& C2-400 & 28025 & 173 & 1627074 & 0.00 & 10060 & 18906 \\  
			& C3-400 & 41012 & 152 & 2114911 & 0.20 & 7842 & 43203 \\  
			& C4-400 & 26314 & 182 & 1687079 & 0.00 & 10447 & 10531 \\  
			& C5-400 & 53948 & 151 & 1988361 & 0.00 & 7043 & 12685 \\  
			& C6-400 & 26109 & 173 & 1640438 & 0.00 & 10716 & 10986 \\  \down
			& C7-400 & 68190 & 145 & 1887086 & 0.32 & 7117 & 43202 \\ \hline
		\end{tabular}%
	}
{}
\end{table}

\clearpage

\begin{table}[]
	\TABLE
	{Results of BPA Scalability with Increasing Cluster Size ($K = 4$, $\Delta = 10$ mins, $R = 50\%$).\label{tbl:bpa_cluster_size_scaling}}
	{\resizebox{\textwidth}{!}{
			\begin{tabular}{crrrrrrrrrrrrrr}
				\hline \up\down
				\multirow{2}{*}{\begin{tabular}[c]{@{}c@{}}Cluster\\ size\end{tabular}} & \multicolumn{1}{c}{\multirow{2}{*}{\begin{tabular}[c]{@{}c@{}}Cluster\\ ID\end{tabular}}} & \multicolumn{1}{c}{\multirow{2}{*}{\begin{tabular}[c]{@{}c@{}}Inbound\\ edge \#\end{tabular}}} & \multicolumn{1}{c}{\multirow{2}{*}{\begin{tabular}[c]{@{}c@{}}Outbound\\ edge \#\end{tabular}}} & \multicolumn{1}{c}{\multirow{2}{*}{\begin{tabular}[c]{@{}c@{}}Tree\\ node\\ \#\end{tabular}}} & \multicolumn{1}{c}{\multirow{2}{*}{\begin{tabular}[c]{@{}c@{}}Column\\ \#\end{tabular}}} & \multicolumn{1}{c}{\multirow{2}{*}{\begin{tabular}[c]{@{}c@{}}Vehicle\\ \#\end{tabular}}} & \multicolumn{1}{c}{\multirow{2}{*}{\begin{tabular}[c]{@{}c@{}}Total\\ distance\\ (m)\end{tabular}}} & \multicolumn{2}{c}{Optimality gap (\%)} & \multicolumn{1}{c}{\multirow{2}{*}{\begin{tabular}[c]{@{}c@{}}Integrality\\ gap (\%)\end{tabular}}} & \multicolumn{4}{c}{Wall time (s)} \\ \cline{9-10} \cline{12-15} 
				& \multicolumn{1}{c}{} & \multicolumn{1}{c}{} & \multicolumn{1}{c}{} & \multicolumn{1}{c}{} & \multicolumn{1}{c}{} & \multicolumn{1}{c}{} & \multicolumn{1}{c}{} & \multicolumn{1}{c}{\begin{tabular}[c]{@{}c@{}}Root\\node soln.\end{tabular}} & \multicolumn{1}{c}{\begin{tabular}[c]{@{}c@{}}Best\\ feasible\\ soln.\end{tabular}} & \multicolumn{1}{c}{} & \multicolumn{1}{c}{\begin{tabular}[c]{@{}c@{}}RMP\\ conv.\end{tabular}} & \multicolumn{1}{c}{\begin{tabular}[c]{@{}c@{}}Root\\ node\\ soln.\end{tabular}} & \multicolumn{1}{c}{\begin{tabular}[c]{@{}c@{}}Best\\ feasible\\ soln.\end{tabular}} & \multicolumn{1}{c}{Total} \\ \hline \up
				\multirow{11}{*}{200} & C0-200 & 10281 & 5900 & 113 & 1747 & 129 & 702670 & 1.55 & 0.00 & 1.55 & 13 & 13 & 35 & 40 \\  
				& C1-200 & 17096 & 10087 & 5870 & 6121 & 84 & 1122544 & 3.50 & 2.37 & 2.37 & 59 & 96 & 6951 & 43201 \\  
				& C2-200 & 13922 & 9341 & 7479 & 7506 & 85 & 1067640 & 2.35 & 2.34 & 2.34 & 55 & 74 & 6263 & 43203 \\  
				& C3-200 & 10551 & 7643 & 99 & 2272 & 113 & 824730 & 0.00 & 0.00 & 0.00 & 19 & 20 & 20 & 78 \\  
				& C4-200 & 9806 & 8163 & 5901 & 2616 & 114 & 684981 & 1.75 & 0.00 & 1.75 & 17 & 31 & 2896 & 4149 \\  
				& C5-200 & 18385 & 10756 & 3695 & 9115 & 80 & 937514 & 1.26 & 1.25 & 1.25 & 111 & 121 & 14244 & 43203 \\  
				& C6-200 & 16047 & 9533 & 13563 & 9572 & 86 & 1136870 & 2.29 & 1.16 & 1.16 & 45 & 202 & 14218 & 43201 \\  
				& C7-200 & 15672 & 11696 & 2279 & 8993 & 76 & 1367209 & 3.88 & 2.62 & 2.62 & 140 & 6148 & 23012 & 43213 \\  
				& C8-200 & 16483 & 12874 & 1 & 6027 & 79 & 928401 & 2.53 & 2.53 & 2.53 & 319 & 43209 & 43209 & 43209 \\  
				& C9-200 & 14587 & 9841 & 15783 & 4425 & 99 & 1055393 & 1.01 & 0.00 & 1.00 & 25 & 37 & 6671 & 20722 \\ \down 
				& C10-200 & 15459 & 11600 & 4244 & 6328 & 87 & 1249965 & 2.29 & 2.29 & 2.29 & 45 & 150 & 9471 & 43202 \\ \hline \up
				\multirow{13}{*}{300} & C0-300 & 38471 & 31409 & 242 & 11624 & 118 & 1665487 & 0.85 & 0.85 & 0.85 & 900 & 1018 & 1018 & 43265 \\  
				& C1-300 & 36931 & 22949 & 496 & 12533 & 119 & 1938007 & 2.52 & 2.52 & 2.52 & 944 & 4151 & 4151 & 43229 \\  
				& C2-300 & 26167 & 20901 & 11208 & 6629 & 150 & 1096341 & 1.33 & 1.33 & 1.33 & 81 & 124 & 7967 & 43200 \\  
				& C3-300 & 35614 & 23695 & 762 & 13192 & 114 & 1510444 & 1.75 & 1.75 & 1.75 & 426 & 12864 & 12864 & 43212 \\  
				& C4-300 & 30651 & 26297 & 1091 & 8757 & 134 & 1478291 & 1.49 & 1.49 & 1.49 & 278 & 519 & 39914 & 43243 \\  
				& C5-300 & 32453 & 25677 & 400 & 9653 & 125 & 1772096 & 0.80 & 0.80 & 0.80 & 622 & 14460 & 14460 & 43259 \\  
				& C6-300 & 35954 & 20702 & 483 & 12056 & 119 & 1934456 & 2.52 & 2.52 & 2.52 & 734 & 1914 & 1914 & 43218 \\  
				& C7-300 & 24094 & 19880 & 3185 & 5592 & 157 & 1105570 & 1.26 & 0.00 & 0.64 & 70 & 1027 & 3841 & 10016 \\  
				& C8-300 & 38059 & 22799 & 255 & 12441 & 115 & 1423615 & 2.60 & 2.60 & 2.60 & 480 & 29027 & 29027 & 43201 \\  
				& C9-300 & 33512 & 23617 & 779 & 13090 & 114 & 1934405 & 2.64 & 2.64 & 2.64 & 515 & 666 & 666 & 43252 \\  
				& C10-300 & 22339 & 16909 & 23037 & 5367 & 156 & 1085482 & 0.64 & 0.64 & 0.64 & 63 & 70 & 7449 & 43202 \\  
				& C11-300 & 34544 & 21130 & 698 & 12376 & 118 & 1910713 & 3.38 & 3.38 & 3.38 & 568 & 988 & 988 & 43227 \\ \down 
				& C12-300 & 24520 & 17238 & 12001 & 6091 & 153 & 1067479 & 0.65 & 0.65 & 0.65 & 88 & 124 & 4043 & 43200 \\ \hline \up
				\multirow{8}{*}{400} & C0-400 & 49218 & 39758 & 1 & 10357 & 181 & 1613436 & 1.65 & 1.65 & 1.65 & 528 & 43205 & 43205 & 43205 \\  
				& C1-400 & 53747 & 35507 & 1 & 13136 & 147 & 1926300 & 2.72 & 2.72 & 2.72 & 751 & 43207 & 43207 & 43207 \\  
				& C2-400 & 54300 & 35489 & 675 & 12751 & 173 & 1631283 & 1.16 & 1.16 & 1.16 & 406 & 540 & 540 & 43223 \\  
				& C3-400 & 55542 & 42004 & 1 & 12834 & 153 & 2124836 & 1.96 & 1.96 & 1.96 & 1596 & 43204 & 43204 & 43204 \\  
				& C4-400 & 49195 & 31033 & 3523 & 11881 & 182 & 1688282 & 1.10 & 1.10 & 1.10 & 222 & 239 & 25641 & 43201 \\  
				& C5-400 & 57633 & 34496 & 1 & 13557 & 152 & 1983300 & 1.97 & 1.97 & 1.97 & 684 & 43204 & 43204 & 43204 \\  
				& C6-400 & 46155 & 30320 & 3708 & 12494 & 173 & 1640472 & 1.16 & 1.15 & 1.15 & 240 & 289 & 11566 & 43210 \\ \down 
				& C7-400 & 61393 & 38346 & 1 & 14833 & 146 & 1886534 & 2.73 & 2.73 & 2.73 & 960 & 43206 & 43206 & 43206 \\ \hline
			\end{tabular}
	}}
	{}
\end{table}

\clearpage

\begin{table}[]
	\TABLE
	{Results of Root-Node Heuristics with $t_{\text{RMP}} = 8$ mins and $t_{\text{MIP}} = 2$ mins ($K = 4$, $\Delta = 10$ mins, $R = 50\%$).\label{tbl:root_node_heuristic}}
	{\resizebox{\textwidth}{!}{
		\begin{tabular}{cr|rrrrr|rrrrrr}
			\hline \up\down
			\multirow{3}{*}{\begin{tabular}[c]{@{}c@{}}Cluster\\ size\end{tabular}} & \multicolumn{1}{c|}{\multirow{3}{*}{\begin{tabular}[c]{@{}c@{}}Cluster\\ ID\end{tabular}}} & \multicolumn{5}{c|}{Enforce forbidden paths} & \multicolumn{6}{c}{Relax forbidden paths} \\ \cline{3-13} \up\down
			& \multicolumn{1}{c|}{} & \multicolumn{1}{c}{\multirow{2}{*}{\begin{tabular}[c]{@{}c@{}}Column\\ \#\end{tabular}}} & \multicolumn{1}{c}{\multirow{2}{*}{\begin{tabular}[c]{@{}c@{}}Vehicle\\ \#\end{tabular}}} & \multicolumn{1}{c}{\multirow{2}{*}{\begin{tabular}[c]{@{}c@{}}Optimality\\ gap (\%)\end{tabular}}} & \multicolumn{2}{c|}{Wall time (s)} & \multicolumn{1}{c}{\multirow{2}{*}{\begin{tabular}[c]{@{}c@{}}Column\\ \#\end{tabular}}} & \multicolumn{1}{c}{\multirow{2}{*}{\begin{tabular}[c]{@{}c@{}}Infeasible\\ column \#\end{tabular}}} & \multicolumn{1}{c}{\multirow{2}{*}{\begin{tabular}[c]{@{}c@{}}Vehicle\\ \#\end{tabular}}} & \multicolumn{1}{c}{\multirow{2}{*}{\begin{tabular}[c]{@{}c@{}}Optimality\\ gap (\%)\end{tabular}}} & \multicolumn{2}{c}{Wall time (s)} \\ \cline{6-7} \cline{12-13} \up\down
			& \multicolumn{1}{c|}{} & \multicolumn{1}{c}{} & \multicolumn{1}{c}{} & \multicolumn{1}{c}{} & \multicolumn{1}{c}{RMP} & \multicolumn{1}{c|}{MIP} & \multicolumn{1}{c}{} & \multicolumn{1}{c}{} & \multicolumn{1}{c}{} & \multicolumn{1}{c}{} & \multicolumn{1}{c}{RMP} & \multicolumn{1}{c}{MIP} \\ \hline \up
			\multirow{22}{*}{100} & C0-100 & 720 & 63 & 0.00 & 2 & 0 & 720 & 0 & 63 & 0.00 & 2 & 0 \\ 
			& C1-100 & 422 & 75 & 0.00 & 2 & 0 & 423 & 0 & 75 & 0.00 & 2 & 0 \\ 
			& C2-100 & 1828 & 46 & 2.17 & 5 & 0 & 1792 & 6 & 46 & 4.35 & 4 & 0 \\ 
			& C3-100 & 1735 & 46 & 2.17 & 4 & 0 & 1729 & 6 & 46 & 2.17 & 4 & 0 \\ 
			& C4-100 & 1653 & 44 & 2.27 & 6 & 5 & 1624 & 4 & 44 & 2.27 & 5 & 1 \\ 
			& C5-100 & 537 & 70 & 0.00 & 1 & 0 & 541 & 0 & 70 & 0.00 & 2 & 0 \\ 
			& C6-100 & 1637 & 47 & 0.00 & 5 & 0 & 1628 & 3 & 47 & 0.00 & 4 & 0 \\ 
			& C7-100 & 1130 & 53 & 0.00 & 2 & 0 & 1123 & 5 & 53 & 0.00 & 2 & 0 \\ 
			& C8-100 & 1198 & 52 & 1.92 & 2 & 0 & 1184 & 3 & 51 & 0.00 & 2 & 0 \\ 
			& C9-100 & 2360 & 41 & 4.88 & 29 & 10 & 2467 & 10 & 41 & 4.88 & 37 & 2 \\ 
			& C10-100 & 891 & 58 & 1.72 & 2 & 0 & 895 & 1 & 58 & 1.72 & 2 & 0 \\ 
			& C11-100 & 1272 & 51 & 0.00 & 3 & 0 & 1252 & 0 & 51 & 0.00 & 2 & 0 \\ 
			& C12-100 & 767 & 62 & 0.00 & 2 & 0 & 776 & 2 & 62 & 0.00 & 2 & 0 \\ 
			& C13-100 & 708 & 61 & 0.00 & 2 & 0 & 708 & 0 & 61 & 0.00 & 1 & 0 \\ 
			& C14-100 & 2003 & 41 & 4.88 & 11 & 2 & 1994 & 12 & 41 & 4.88 & 8 & 2 \\ 
			& C15-100 & 1548 & 48 & 4.17 & 3 & 0 & 1580 & 9 & 48 & 6.25 & 4 & 0 \\ 
			& C16-100 & 992 & 55 & 1.82 & 2 & 0 & 988 & 2 & 55 & 1.82 & 2 & 0 \\ 
			& C17-100 & 817 & 60 & 0.00 & 2 & 0 & 810 & 3 & 60 & 0.00 & 2 & 0 \\ 
			& C18-100 & 808 & 58 & 1.72 & 2 & 0 & 808 & 0 & 58 & 1.72 & 1 & 0 \\ 
			& C19-100 & 1337 & 51 & 1.96 & 3 & 0 & 1334 & 9 & 52 & 3.85 & 3 & 0 \\ 
			& C20-100 & 1672 & 46 & 2.17 & 4 & 1 & 1643 & 2 & 46 & 2.17 & 4 & 1 \\ \down
			& C21-100 & 500 & 78 & 0.00 & 2 & 0 & 503 & 0 & 78 & 0.00 & 1 & 0 \\ \hline \up
			\multirow{11}{*}{200} & C0-200 & 1537 & 129 & 1.55 & 12 & 0 & 1539 & 2 & 129 & 1.55 & 12 & 0 \\ 
			& C1-200 & 4216 & 85 & 3.53 & 35 & 1 & 4146 & 19 & 84 & 2.38 & 35 & 2 \\ 
			& C2-200 & 4019 & 85 & 2.35 & 34 & 2 & 4018 & 16 & 85 & 3.53 & 33 & 3 \\ 
			& C3-200 & 2029 & 113 & 0.00 & 13 & 0 & 2022 & 4 & 113 & 0.00 & 13 & 0 \\ 
			& C4-200 & 2017 & 114 & 1.75 & 14 & 0 & 2039 & 4 & 114 & 1.75 & 13 & 0 \\ 
			& C5-200 & 5161 & 80 & 1.25 & 82 & 9 & 5010 & 29 & 81 & 2.47 & 48 & 35 \\ 
			& C6-200 & 4294 & 87 & 2.30 & 28 & 1 & 4315 & 18 & 87 & 2.30 & 25 & 3 \\ 
			& C7-200 & 5516 & 76 & 2.63 & 81 & 67 & 5563 & 38 & 77 & 3.90 & 63 & 12 \\ 
			& C8-200 & 5925 & 79 & 2.53 & 245 & 120 & 5986 & 29 & 79 & 2.53 & 209 & 120 \\ 
			& C9-200 & 3142 & 99 & 1.01 & 16 & 2 & 3165 & 8 & 99 & 1.01 & 17 & 1 \\ \down
			& C10-200 & 4118 & 87 & 2.30 & 38 & 5 & 4150 & 10 & 87 & 2.30 & 30 & 3 \\ \hline \up
			\multirow{13}{*}{300} & C0-300 & 9312 & 119 & 1.68 & 449 & 120 & 9366 & 58 & 119 & 1.68 & 235 & 120 \\ 
			& C1-300 & 9928 & 119 & 2.52 & 376 & 120 & 9795 & 44 & 120 & 4.17 & 171 & 114 \\ 
			& C2-300 & 4879 & 150 & 1.33 & 63 & 5 & 4894 & 3 & 150 & 1.33 & 62 & 6 \\ 
			& C3-300 & 9465 & 114 & 1.75 & 286 & 62 & 9530 & 43 & 116 & 3.45 & 264 & 120 \\ 
			& C4-300 & 7002 & 134 & 1.49 & 192 & 8 & 7021 & 24 & 134 & 1.49 & 113 & 14 \\ 
			& C5-300 & 7920 & 125 & 0.80 & 296 & 120 & 7832 & 30 & 125 & 1.60 & 142 & 19 \\ 
			& C6-300 & 9495 & 119 & 2.52 & 433 & 28 & 9531 & 70 & 119 & 3.36 & 157 & 14 \\ 
			& C7-300 & 4311 & 158 & 1.27 & 63 & 4 & 4298 & 3 & 158 & 1.27 & 55 & 1 \\ 
			& C8-300 & 9726 & 116 & 3.45 & 243 & 120 & 9756 & 31 & 115 & 2.61 & 219 & 120 \\ 
			& C9-300 & 9796 & 115 & 3.48 & 401 & 120 & 9835 & 55 & 115 & 3.48 & 304 & 120 \\ 
			& C10-300 & 3884 & 156 & 0.64 & 50 & 2 & 3915 & 8 & 156 & 1.28 & 49 & 1 \\ 
			& C11-300 & 9354 & 119 & 4.20 & 303 & 120 & 9685 & 51 & 118 & 4.24 & 196 & 26 \\ \down
			& C12-300 & 4207 & 153 & 0.65 & 59 & 1 & 4285 & 6 & 153 & 0.65 & 53 & 1 \\ \hline \up
			\multirow{8}{*}{400} & C0-400 & 10289 & 181 & 1.66 & 356 & 120 & 10355 & 27 & 182 & 2.20 & 318 & 120 \\ 
			& C1-400 & 12853 & 152 & 5.92 & 417 & 120 & 12965 & 42 & 151 & 5.30 & 419 & 120 \\ 
			& C2-400 & 9989 & 173 & 1.16 & 185 & 42 & 10112 & 19 & 173 & 1.16 & 203 & 21 \\ 
			& C3-400 & 12615 & 155 & 11.61 & 486 & 120 & 12714 & 44 & 157 & 4.46 & 472 & 120 \\ 
			& C4-400 & 8872 & 182 & 1.10 & 153 & 18 & 8948 & 22 & 183 & 2.19 & 146 & 26 \\ 
			& C5-400 & 13162 & 154 & 3.25 & 415 & 120 & 13248 & 45 & 153 & 2.61 & 443 & 120 \\ 
			& C6-400 & 8971 & 173 & 1.16 & 155 & 70 & 9130 & 28 & 173 & 1.73 & 159 & 25 \\ \down
			& C7-400 & 14322 & 150 & 6.67 & 504 & 120 & 14412 & 64 & 152 & 6.58 & 461 & 120 \\ \hline
		\end{tabular}	
	}}
	{}
\end{table}

\clearpage

\end{APPENDIX}
%
%


\bibliographystyle{ormsv080} 
\bibliography{references} 


\end{document}